\documentclass[a4paper]{article}
\usepackage[english]{babel}
\usepackage[utf8]{inputenc}
\usepackage{amsmath,amsthm}
\usepackage{amsfonts,amssymb}

\usepackage{mathtools}

\usepackage{tikz,tikz-cd} 

\usepackage{xifthen} 

\newcommand{\T}{\mathbb{T}}
\newcommand{\J}{\mathbb{J}}

\newcommand{\Hb}{\mathbb{H}}

\newcommand{\Acal}{\mathcal{A}}

\newcommand{\Ycal}{\mathcal{Y}} 
 
\newcommand{\Ical}{\mathcal{I}}

\newcommand{\Dcal}{\mathcal{D}} 
\newcommand{\Fcal}{\mathcal{F}}

\newcommand{\Jcal}{\mathcal{J}}

\newcommand{\Wcal}{\mathcal{W}} 
\newcommand{\Xcal}{\mathcal{X}} 
\newcommand{\Ccal}{\mathcal{C}} 
\newcommand{\Vcal}{\mathcal{V}} 
\newcommand{\Bcal}{\mathcal{B}}

\newcommand{\setop}{\textsc{Set}^{op}}
\newcommand{\cof}{ \textsc{Cof}}
\newcommand{\fgpd}{\textsc{Gpd}^f}
\newcommand{\fglob}{\textsc{Glob}^f}

\usepackage{titlesec}
\titleformat{\subsubsection}[runin]{\normalfont}{\thesubsubsection}{0pt}{}[.]

\renewcommand{\thesubsubsection}{\arabic{section}.\arabic{subsection}.\arabic{subsubsection}}
\csname @addtoreset\endcsname{subsubsection}{section}

\newcommand{\block}[1]
{

\par \subsubsection{} #1

\bigskip}

\newcommand{\blockn}[1]{\par #1 \bigskip}
\newcommand{\blockp}[1]{\par #1}

\newcommand{\Th}[1]
	{
	\bigskip	
	\textbf{Theorem : }{\itshape #1}
		
	\bigskip
	}

\newcommand{\Prop}[1]
	{

	\bigskip
	
	\textbf{Proposition : }{\itshape #1}
		
	\bigskip
	
	}

\newcommand{\Conjecture}[1]
	{

	\bigskip
	
	\textbf{Conjecture : }{\itshape #1}
		
	\bigskip
	
	}

\newcommand{\Cor}[1]
	{

	\bigskip
	
	\textbf{Corollary : }{\itshape #1}	
		
	\bigskip

	}

\newcommand{\Lem}[1]
	{

	\bigskip
	
	\textbf{Lemma : }{\itshape #1}
		
	\bigskip
	
	}

\newcommand{\Def}[1]
	{
	
	\bigskip
	
	\textbf{Definition : }{\itshape #1}
	
	\bigskip
	
	}

\newcommand{\Dem}[1]{
	
	\smallskip
	
	\textbf{Proof : } \par
	 {#1} $\square$
	 
	 \bigskip
}

\begin{document}

\pagestyle{plain}
\title{Algebraic models of homotopy types and the homotopy hypothesis}

\author{Simon Henry}

\maketitle

\begin{abstract}

We introduce and study a notion of cylinder coherator similar to the notion of Grothendieck coherator which define more flexible notion of weak $\infty$-groupoids. We show that each such cylinder coherator produces a combinatorial semi-model category of weak $\infty$-groupoids, whose objects are all fibrant and which is in a precise sense ``freely generated by an object''. We show that all those semi model categories are Quillen equivalent together and Quillen to the model category of spaces.
A general procedure is given to produce such coherator, and several explicit examples are presented: one which is simplicial in nature and allows the comparison to the model category for spaces. A second example can be describe as the category of globular sets endowed with ``all the operations that can be defined within a weak type theory''. This second notion seem to provide a definition of weak infinity groupoids which can be defined internally within type theory and which is classically equivalent to homotopy types. Finally, the category of Grothendieck $\infty$-groupoids for a fixed Grothendieck coherator would be an example of this formalism under a seemingly simple conjecture whose validity is shown to imply Grothendieck homotopy hypothesis. This conjecture seem to sum up what needs to be proved at a technical level to ensure that the theory of Grothendieck weak $\infty$-groupoid is well behaved.
\end{abstract}

\renewcommand{\thefootnote}{\fnsymbol{footnote}} 
\footnotetext{\emph{Keywords.} $\infty$-goupoids, cylinder categories, homotopy hypothesis}
\footnotetext{\emph{2010 Mathematics Subject Classification.} 55U35,18G30 }
\footnotetext{\emph{email:} simon.henry@college-de-france.fr}
\renewcommand{\thefootnote}{\arabic{footnote}}


\tableofcontents

\section{introduction}

\subsection{Motivations}

\blockn{In the unpublished manuscript \cite{grothendieck1983pursuing}, A.Grothendieck introduced a very natural (although quite ingenious) notion of ``weak\footnote{Weak means that all the axioms involved in the definition are never stated as equalities between two operations but always as existence of isomorphisms satisfying new axioms, themselve stated as existence of isomorphisms.} $\infty$-groupoid'', a presentation of this notion can be found in \cite{maltsiniotis2010grothendieck}. Grothendieck showed in particular that to any topological space one can associate such a ``fundamental'' weak $\infty$-groupoids, encoding paths, homotopy between path, homotopy between homotopy and so on, and he stated a precise conjecture asserting that this fundamental $\infty$-groupoid construction induces an equivalence between the homotopy category of spaces and the homotopy category of his weak $\infty$-groupoids. This conjecture is known as Grothendieck Homotopy Hypothesis (conjecture $2.8$ of \cite{maltsiniotis2010grothendieck}). It is well known that strict $\infty$-groupoids, while considerably easier to define, are not sufficient for this kind of results, for example it is proved in \cite{simpson1998homotopy} that the fundamental $3$-groupoids of the $2$-sphere cannot be represented in a reasonable way by a strict $3$-groupoids, and other obstruction results in this spirit where known long before.}

\blockn{Grothendieck's definition of weak $\infty$-groupoids has also been extended by G.Mal\-tsi\-nio\-tis (in \cite{maltsiniotis2010grothendieck}) to a definition of weak $\infty$-categories (often called Grothendieck-Maltsiniotis $\infty$-categories). A related definition of weak $\infty$-categories has been proposed independently by M.Batanin in \cite{batanin1998monoidal} and improved by T.Leinster in \cite{leinster2004higher}. These are called Batanin or Batanin-Leinster $\infty$-categories. The precise relationship between Grothendieck-Maltsiniotis and Batanin-Leinster $\infty$-category is explained in D.Ara p.h.d. thesis, \cite{ara2010infini}. All these approaches are ``globular'' in the sense that an $\infty$-category is represented by a globular set with operations, i.e. one has a set of objects, a set of morphisms between any two objects, a set of $2$-morphisms between any two parallel morphisms, $3$- morphisms between any pair of parallel $2$-morphisms, etc. with all the compositions operations and coherence isomorphisms that one can expect in a $\infty$-categories.
}

\blockn{Nowadays, it seems that, while the fact that weak $\infty$-groupoids and spaces up to weak homotopy equivalence are the same is considered as a very important and deep fact, Grothendieck homotopy hypothesis seem to be considered more as a marginal problem. Indeed it has generally been accepted that instead of trying to define what is a $\infty$-groupoid one should take the homotopy hypothesis as their definition.
One can for example replace spaces by a more algebraic category which is known to be equivalent and define an $\infty$-groupoid to be a simplicial Kan complex.

This has been taken as the starting point of higher category theory. From this definition of $\infty$-groupoid as the Kan complex, one have then been able to define $(\infty,1)$-categories , i.e. $\infty$-categories where all the $2$-arrows and above are invertible, as simplicially enriched categories, complete Segal spaces or as quasi-categories, and those have been studied extensively. In fact large portion of ordinary category theory can be extended to those $(\infty,1)$-categories (see for examples \cite{lurie2009higher} and \cite{luriehigher}). 
}

\blockn{Moreover, the fact that the homotopy hypothesis has not been proved is not the only problem with the ``globular'' approaches to higher categories mentioned above:

Their definition depends on certain choices (a coherator for Grothendieck-Maltsi\-niotis $\infty$-categories, and a weakly initial contractible globular operad for Batanin-Leinster $\infty$-categories), different choices produces different categories of $\infty$-groupoids (and strict functors) and it is unknown if these choices are unimportant ``up to homotopy''. Also defining what are weak functors and weak natural transformations between $\infty$-categories is problematic, there has been approaches proposed (see for example \cite[8.8 and 8.9]{batanin1998monoidal} ) but it is not always clear how to compose two such functors and natural transformation and in any case we do not know how to show that one can form an $\infty$-category of $\infty$-groupoids or of $\infty$-categories. We do not know either how to form an $\infty$-groupoids of functors between two $\infty$-groupoids, in fact we do not even know how to form the $\infty$-groupoid of arrows of a given $\infty$-groupoid\footnote{In fact, we think that being able to do this would allow to prove our conjecture \ref{conjecture}, which would in turn solve all the problems mentioned above.} or slices of an $\infty$-category. All these problems prevent us to develop a nice theory of globular higher categories as it has been done for the simplicial approaches (for example in \cite{lurie2009higher} and \cite{luriehigher}).}

\blockn{This might sound very discouraging and a good reason to completely forget about globular approaches to higher category theory, but the recent development of homotopy type theory (see \cite{voevodsky2013homotopy}) has given us very good reasons to be interested in these globular approaches again.

Indeed, one thing that seem to be missing in present higher category theory is a good notion of internal logic of a $(\infty,1)$-category, and it is exactly what intentional type theory is supposed to be: intentional type theory with just identity type and $\Sigma$-type should be the internal logic of a $(\infty,1)$-category with finite projective limits. Strongest type theories ( with functions type, functional extensionality , $\Pi$-type, higher inductive type, univalence axiom etc... ) should corresponds to more structured higher categories (like locally cartesian closed or co-complete) and full homotopy type theory as presented in \cite{voevodsky2013homotopy} should corresponds to the internal logic of an $(\infty,1)$-topos. 

This already gives us a first reason to be interested in Grothendieck $\infty$-groupoids and more generally in globular Higher categories: type theory is deeply globular, indeed, any naive notion of model of type theory produces globular sets and even globular higher groupoids (see \cite{van2011types}, \cite{lumsdaine2009weak} and \cite{bourke2016note}). One could says that type theory should be in fact the internal logic of a certain type of \emph{globular} $(\infty,1)$-categories.

A second reason is that it is of prime importance for ordinary categorical logic to be able to discuss category theory itself ``internally''. One should hence expect that the internal logic of a Higher categories should be able to discuss higher category theory as well.
It appears that intentional type theory has been unable at the present time to discuss simplicial or semi-simplicial sets, in the sense that we are not able to write down a definition of ``simplicial types'' making sense in type theory, hence there is no chance (at least at the present time) that we can use quasi-categories or simplicially enriched categories within type theory. Being able to discuss higher categories under one form or another is one of the most important open problem in homotopy type theory.

Surprisingly, it appears that it is not too difficult to define a globular type (it is generally done using co-induction) and it has been shown (see \cite{brunerie2013syntactic},\cite{altenkirch2012syntactical}) that it is possible to define a notion of globular weak $\infty$-groupoid meaningful in type theory, and very close, if not equivalent, to Grothendieck's original definition. Of course as we are not able to manipulate globular higher groupoids correctly in classical mathematics we do not know how to do it in type theory either so this does not told us how to discuss higher categories correctly within type theory.
}

\blockn{This work started as an attempt to make the theory of globular higher categories work as well as the other approaches  (for example simplicial) to higher categories at least in classical mathematics.

We failed to completely solve any of the problems mentioned above for either the Gorthendieck-Maltsiniotis or the Batanin-Leinster approaches, but one has obtain some results in that direction, as well as ways to avoid those difficult problems. completely

For examples, one of the consequences of the following work is a new globular approach to weak $\infty$-groupoids based on type theory for which we can prove in classical mathematics the analogue of the homotopy hypothesis and which seems to be sufficient for the two motivations mentioned above: it should be definable within type theory and any type in intentional type theory caries the structure of such an $\infty$-groupoid. }

\blockn{Another problem is that the definitions of higher categorical structures that are classically used (like quasi-categories) seems to all be deeply ``non-constructive''. We think that the models we are producing here for $\infty$-groupoids will be more suited for studying higher categorical structures within constructive framework.}

\blockn{Although the long terms goal of this works is to obtain model of Higher categorical structures that will work within constructive mathematics and intentional type theory, the present paper focus on the non-constructive point of view: our main goal is to prove that our new definitions of higher groupoids are equivalent to more classical (and generally non-constructive) notion of higher groupoids (space or Kan complex).

Also we will entirely focus on $\infty$-groupoids, but we hope that the methods developed in this paper might allows us to discuss more general higher structures in future work (like $(\infty,1)$ and $\infty$-categories, monoidal $(\infty,1)$-categories, $\infty$-operads etc...) }

\subsection{Introduction}

\blockn{At a more technical level, the main contribution of the present work can be described as a construction that allows to produce semi\footnote{See \ref{DefSemiModelCat}}-model categories whose objects are all fibrant and which are ``freely generated by one object'' among such category in a correct homotopy theoretical sense.

We will also set up a formalism that will allow to discuss and study more general sort of ``freely generated'' semi-model categories that we hope might be of help to discuss more general higher categorical structures latter (like $\infty$-categories, $\infty$-operads, etc.).}

\block{This follows the idea that the $(\infty,1)$-category of $\infty$-groupoids is the co-complete $(\infty,1)$-category freely generated by one object and that (combinatorial) model categories are model for locally presentable (hence co-complete) $(\infty,1)$-categories. For this reasons those semi-model categories freely generated by one object will be thought of as semi-model categories of weak $\infty$-groupoids, and as they share a (homotopical) universal property they will automatically be all Quillen equivalent.}

\blockn{In order to perform this construction we will rely heavily on a notion weaker than model category: ``categories with cylinder object'' or in short ``cylinder categories''. They are the categorical opposite of the path categories introduced recently by  B.Van der Berg and I.Moerdijk in \cite{berg2016exact}, which in turn are a strengthening of K.Brown category of fibrant objects (introduced in \cite{brown1973abstract}). See section \ref{SecCylinderCat}.}

\blockn{If (combinatorial) model categories are supposed to model locally presentable $(\infty,1)$-categories, cylinder categories are supposed to model finitely co-complete $(\infty,1)$-categories. Our construction will be done in two step: first we construct a cylinder category ``homotopically freely generated by one object'' and then we ``complete'' it into a semi-model category.

The first step will be achieved by constructing (in section \ref{SecWeakModelStr}) a sort of weak model category whose fibrant objects are the cylinder categories. This allows us to perform all sort of free construction by first considering a genuine algebraic free construction and then take a fibrant replacement, hence producing a cylinder category which satisfies a homotopy theoretic version of the universal property of an actual freely generated object.

The second step is more rigid and not homotopy theoretical in nature: the functor from the category of semi-model categories with fibrant object and left Quillen functors to the category of cylinder categories, which send a semi-model categories to its category of cofibrant object have (up to standard size problems) a left adjoint studied in subsection \ref{SubSecCompletion}. 

Combining the two steps one obtains semi-model categories whose objects are all fibrant which satisfies a weak universal property saying that it is freely generated by one object among semi-model category whose objects are all fibrant with left Quillen functors as morphisms.
}

\block{
The first step should be thought of an analogue of the construction of a Grothendieck coherator and the cylinder category obtained this way will be called a cylinder coherator (see \ref{SubSecCylindCoh}) for this reason. The analogy is made very precise in subsection \ref{SubSecCylindCoh} where we prove that under a technical conjecture a Gorthendieck coherator actually produces a cylinder coherator which defines the same notion of weak $\infty$-groupoids.

The second step will be thought of as looking at the category of groupoids with respect to this cylinder coherator.

In some sense the objects of the cylinder coherator are all the finite ``generalized diagrams'' that are meaningful for our notion of groupoids, the cofibrations between them are the inclusions of diagrams and the more general map are the various composition and coherence operations expected for our notion of weak $\infty$-groupoids.
}

\blockn{While ideas from type theory have been extremely influential in this work,  explicit mention of type theory will be kept to a minimum in the present paper. It will only reappears as an example in the last subsection \ref{SubSecSyntactic}. But in some sense, type theory is hidden in the notion of cylinder category, as it has been shown in \cite{berg2016path} that path categories (the dual of cylinder categories) are in some sense the categorical semantics of a weak intentional type theory where the identity elimination principle has been replaced by a propositional equality instead of a definitional equality. In some sense a cylinder coherator can be thought of as the opposite of the syntactic category of a certain kind of intentional type theory ``freely generated'' by one type (this analogy will be made a little more precise in subsection \ref{SubSecSyntactic}).
}

\blockn{In section \ref{SecCylinderCat}, we will first recall\footnote{The terminology is new and due to the author, but the notion is the exact dual of the notion of path category introduced by I.Moerdijk and B.Van der Berg in \cite{berg2016exact}.} the notion of \emph{cylinder category} and its basic properties. In \ref{SubSecSyntactic} we introduced a notion of equivalence between cylinder categories called \emph{acyclic morphism}. This notion is new but very similar results were already known for Brown category of fibrant objects and due to D-C.Cisinski (see theorem $3.19$ of \cite{cisinski2010categories}).

Finally in subsection \ref{SubSecCompletion} we construct the ``completion'' of a cylinder category into a semi-model category, we prove the universal property of this completion and show that acyclic morphisms of cylinder categories induce Quillen equivalences between the completions.  }

\blockn{In section \ref{SecPreCylinder} we introduced the category of ``pre-cylinder categories'' and study its algebraic property. The underlying idea is that we have only kept the properties of cylinder categories that are ``algebraic'' in nature and forget the weak structures (existence of factorization of morphisms and of retractions of trivial cofibrations) of cylinder categories that were not preserved by morphisms of cylinder categories. Doing so have the effect of providing a category (of pre-cylinder categories) which is much more algebraic and admit all sort of categorical construction, in particular free constructions discussed in \ref{SectionFreeCofCat} and a symmetric closed monoidal structure discussed in \ref{SubSecmonoidalStr}.}

\blockn{In section \ref{SecWeakModelStr} we construct our ``weak model structure'' on the category of pre-cylinder categories whose fibrant object are the cylinder categories and whose weak equivalences between fibrant objects are the acyclic morphisms. A precise summary of the sense in which it is a weak model structure can be found in \ref{ModelStrSumUp}. In subsection \ref{SubSecFibCofib} we introduce the two (partial) weak factorization systems of our model category structure, in subsection \ref{SubSecReedy} we use our notion of fibration of pre-cylinder categories and the monoidal structure on the category pre-cylinder categories to give a new presentation of the construction of ``Reedy type'' cylinder categories (but much of the results here were already known and can be found in \cite{cisinski2010categories} or \cite{radulescu2006cofibrations} with more down to earth proof), we then use this to construct some examples of fibrations and trivial fibrations of cylinder categories that will play the role of a path objects for cylinder categories and are a key component in the proof of the existence of our weak model structure.
In subsection \ref{SubSecWeakModelStrs} we finish the construction of the weak model structure. In subsection \ref{SubSecSlice} we discuss a construction of a ``homotopy slice'' of a cylinder categories that produces new examples of fibrations of cylinder categories, this construction will be used in subsection \ref{SubSecGrothGroupoid}.
}

\blockn{In the final section (\ref{SecCylindCoh}) of the paper we introduce the notion of cylinder coherator. A cylinder coherator (definition \ref{DefCylindCoh}) is a fibrant replacement of the free pre-cylinder categories $F_*$ generated by one object (see \ref{PropFstar}), if $\Ccal$ is a cylinder coherator a $\Ccal$-groupoids is an object of the semi-model categories $\widetilde{\Ccal}$ obtained using the completion procedure of subsection \ref{SubSecCompletion}.
The main result of the paper is theorem \ref{MainThCylindCoh} which gives all the good properties of this notion of weak groupoids, at this point all those properties except the homotopy hypothesis follows trivially from the framework developed in the rest of the paper. In order to prove the analogue of the homotopy hypothesis for $\Ccal$-groupoid we just need to prove it for a single cylinder coherator $\Ccal$, which is done in the next subsection. The small object argument provides examples of cylinder coherator, but those tends not to be very convenient. The last three sub-sections are each devoted to deals with specific cylinder coherators which are simpler than those produced by a naive application of the small object argument:

In subsection \ref{SubSecSimplicialCoh} we construct a cylinder coherator whose associated notion of weak $\infty$-groupoids are semi-simplicial Kan complex for which we are able to prove a comparison theorem with the homotopy category of spaces and hence proving the homotopy hypothesis for all cylinder coherator.

In subsection \ref{SubSecGrothGroupoid} we show that under a simple looking technical conjecture (\ref{conjecture}) the theory of Grothendieck $\infty$-groupoids fits into our framework and hence that conjecture \ref{conjecture} implies Grothendieck homotopy hypothesis. This conjecture is a very concrete statement about Grothendieck's definition which seem considerably weaker than the homotopy hypothesis. A failure of this conjecture would in some sense means that Grothendieck's definition is missing some operations that should be present for weak $\infty$-groupoids (for example, it would miss certain operations that are definable in intentional type theory).
We also show (without any assumption) that any object in a path category (the opposite of a cylinder categories) caries the structure (on its iterated path objects) of a Grothendieck $\infty$-groupoid. Combining this with the result of \cite{berg2016path} that syntactic categories of type theory with propositional identity types are path categories this allows to extend the result of \cite{bourke2016note} that type in intentional type theory are Grothendieck groupoids to the case of a type theory with propositional identity types (the strategy of \cite{bourke2016note} does not extend to this case).

Finally in subsection \ref{SubSecSyntactic} we show that the opposite of the syntactic category of type theory with just propositional identity type and one ``free'' type $*$ is a cylinder coherator. This produces a globular and syntactical definition of weak $\infty$-groupoids which is very similar to the definition given in \cite{brunerie2013syntactic} and \cite{altenkirch2012syntactical} but for which, contrary to Grothendieck's definition, one can prove the homotopy hypothesis, the existence of a semi-model structure of these weak $\infty$-groupoids and one has a definition of the $\infty$-groupoids of weak functors between two $\infty$-groupoids.
We are convinced that this notion of weak $\infty$-groupoids will make sense in intentional type theory, but we think that only a computer implementation would be a convincing proof of this.
}

\subsection{Preliminaries}

\block{\label{DefSemiModelCat}By a semi-model category we mean (following Spitzweck \cite{spitzweck2001operads} definition $2.3$) a structure similar to Quillen model category but where the factorization of an arrow as a trivial cofibration followed by a fibration is only possible for arrow whose domain is cofibrant (all the other axioms being the same). 

A semi-model category is said to be combinatorial if the underlying category is locally presentable and if there is two sets of morphisms $I$ and $J$ (called the generating cofibration and generating trivial cofibration) such that fibrations are the arrow with the right lifting property with respect to $J$ and trivial fibrations are the arrow with the right lifting property with respect to $I$. In this situation we will call ``trivial cofibration'' the map in the class generated by $J$, i.e. the map that have the left lifting property with respect to all fibration and ``acyclic cofibrations'' the map that are cofibration and weak equivalences, those two class agree for arrow whose domain is cofibrant but we cannot say anything for general arrows.

Most of the notion and construction available for ordinary model category carries over easily to semi-model category, we just need to take cofibrant replacement a little more often to avoid the pathological behavior of trivial cofibration with a non-cofibrant domain.
In particular, one still have the same description of the homotopy category in terms of fibrant/cofibrant objects, the notion of Quillen adjunction and Quillen equivalences are still the same and a Quillen equivalence still induces an equivalence between the homotopy categories.
}

\blockn{In diagrams, hooked arrow denoted ($\hookrightarrow$) always denote cofibration (we will not use this notion for monomorphism) and the symbol $\sim$ on a arrow mean that is is a weak equivalence and fibrations are denoted by arrow of the form $\twoheadrightarrow$.}

\section{Cylinder categories}
\label{SecCylinderCat}
\subsection{Definition}

\block{
\Def{A cylinder category $(\Ccal, \cof, \Wcal)$, is a category $\Ccal$, together with two classes of morphisms $\cof$ and $\Wcal$ called respectively the cofibrations and the weak equivalences (and elements of $\Wcal \wedge \cof$ are called the trivial cofibrations) such that:

\begin{enumerate}

\item $\Wcal$ and $\cof$ contains all isomorphisms and are stable under composition.

\item $\Wcal$ satisfies the $2$-out-of-$6$ property: If $f,g,h$ are three composables arrows such that $f \circ g$ and $g \circ h$ are in $\Wcal$ then $f,g,h$ and $f \circ g \circ h$ are in $\Wcal$.

\item $\Ccal$ has an initial object $0$ and for all object $X$ of $\Ccal$ the unique map $0 \rightarrow X$ is a cofibration.

\item Pushout of cofibrations exist and are cofibrations.
\item Pushout of trivial cofibrations are trivial cofibrations.

\item For every object $X \in \Ccal$ there exists a fatorization of the codiagonal:

\[ X \coprod X \hookrightarrow IX \overset{\sim}{\rightarrow} X \]

Where the first map is a cofibration and the second is a weak equivalence. 

\item Every trivial cofibration admit a retraction.

\end{enumerate}

}

}

\blockn{Cylinder categories are exactly the opposites of the ``Path categories'' introduced recently by I.Moerdijk and B.van den Berg in \cite{berg2016exact}. They are also closely related to Brown categories of co-fibrant objects \cite{brown1973abstract}.}

\block{An object $IX$ in a factorization as a cofibration followed by a weak equivalence:

\[X \coprod X \hookrightarrow IX \overset{\sim}{\rightarrow} X \]

is called a cylinder object for $X$. ``$IX$'' will always denotes such a cylinder object for $X$.
}

\block{The $2$-out-of-$6$ property easily implies the more usual $2$-out-of-$3$ property saying that if any two of the arrows $f,g$ and $f\circ g$ are weak equivalences then the third also is. A discussion of the role of this stronger $2$-out-of-$6$ property in comparison with the $2$-out-of-$3$ property can be found in \ref{prop_WeAreHiso}.}

\blockn{Here are some examples of cylinder categories:

\begin{itemize}

\item If $\Ccal$ is a semi-model category in which every object is fibrant, then the full subcategory of co-fibrant objects (with the induced notion of weak equivalences and cofibration) is a cylinder category.

\item The category of $CW$-complexes with homotopy equivalences and the relative $CW$-map as cofibration is a cylinder category (It is almost a special case of the above).

\item The category of finite $CW$-complexes (with the same weak equivalences and cofibratons) is again a cylinder category: indeed cylinder objects for finite $CW$-complex can be chosen to be finite and pushout along cofibrations between finite $CW$-complexes produces finite $CW$-complexes.

\end{itemize}

}

\block{\Def{A morphism $F$ between cylinder categories is a functor $F : \Ccal \rightarrow \Ccal'$ such that:

\begin{itemize}
\item $F$ sends weak equivalences to weak equivalences.
\item $F$ sends cofibrations to cofibrations.
\item $F$ preserve the initial object and the pushout along cofibrations.
\end{itemize}

}

Such morphisms are also called ``exact functor'' in the literature.

Note that there is a terminal cylinder category: the terminal category where the unique arrow is a trivial cofibration, it will be denoted by $0$.
}

\blockn{We conclude this section by introducing the notion of class of generating cofibrations for a cylinder category. It appears that this notion only depends on a small part of the structure of cylinder category and that it will be extremely useful in the rest of the paper to explicit this structre. For this reason we introduce immediately the following notion:}

\block{\label{DefCofibCat}\Def{A cofibration category is a category $\Ccal$ together with a class $\cof_{\Ccal}$ of maps called cofibrations such that:

\begin{itemize}

\item isomorphisms are cofibrations and cofibrations are stable under composition.

\item $\Ccal$ has a initial object $0$ and for all object $X \in \Ccal$ the map $O \rightarrow X$ is a cofibration.

\item Pushout of cofibrations exists and are cofibrations.

\end{itemize}

A morphism of cofibrations categories is a functor that send cofibrations to cofibrations and preserves the initial object and the pushout along cofibrations.
}
Obviously any cylinder category is a cofibration category when forgetting the weak equivalences. The dual notion, ``fibration category'', have also been called ``tribe'' by A.Joyal. 
}

\block{\Def{\begin{itemize}

\item If $\Fcal$ is some class of maps in a category $\Ccal$ one call a iterated pushout of maps in $\Fcal$ a map which is a (finite) compositions of maps that are pushout of maps in $\Fcal$.

\item An object is said to be an iterated pushout of maps in $\Fcal$ if the map $0 \rightarrow X$ is such an iterated pushout.

\item If $\Ccal$ is a cofibration category one say that $\Fcal$ is a class of generating cofibrations if the cofibrations of $\Ccal$ are exactly the iterated pushout of maps in $\Fcal$.

\item A transfinite composition of iterated pushout will be called a transfinite iterated pushout.

\end{itemize}}}

\blockn{For example, in the category of finite $CW$-complexes the inclusion of the $(n-1)$-sphere into the $n$-ball for all $n \geqslant 0$ form a set of generating cofibrations.}

\block{\label{MorphFromGenCofib}\Prop{Let $\Ccal$ be a cofibration category with a set of generating cofibrations, a functor $F : \Ccal \rightarrow \Dcal$ to another cofibration category is a morphism if and only if: it preserves the initial object, it sends the generating cofibrations to cofibrations and it preserves pushout along generating cofibrations.}

\Dem{Such a functor automatically send any pushout of a generating cofibrations to a pushout of a cofibrations. As all cofibrations in $\Ccal$ are composite of pushout of generating cofibrations, $F$ sends all cofibrations to cofibrations. A pushout along a pushout of a generating cofibration is also a pushout along a generating cofibration, hence these pushout are preserved by $F$. A general pushout along a cofibration is obtained as a series of such pushout so they are also preserved by $F$ and this concludes the proof.}
}

\subsection{Homotopy category and generality on cylinder categories}

\blockn{In this subsection, we review for completeness the basic theory of cylinder categories and the construction of the homotopy category. All the material here originally comes from \cite{berg2016exact} and \cite{brown1973abstract} except for some proof that have been slightly improved. }

\block{\Prop{Any morphism in a cylinder category can be factored as a cofibration followed by a weak equivalence. More precisely, any morphism in a cylinder category can be factored as a cofibration followed by a retract of a trivial cofibration. }

This results and the proof that follows are due to K.Brown \cite{brown1973abstract}.

\Dem{Let $F : A \rightarrow B$ be a morphism in $\Ccal$ a cylinder category, we form the co-product $P=IA \coprod_A B$ for $IA$ some cylinder object $A$. $P$ can be alternatively written as the pushout:

\[ P = IA \coprod_{A \coprod A} \left( A \coprod B \right) \]

In particular, the map $A \coprod B \hookrightarrow P$ is a cofibration, and hence one has a cofibration $A \hookrightarrow P$, the map $B \hookrightarrow P$ is a trivial cofibration because it is a pushout of $A \hookrightarrow IA$. There is a map from $P$ to $B$ which is $IA \rightarrow A \rightarrow B$ on $IA$, the identity on $B$, this map is a retract of the trivial cofibration $B \hookrightarrow O$ (hence a weak equivalence by $2$-out-of-$3$), and this form the factorization $A \hookrightarrow P \overset{\sim}{\rightarrow} B$ claimed in the proposition.
}
}

\block{\label{Cor_functorsOnWE}\Cor{

\begin{itemize}

\item Any weak equivalence factors as a trivial cofibration followed by a retract of a trivial cofibration.

\item Any functor between cylinder categories that send trivial cofibrations to trivial cofibrations also send weak equivalences to weak equivalences.

\item A pushout of a weak equivalence along a cofibration is again a weak equivalence.
\end{itemize}
}

Note that as pushout of cofibrations exist and are cofibrations, saying that ``any pushout of a weak equivalence along a cofibration is a weak equivalence'' also means that if $f:A\rightarrow B$ is a weak equivalence between two objects under an object $C$ and $f :C \rightarrow C'$ is a cofibration then the map induced by $f$, $A \coprod_C C' \rightarrow B \coprod_C C'$ is also a weak equivalence

\Dem{If $f$ is a weak equivalence, then as any map, it can be factored as a cofibration followed by a retract of a trivial cofibration, but, by $2$-out-of-$3$ the first cofibration will be a trivial cofibration. Any functor that preserves trivial cofibrations will send this factorization to a factorization as a trivial cofibration followed by a retract of a trivial cofibration, hence to a factorization into two weak equivalences and hence to a weak equivalence.

The same argument applies to pushout along cofibration: pushout of trivial cofibrations are again trivial cofibrations by the fifth axiom of cylinder categories and the pushout of a retract of a trivial cofibration is again a retract of a trivial cofibration.
}

}

\block{\Cor{Let $\Ccal$ be a cylinder category. Let $A$ be an object of $\Ccal$. Then the category $\Ccal(A)$ of cofibrations $A \hookrightarrow B$ (with commutative triangle as morphisms) is a cylinder category with weak equivalences (resp. cofibrations) being the arrow that are weak equivalences (resp. cofibrations) in $\Ccal$.
}

\Dem{All the axioms about pushout, $2$-out-of-$6$, isomorphisms, stability under compositions, and the map $0 \rightarrow X$ being cofibrations are trivial. The existence of retract for trivial cofibration follow from the fact that a retraction of a map under an object $A$ is automatically also under $A$. The only non trivial axiom is the existence of cylinder object, but a cylinder object for $A \hookrightarrow B \in \Ccal(A)$ is obtained simply by taking a factorisation in $\Ccal$ as a cofibration followed by a weak equivalence of the map $B \coprod_A B \rightarrow B$.}
}

\block{\Def{If $A \hookrightarrow B$ is a cofibration in a cylinder category, a relative cylinder object for $B$ (relative to $A$ or under $A$) is an object $I_A B$ such that one has a cofibration/weak equivalence factorization:
\[ B \coprod_A B \hookrightarrow I_A B \overset{\sim}{\rightarrow} B \]

This is the same as a cylinder object for $B$ in the cylinder category $\Ccal(A)$.
 }}

\block{If $A$ and $B$ are two objects of a cylinder category $\Ccal$, and if $v:A\rightarrow B$ is any morphism in $\Ccal$, then one has a pushout functor:

\[ v_{\sharp} : \Ccal(A) \rightarrow \Ccal(B) \]

Which send any object $A \hookrightarrow D$ of $\Ccal(A)$ to $B \hookrightarrow B \coprod_A D$.

$v_{\sharp}$ is a morphism of cylinder categories.

}

\block{\label{Lemsemilifting}\Lem{If in a cylinder category $\Ccal$ one has a square:

\[
\begin{tikzcd}[ampersand replacement=\&]
A \arrow[hookrightarrow]{d}{i} \arrow{r} \&  C \arrow{d}{\sim} \\
B \arrow{r} \& D \\
\end{tikzcd}
\]

where the left vertical arrow is a cofibration and the right vertical arrow is a weak equivalence, then there exists a map $B \rightarrow C$ which makes the upper triangle:

\[
\begin{tikzcd}[ampersand replacement=\&]
A \arrow[hookrightarrow]{d}{i} \arrow{r} \&  C \\
B \arrow{ur} \& \\
\end{tikzcd}
\]

commutes.
}

\Dem{We form the pushout $P:=B \coprod_A C$, let $k:P \rightarrow D$ the natural map and we factors $k$ as a cofibration followed by a weak equivalence:

\[
\begin{tikzcd}[ampersand replacement=\&]
A \arrow[hookrightarrow]{d}{i} \arrow{r} \&  C \arrow[hookrightarrow]{d} \arrow[bend left = 15]{dddrr}{\sim}\& \& \\
B \arrow{r} \arrow[bend right = 15]{rrrdd} \& P \arrow[hookrightarrow]{dr} \& \& \\
\& \& X \arrow{dr}[description]{\sim} \& \\
\& \& \& D \\
\end{tikzcd}
\]

By $2$-out-of-$3$ the cofibration $C \hookrightarrow X$ is a weak equivalence, hence a trivial cofibration and hence admit a retraction. Composing this retraction with the map from $B$ to $X$ produces the map whose existence is claimed in the lemma.
}
}

\block{\Def{\begin{itemize}

\item If $f,g : X \rightrightarrows Y$ are two parallel maps in a cylinder category, one says that $f$ is homotopic to $g$ and one writes $f \sim g$ if the map $(f,g): X \coprod X \rightarrow Y$ can be extended to a map $IX \rightarrow Y$ from some cylinder object, such an extension is called an homotopy between $f$ and $g$.

\item If $i:A \hookrightarrow X$ is a cofibration, and $f,g : X \rightrightarrows Y$ are two parallel maps such that $f\circ i = g \circ i$ one says that $f$ is homotopic to $g$ relative to $A$ and one writes $f \sim_A g$ if the map: $(f,g): X \coprod_A X \rightarrow Y$ can be extended to a map $I_A X \rightarrow Y$ for some relative cylinder object, such an extension is called an homotopy between $f$ and $g$ relative to $A$.

\end{itemize}
}}

\block{\label{PropHomotopyGen}\Prop{ \begin{enumerate}

\item The notion of (relative) homotopy between maps does not depends on the choice of the (relative) cylinder object.

\item Homotopy and relative homotopy are equivalence relations.

\item If $A \hookrightarrow B$ is a trivial cofibration, any two maps from $B$ to $X$ which agrees on $A$ are homotopy equivalent relative to $A$.

\item Homotopy and relative homotopy are compatible with post-compositon in the sense that if $ f \sim_{(A)} g$ then $ u \circ f \sim_{(A)} u \circ g$ for any $u$ such that the composition exists.

\item The homotopy relation is compatible with pre-composition as well, the relative homotopy relation is compatible with pre-composition in the sense that:

if one has a square:

\[
\begin{tikzcd}[ampersand replacement=\&]
 A' \arrow[hookrightarrow]{d} \arrow{r} \& A \arrow[hook]{d} \\
B' \arrow{r}{v} \& B \\
\end{tikzcd}
\]

whose vertical arrows are cofibrations, two morphisms $f,g: B \rightrightarrows X$ which agree on $A$, and such that $f \sim_A g$ then $f \circ v \sim_{A'} g \circ v$.

\end{enumerate}

}

\Dem{
\begin{enumerate}

\item Let $I_A X$ and $I'_A X$ be two (relative) cylinder objects, then considering the square:

\[
\begin{tikzcd}[ampersand replacement=\&]
 X \coprod_A X \arrow[hookrightarrow]{d} \arrow[hookrightarrow]{r} \& I_A X \arrow{d}{\sim} \\
I'_A X \arrow{r}{\sim} \& X \\
\end{tikzcd}
\]

One has by lemma \ref{Lemsemilifting} a map $I'_A X \rightarrow I_A X$ under $X \coprod_A X$ showing that the homotopy relation in the sense of $I_A X$ implies the homotopy relation in the sense of $I'_A X$ and (as it is symmetric) that they are equivalent.

\item The fact that the homotopy relation is reflexive follows from the fact that the codiagonal map can be extended to a map $IX \rightarrow X$ by definition of $IX$. The symmetry follows from the same argument as the point $1.$ with $I'_A X$ being $I_A X$ with the two maps $X \hookrightarrow I_A X$ exchanged. For the transitivity, if $f \sim_A g \sim_A h$ then one has a map: $I_A X \coprod_X I_A X \rightarrow Y$ attesting those two homotopy, but if one consider the diagram:

\[
\begin{tikzcd}[ampersand replacement=\&]
 X \coprod_A X \arrow[bend left = 30]{rr}{(f,h)} \arrow[hookrightarrow]{d} \arrow{r} \& I_A X \coprod_X I_A X \arrow{d}{\sim} \arrow{r} \& Y \\
I_A X \arrow{r} \& X \\
\end{tikzcd}
\]

where the first horizontal map corresponds to the two map $X \rightarrow I_A X$ that are not used in the formation of the coproduct $I_A X \coprod_X I_A X$, then the diagonal map composed with the map to $Y$ produces an homotopy between $f$ and $h$.

\item If $A \hookrightarrow B$ is a trivial cofibration, then the map $B \hookrightarrow B \coprod_A B$ is again a trivial cofibration( pushout of a trivial cofibration) hence by $2$-out-of-$3$ the co-diagonal map $B \coprod_A B \rightarrow B$ is a weak equivalence. Hence $B \coprod_A B$ is already a relative cylinder object for $B$ under $A$, and hence any two compatible maps defines a map $(f,g): B \coprod_A B \rightarrow X$ and hence are homotopy equivalent.

\item An homotopy between two maps is just a map $I_A B \rightarrow X$ hence it can be post-composed with any map.

\item It is enough to prove the relative case, the other case follow by taking $A=A'=0$.

One can then consider the diagram:

\[
\begin{tikzcd}[ampersand replacement=\&]
B' \coprod_{A'} B' \arrow[bend left = 30]{rrr}{(f\circ v, g \circ v)} \arrow[hookrightarrow]{d} \arrow{r}{(v,v)}\& B \coprod_A B  \arrow[hook]{r} \& I_A B \arrow{d}{\sim} \arrow{r}[below]{f \sim_A g } \& X \\
I_{A'} B' \arrow{r} \& B' \arrow{r}{v} \& B \\
\end{tikzcd}
\]

any diagonal map making the upper triangle commutes composed with the map to $X$ produces the desired homotopy. 

\end{enumerate}

}
}

\block{\label{Prop_Semilifting2}\Prop{If in a cylinder category $\Ccal$ one has a square:

\[
\begin{tikzcd}[ampersand replacement=\&]
A \arrow[hookrightarrow]{d}{i} \arrow{r} \&  C \arrow{d}{\sim} \\
B \arrow{r} \& D \\
\end{tikzcd}
\]

where the left vertical arrow is a cofibration and the right vertical arrow is a weak equivalence, then there exists a map $B \rightarrow C$ which makes the upper triangle commutes and the lower triangle commute up to homotopy relative to $A$.
}

This proposition is due to I.Moerdijk and B.Van den Berg, it can be found in \cite{berg2016exact} with a different proof.

\Dem{We use the same construction as in the proof of lemma, one gets a diagram:

\[
\begin{tikzcd}[ampersand replacement=\&]
A \arrow[hookrightarrow]{d}{i} \arrow{r} \&  C \arrow[hookrightarrow]{d}{j} \arrow[bend left = 15]{dddrr}{\sim}\& \& \\
B \arrow{r}{k} \arrow[bend right = 15]{rrrdd} \& P \arrow[hookrightarrow]{dr}{c} \& \& \\
\& \& X \arrow{dr}[description]{\sim} \& \\
\& \& \& D \\
\end{tikzcd}
\]

By $2$-out-of-$3$ the cofibration $c \circ j:C \hookrightarrow X$ is a weak equivalence, hence a trivial cofibration and hence admit a retraction $r$. One has $r \circ c \circ j = Id_C$ and by proposition \ref{PropHomotopyGen} $(3.)$ one has: $ c \circ j \circ r \sim_{C} Id_X$. The map we are trying to construct is $r \circ c \circ k : C \rightarrow B$, the fact that $r \circ c \circ j = Id_C$ gives the commutation of the upper triangle, and the fact that $ c \circ j \circ r \sim_{C} Id_X$ together with proposition \ref{PropHomotopyGen} $(5.)$ and $(4.)$ gives the commutation of the lower triangle up to homotopy relative to $A$.
}
}

\block{\label{PropDefhC}\Prop{If $\Ccal$ is a cylinder category, then there is a category $h\Ccal$ whose objects are the objects of $\Ccal$ and whose arrows are the equivalence class of arrows in $\Ccal$ for the homotopy relation.

Moreover $h \Ccal$ is the localization of $\Ccal$ at weak equivalences, for the natural quotient functor $\Ccal \rightarrow h \Ccal$.
}

\Dem{The fact that $h \Ccal$ is well defined and form a category follows directly from  proposition \ref{PropHomotopyGen}.

Trivial cofibrations in $\Ccal$ are send to isomorphisms in $h\Ccal$: any trivial cofibration $j: A\hookrightarrow B$ admit a retraction $r$, hence one has $r \circ j = Id_A$ and $j \circ r \sim_A Id_B$ because of proposition \ref{PropHomotopyGen} $(3.)$, which implies that $j \circ r \sim B$ (it follows from \ref{PropHomotopyGen} $(5.)$) hence in $h \Ccal$, $r$ is an inverse for $j$.

As any weak equivalence in $\Ccal$ is the composite of a trivial cofibration followed by a retract of a trivial cofibration it implies that any weak equivalences in $\Ccal$ is an isomorphism in $h \Ccal$.

If $F : \Ccal \rightarrow \Dcal$ is any functor which send weak equivalences to isomorphisms, then for any object $X$, the map $IX \rightarrow X$ is sent to an isomorphism, as the two maps $X \rightrightarrows IX$ are section of this map their images are equal in $\Dcal$, and hence any pair of maps that are homotopic in $\Ccal$ have equal images in $\Dcal$, hence the functor from $\Ccal$ to $\Dcal$ can be factored (essentially uniquely) into $h \Ccal$ and this concludes the proof.
}
}

\block{\label{prop_WeAreHiso}\Prop{A morphism in $\Ccal$ is a weak equivalence if and only if its image in $h \Ccal$ is an isomorphisms.

}

Note that the proof of this proposition will be the first place were we use the $2$-out-of-$6$ property (and not just the $2$-out-of-$3$ property). If one starts with a category that satisfies all the axioms of a cylinder category except the $2$-out-of-$6$ property which is replaced by the $2$-out-of-$3$ property, then one can still construct the category $h\Ccal$ without problems. But the present proposition become equivalent to the $2$-out-of-$6$ property and if one call ``homotopy equivalences'' the maps that are invertible in $h\Ccal$ then this endows $\Ccal$ with a new structure of cylinder category (satisfying the $2$-out-of-$6$ property) which is minimal for the $2$-out-of-$6$ property and still defines the same homotopy relation and the same homotopy category. We learned this remarks from \cite[section 7.2]{radulescu2006cofibrations} and it is apparently due to D-C.Cisinski. 

\Dem{A map $f :X \rightarrow Y$ in $\Ccal$ is invertible in $h\Ccal$ if and only if there exists a map $g:Y \rightarrow X$ together with homotopies $f \circ g \sim Id_Y$ and $g \sim f \sim Id_X$. A map homotopic to a weak equivalence is a weak equivalence: indeed if one has a homotopy $a \sim b : IX \rightarrow Y$ with $b$ a weak equivalence then as $b$ is the composite $X \overset{\sim}{\hookrightarrow} IX \rightarrow Y$, the map $IX \rightarrow Y$ is a weak equivalence by the property $2$-out-of-$3$ and hence $a$ is an equivalence as a composite of weak equivalences. Hence in our initial situation, $f \circ g$ and $g\circ f$ are weak equivalences, and hence $f$ and $g$ are also weak equivalences by the $2$-out-of-$6$ property.
}
}

\blockp{The following proposition is often called the gluing lemma or the cube lemma:}

\block{\label{gluingLemma}\Prop{In a cylinder category, if one has a diagram of the form:

 \[
\begin{tikzcd}[ampersand replacement=\&]
B \arrow{d}{\sim}[swap]{v_B} \&  A \arrow{d}{\sim}[swap]{v_A}  \arrow[hookrightarrow]{l}[swap]{i} \arrow{r}{f} \& C \arrow{d}{\sim}[swap]{v_c}  \\
B'  \&  A' \arrow[hookrightarrow]{l}[swap]{i'} \arrow{r}{f'} \& C' \\
\end{tikzcd}
\]

Then the comparison map:

\[ B \coprod_A C \rightarrow B' \coprod_{A'} C' \]

is also a weak equivalence.
}

This results is well known for model category. We learned that it is true in weaker context from \cite[Lemma 1.4.1]{radulescu2006cofibrations}. Our proof is slightly different from the one given in this reference.

\Dem{We first fact $f$ into a cofiration followed by a weak equivalence:

 \[
\begin{tikzcd}[ampersand replacement=\&]
B \arrow{d}{\sim}[swap]{v_B} \&  A \arrow{d}{\sim}[swap]{v_A}  \arrow[hookrightarrow]{l}[swap]{i} \arrow[hook]{r} \& D \arrow{r}{\sim} \& C \arrow{d}{\sim}[swap]{v_c}  \\
B'  \&  A' \arrow[hookrightarrow]{l}[swap]{i'} \arrow{rr}{f'} \& \& C' \\
\end{tikzcd}
\]

And we factor $f'$ into the pushout $\displaystyle D' := A' \coprod_A D$:

 \[
\begin{tikzcd}[ampersand replacement=\&]
B \arrow{d}{\sim}[swap]{v_B} \&  A \arrow{d}{\sim}[swap]{v_A}  \arrow[hookrightarrow]{l}[swap]{i} \arrow[hook]{r} \& D \arrow{d}{\sim} \arrow{r}{\sim} \& C \arrow{d}{\sim}[swap]{v_c}  \\
B'  \&  A' \arrow[hookrightarrow]{l}[swap]{i'} \arrow[hook]{r} \& D' \arrow{r}{\sim} \& C' \\
\end{tikzcd}
\]

The map $D \rightarrow D'$ is a weak equivalence because it is a pushout of a weak equivalence along a cofibration (corollary \ref{Cor_functorsOnWE}) and the map $D' \rightarrow C'$ is a weak equivalence by $2$-out-of-$3$.

In particular one obtains a square:

 \[
\begin{tikzcd}[ampersand replacement=\&]
 D \coprod_A B \arrow{r}{\sim} \arrow{d} \& C \coprod_A B \arrow{d} \\
 D' \coprod_{A'} B' \arrow{r}{\sim} \& C' \coprod_{A'} B' \\
\end{tikzcd}
\]
where the two horizontal arrows are weak equivalences because they are both pushout of a weak equivalence along a cofibration (\ref{Cor_functorsOnWE}). Hence it is enough to prove that the left most arrow is a weak equivalence, but as $D' = D \coprod_A A'$, this leftmost map is in fact isomorphic to the map:

\[ D \coprod_A B \rightarrow D \coprod_A B' \]

Hence the result follows from the fact that $A \hookrightarrow D$ is a cofibration and $B \rightarrow B'$ is a weak equivalence (by a third application of corollary \ref{Cor_functorsOnWE}).
}
}

\subsection{Acyclic morphisms of cylinder categories}
\label{SubSecAcyclic}

\blockn{Results in this sub-section are new, but very similar results where already known for Brown categories of fibrant objects and are due to D-C.Cisinski in \cite{cisinski2010categories} (For example theorem $3.19$). }

\block{\label{DefAcyclic}\Def{A morphism $F : \Ccal \rightarrow \Dcal$ between cylinder categories is said to be:

\begin{itemize}

\item Homotopy surjective if for every object $D$ in $\Dcal$ there exists an object $C$ in $\Ccal$ and a weak equivalence $D\rightarrow F(C)$ (or equivalently $F(C) \rightarrow D$). Equivalently, $F$ is homotopy surjective if $F:h\Ccal \rightarrow h \Dcal$ is essentially surjective.

\item Homotopy fully faithful if for any cofibration $i:A \hookrightarrow B$ in $\Ccal$, for any arrow $x: A \rightarrow X$ and any triangle:

 \[
\begin{tikzcd}[ampersand replacement=\&]
 F(A) \arrow[hookrightarrow]{d}{F(i)} \arrow{r}{F(x)} \& F(X) \\
F(B) \arrow{ur}{v} \& \\
\end{tikzcd}
\]

in $\Dcal$, there exists an arrow $v' : B \rightarrow X$ which makes the triangle:

 \[
\begin{tikzcd}[ampersand replacement=\&]
 A \arrow[hookrightarrow]{d}{i} \arrow{r}{x} \& X \\
B \arrow{ur}{v'} \&  \\
\end{tikzcd}
\]

commutes and such that $F(v')$ is homotopic to $v$ relative to $A$.

\item Acyclic\footnote{We will also sometimes say ``homotopy equivalence'', but we prefer the term acyclic which prevent confusion with the categorical equivalences.} if it is both homotopy surjective and homotopy fully faithful.

\end{itemize}

}}

\block{\label{LemFullyFaithfulByRetract}\Lem{A morphism $F:\Ccal\rightarrow \Dcal$ between cylinder categories is homotopy fully faithful if and only if, for any cofibration $i:A \hookrightarrow B$ in $\Ccal$, if $i$ admit a retraction $r$ in $\Dcal$, then $i$ admit a retraction $r'$ in $\Ccal$ such that $r$ is homotopy equivalent to $F(r')$ relative to $A$.
}

\Dem{The condition stated in the lemma correspond to the definition of homotopy fully faithful given in \ref{DefAcyclic} restricted to the case where the horizontal map $A \rightarrow X$ is the identity. In particular the condition of the lemma is trivially implied by the definition of \ref{DefAcyclic}. But conversely, finding a dashed arrow in a diagram of the form:

 \[
\begin{tikzcd}[ampersand replacement=\&]
 A \arrow[hookrightarrow]{d}{i} \arrow{r}{u} \& X \\
B \arrow[dashed]{ur}{v} \&  \\
\end{tikzcd}
\]

Is exactly the same as finding a retraction to the cofibration:

 \[ X \hookrightarrow X \coprod_A B \]

Moreover, an homotopy relative to $A$ between two such dashed arrow is the same as a factorization:

\[ X \coprod_A B \coprod_A B \hookrightarrow X \coprod_A I_A B \rightarrow X \]

But $X \coprod_A I_A B$ is a relative cylinder object for $X \hookrightarrow X \coprod_A B$ hence this is exactly the same as a homotopy relative to $X$ of the two retractions. Hence this proves the lemma.
}
}

\block{\label{Charac_FullyFaithful}\Cor{For a morphism $F:\Ccal \rightarrow \Dcal$ the following conditions are equivalent:

\begin{enumerate}

\item $F$ is homotopy fully faithful
\item For all object $A$ in $\Ccal$ the functor $hF_A: h\Ccal(A) \rightarrow h \Dcal(A)$ is full.
\item For all object $A \in \Ccal$ the functor $F_A:\Ccal(A) \rightarrow \Dcal(A)$ is homotopy fully faithful.
\item For all object $A \in \Ccal$ the functor $hF_A:h\Ccal(A) \rightarrow h\Dcal(A)$ is fully faithful.
\end{enumerate}
}

\Dem{Unfolding the construction of $h\Ccal(A)$ and $h\Dcal(A)$ in the condition $2.$, show that it is equivalent to the definition \ref{DefAcyclic} of ``homotopy fully faithful'' in the special case where the horizontal map $A \rightarrow X$ is a cofibration: Indeed, they both says that for all object $B$ and $X$ for $\Ccal(A)$ if there is a map between them in $h\Dcal$ then it is homotopy equivalent relative to $F(A)$ to a map in $\Ccal(A)$.

But because of the lemma \ref{LemFullyFaithfulByRetract}, it is enough to check the definition of homotopy fully faithful when $A \rightarrow X$ is the identity, hence in particular a cofibration. Hence conditions $1.$ and $2.$ are equivalent.

\bigskip

Using condition $2.$ one see that being homotopy fully faithful is stable under co-slice, hence these conditions are also equivalent to condition $3.$.

\bigskip

In order to prove that they are also equivalent to condition $4.$ we just need to prove that if $F$ is homotopy fully faithful then $hF$ is faithful, but this follow immediately from the definition of homotopy fully faithful applied to:

 \[
\begin{tikzcd}[ampersand replacement=\&]
 A \coprod A \arrow[hookrightarrow]{d}{i} \arrow{r}{(v,v')} \& X \\
IA \arrow[dashed]{ur}{h} \&  \\
\end{tikzcd}
\]

Indeed, if $v$ and $v'$ (two maps in $\Ccal$) are homotopy equivalent in $\Dcal$ (hence equal in $h\Dcal$), then there is such a dashed arrow in $\Dcal$, hence it can be lifted to $\Dcal$ and hence $v$ and $v'$ are homotopy equivalent in $\Ccal$ (hence equal in $h\Ccal$).

}
}

\block{\label{Charac_homotopyequiBetweenCylinder}\Prop{For a morphism $F: \Ccal \rightarrow \Dcal$ of cylinder categories, the following conditions are equivalent:

\begin{enumerate}

\item $F$ is acyclic.

\item $F$ detect weak equivalences\footnote{i.e. if $F(v)$ is a weak equivalence is and only if $v$ is. } and for all object $A$ of $\Ccal$ the functor $ F_A :\Ccal(A) \rightarrow  \Dcal(F(A))$ is homotopy surjective.

\item The functor $hF : h\Ccal \rightarrow h\Dcal$ is an equivalence of category.

\item For all object $A$ of $\Ccal$, the functor: $h F_A: h \Ccal(A) \rightarrow h \Dcal (F(A)) $ is an equivalence of category.

\end{enumerate}
}

\Dem{We start with $3. \Rightarrow 2.$ 
If $hF$ is an equivalence of categories, it in particular detect isomorphisms, but as weak equivalences are exactly the isomorphisms of the homotopy category, this implies that $F$ detect weak equivalences.

So we need to prove that $F_A$ is homotopy surjective.

Let $F(A) \hookrightarrow B$ be an object of $\Dcal(F(A))$.
As $hF$ is essentially surjective, there exists an object $B'$ of $\Ccal$ such that $F(B')$ is weakly equivalent to $B$. One can form in $h\Dcal$ a triangle:

 \[
\begin{tikzcd}[ampersand replacement=\&]
 F(A) \arrow{d} \arrow{r} \& F(B') \arrow{dl}{\sim} \\
B \& \\
\end{tikzcd}
\]

which hence is a triangle in $\Dcal$ that commutes up to homotopy. Moreover, as $hF$ is fully faithfull The map $F(A) \rightarrow F(B)$ can be freely assumed to be a map of the form $F(v)$, and by replacing this $v$ be a cofibration/weak equivalence factorization, one can moreover assume that $v$ is a cofibration (by changing the object $B'$). Rewriting this diagram, by explicitly writing the homotopy that makes it commute on gets:

 \[
\begin{tikzcd}[ampersand replacement=\&]
 F(A) \arrow[hookrightarrow]{d}{\sim} \arrow[hook]{r}{F(v)} \& F(B') \arrow{dd}{\sim} \\
 F(IA)  \arrow{dr}\& \\
F(A)  \arrow[hookrightarrow]{u}{\sim} \arrow[hookrightarrow]{r} \& B \\
\end{tikzcd}
\]

Then, forming the pushout $P= B' \coprod_A IA$ one gets a diagram:

 \[
\begin{tikzcd}[ampersand replacement=\&]
 F(A) \arrow[hookrightarrow]{d}{\sim}[swap]{i_0} \arrow[hookrightarrow]{rr}{F(h)} \& \& F(B') \arrow[hookrightarrow]{ld}{\sim} \arrow{dd}{\sim} \\
 F(IA)  \arrow{drr} \arrow[hookrightarrow]{r} \& F(P) \arrow{dr}\& \\
F(A)  \arrow[hookrightarrow]{u}{i_1} \arrow[hookrightarrow]{rr}\& \& B \\
\end{tikzcd}
\]

The map $F(A) \overset{i_1}{\hookrightarrow} F(IA) \rightarrow F(P)$ is a cofibration of the form $F( \_ )$ and the map from $F(P)$ to $B$ is a weak equivalence by $2$-out-of-$3$, hence one has constructed an object $A \hookrightarrow P$ in $\Ccal$ such that $F(P)$ is homotopy equivalent to $B$ in $\Dcal(A)$ which concludes the proof of this implication.

\bigskip

We now prove that $2. \Rightarrow 1.$.
We only need to proves that $F$ is homotopy fully faithful. We will use the criterion in Lemma \ref{LemFullyFaithfulByRetract}.

Let $i: A \hookrightarrow B$ a cofibration in $\Ccal$ and assume that it has a retraction $r:F(B) \rightarrow F(A)$ in $\Dcal$. We consider a cofibration/weak equivalences factorization of $r$ in order to a get a factorization of the identity of $F(A)$ as:

\[F(A) \hookrightarrow F(B) \hookrightarrow C \overset{\sim}{\rightarrow} A \]

As $F_B$ is homotopy surjective, there exists a cofibration $j:B \hookrightarrow C'$ such that one has a diagram:

 \[
\begin{tikzcd}[ampersand replacement=\&]
 F(A) \arrow[hook]{r}{F(i)} \& F(B) \arrow[hook]{dr}[swap]{F(j)} \arrow[hook]{r} \& C \arrow{r}{\sim} \& A \\
  \& \& C' \arrow{u}{\sim}
\end{tikzcd}
\]

In particular $F(j) \circ F(i)$ is a weak equivalence in $\Dcal$ and, as $F$ detect weak equivalences, $j \circ i$ is a trivial cofibration, hence it admit a retraction $r':C' \rightarrow A$. This already produces a retraction $r''=r \circ j:B \rightarrow A$ of $i$ in $\Ccal$, in order to conclude one need to show that $F(r'')$ is homotopy equivalent to $r$ relative to $A$. But in $\Dcal$, $F(j \circ i)$ already has a retraction $r'''$ obtained by composing (in the above diagram) the equivalence $\lambda :F(C') \rightarrow C$ to the equivalence $C \rightarrow F(A)$. When restricted to $F(B)$ this retraction induces the retraction $r:F(B) \rightarrow F(A)$  we started from, and any two retractions of a trivial cofibration $F(A) \hookrightarrow F(C')$ are homotopy equivalent relative to $F(A)$ which proves by restricting\footnote{By restriction of homotopies, we are referring to point $5.$ of \ref{PropHomotopyGen}} to $F(B)$ that $r = r'''|_B$ and $F(r'') = F(r')|_B$ are homotopy equivalent relative to $F(A)$ which concludes the proof of $2. \Rightarrow 1.$.

At this point, all the other implications become trivial: Indeed if $F$ is acyclic then $hF$ is essentially surjective by definition and fully faithful by corollary \ref{Charac_FullyFaithful}, hence $1. \Rightarrow 3.$. So the first three conditions are equivalent. The condition $4.$ tautologically implies $3.$ and follows from $2.$ together with $1.$ because of corollary \ref{Charac_FullyFaithful}.

}

}

\block{\label{vsharpEquiv}\Prop{If $v :A \rightarrow B$ is a weak equivalence in a cylinder category then the morphism $v_{\sharp}: \Ccal(A) \rightarrow \Ccal(B)$ is acyclic.
}

\Dem{We will use characterization $2.$ of proposition \ref{Charac_homotopyequiBetweenCylinder}. For each object of $\Ccal(A)$ there is a weak equivalence from $X$ to $v_{\sharp}(X)$ given by the pushout defining $v_{\sharp}$ hence $v_{\sharp}$ detect weak equivalences. As functors of the form $(v_{\sharp})_X$ corresponds to $w_{\sharp}$ for $w$ the pushout of $v$ along the cofibration $A \hookrightarrow X$, it is enough to prove that $v_{\sharp}$ is homotopy full. For any object $B \hookrightarrow X$ of $\Ccal(B)$, a factorization as cofibration followed by a weak equivalence of $A \rightarrow B \hookrightarrow X$ produces an object of $\Ccal(A)$ which will be sent to an object homotopy equivalent to $X$ in $\Ccal(B)$ hence this concludes the proof of the first implication.
}
}

\block{\Prop{Homotopy equivalences between cylinder categories satisfies the  $2$-out-of-$6$ property.}

\Dem{This follows from characterization $3.$ in proposition \ref{Charac_homotopyequiBetweenCylinder} and the fact that equivalences of categories satisfies this condition.
}
}

\subsection{Completion of cylinder category in semi-model categories}
\label{SubSecCompletion}
\blockn{In this section we show how a \emph{small} cylinder category can be naturally completed into a semi-Model category. In $\infty$-categorical therms, this construction corresponds to freely adding direct co-limit to a finitely co-complete small $(\infty,1)$-category in order to turn it into a co-complete $(\infty,1)$-category, and this is usually done by looking at the category of presheaves sending finite co-limit to finite projective limits. }

\block{\Def{Let $\Ccal$ be a small cofibration category, we denote by $\widetilde{\Ccal}$ the category of presheaves over $\Ccal$ that send the initial object to the singleton and pushout along cofibrations to pullback of sets.}

As $\Ccal$ is small, $\widetilde{\Ccal}$ is locally presentable category, in particular it is complete and co-complete. In fact $\widetilde{\Ccal}$ can be described as the category of models of some finite limits sketches (as in \cite[chapter 4]{barr1985toposes}) which are equivalent to the ``Partial Horn theory'' or the ``quasi-equational theory'' of \cite{palmgren2007partial}, i.e. it is algebraic in some sense.

Moreover, limits and directed colimits in $\Ccal$ are computed ``objectwise'' (i.e. as limits and directed colimits of presheaves), finite (non directed) colimits are considerably harder to compute.

For any cofibration category $\Ccal$, the (covariant) Yoneda embeddings of $\Ccal$ takes values in $\widetilde{\Ccal}$ and $\Ccal$ will always be identified with its image though the Yoneda embeddings $\Ccal \rightarrow \tilde{\Ccal}$.

\Def{If $F:\Ccal \rightarrow \Dcal$ is a morphism between two cofibrations categories, we denote by $R_F:\tilde{\Dcal} \rightarrow \tilde{\Ccal}$ the functor of restriction along $F$ and by $\tilde{F}$ its left adjoint (the left Kan extension of $F$).}

$\tilde{F}$ is always an extension of $F$ when $\Ccal$ is identified with the image of the Yoneda embeddings $\Ccal \rightarrow \widetilde{\Ccal}$.

}

\blockp{Completion of small cofibrations categories comes with a notion of cofibration and even a weak factorization system:}

\block{\label{CofInTildeC}\Def{
An arrow $f$ in $\widetilde{\Ccal}$ will be called:
\begin{itemize}

\item A trivial fibration if it has the right lifting property with respect to all cofibrations in $\Ccal$.

\item A finite cellular map if it is a iterated pushout of cofibrations in $\Ccal$.

\item A cellular map if it is a transfinite iterated pushout of cofibrations in $\Ccal$.

\item A cofibration if it is a retract of a cellular map.

\end{itemize} 

}

By the small object argument, and because $\Ccal$ is assumed to be small, every arrow in $\widetilde{\Ccal}$ can be factored as a cellular map followed by a trivial fibration and cofibrations are exactly the maps which have the left lifting properties with respect to all trivial fibrations.
Cofibrations and trivial fibrations form a weak factorization system.

Moreover it is easy to see that an object $X$ of $\widetilde{\Ccal}$ is in $\Ccal$ if and only if the map $0 \rightarrow X$ is a finite cellular map.
Note that unless we require $\Ccal$ to be Cauchy complete\footnote{i.e. that all idempotent of $\Ccal$ splits.} and that cofibrations of $\Ccal$ are stable under retract cofibrations in $\Ccal$ might be different from cofibrations in $\widetilde{\Ccal}$: cofibrations in $\Ccal$ corresponds in $\widetilde{\Ccal}$ to finite cellular map between representable objects and they are in particular cofibrations in $\widetilde{\Ccal}$.
}

\block{\label{LemCellularObjectDirect}\Lem{Every cellular objects in $\tilde{\Ccal}$, i.e. objects $X$ such that the map $0 \rightarrow X$ is cellular, can be written as a directed co-limit of representable such that all the transition map are cofibrations in $\Ccal$ (i.e. finite cellular maps). }

\Dem{Let $X$ be a cellular object, it means that $X$ can be written as the colimit of a sequence of objects $X_{\alpha}$ induced by an ordinal $\rho$ such that for each limit ordinal $\alpha < \rho$, $X_{\alpha}$ is the colimit of the $X_{\beta}$ for $\beta <\alpha$ and for each successor ordinal $\alpha^+$, $X_{\alpha^+}$ is obtained from $X_{\alpha}$ as a pushout:

\[ X_{\alpha^+} = X_{\alpha} \coprod_A B\]

for $A \hookrightarrow B$ a cofibration in $\Ccal$.

We will prove by induction on $\alpha$ that $X_{\alpha}$ is a directed colimit as claimed in the lemma, and that moreover the map $X_{\alpha'} \rightarrow X_{\alpha}$ for $\alpha' < \alpha$ corresponds to an enlargement of the category indexing this colimits.

The case $\alpha =0$ is tautological (with $X_{0} = 0$ ).
If $\alpha$ is a limit ordinal then the result is also easy: one just need to take the increasing union of all the indexing categories of all the colimit of the $X_{\alpha'}$ for $\alpha' < \alpha$.
Assume the result holds for $X_{\alpha}$, i.e. $X_{\alpha}= \lim_{i \in I} X_i$ where $I$ is a direct category, the $X_i$ are in $\Ccal$ and the transition map between the $X_i$ are cofibrations in $\Ccal$. We will write $X$ for $X_{\alpha}$ and $X^+$ for $X_{\alpha^+}$.

As directed colimit in $\tilde{\Ccal}$ are computed object wise, the map $A \rightarrow X$ defining $X^+ = X \coprod_A B$ factors into one of the $X_i$.
Let us consider the following diagram: it has two type of objects: the $X_i$ for $i \in I$ and the $Y_{a,i}$ for $i\in I$ and $a :A \rightarrow X_i$ an arrow inducing the correct map from $A$ to $X$. The object $Y_{a,i}$ will be sent to $X_i \coprod_{A}B$, morphisms from $X_i$ to $X_j$ stay the same as in the original diagram, there is no morphisms from $Y_{i,a}$ to $X_j$, morphisms from $X_i$ to $Y_{j,a}$ are the morphisms from $X_i$ to $X_j$ in the original diagram and morphisms from $Y_{i,a}$ to $Y_{j,a'}$ are the morphisms from $X_i$ to $X_j$ in the original diagram that are compatible to  the maps from $A$ (in $\Ccal$) and hence induces a cofibration from $X_i \coprod_{A}B$ to $X_j \coprod_{A}B$. It is a routine check to see that this is a directed diagram and that its colimit is indeed $X^+$ (for example by checking the universal property) and this concludes the proof.
}
}

\block{\label{DefFibWeakEq}\Def{If $\Ccal$ is a cylinder category, a morphism $f :X \rightarrow Y $ in $\tilde{\Ccal}$ is called:

\begin{itemize}

\item A fibration if it has the right lifting properties with respect to all trivial cofibrations in $\Ccal$.

\item A trivial cofibration if it has the left lifting property with respect to all fibrations.

\item A weak equivalence if if for every square:

\[
\begin{tikzcd}[ampersand replacement=\&]
A \arrow[hookrightarrow]{d}{i} \arrow{r}  \&  X \arrow{d} \\
B \arrow{r} \& Y \\
\end{tikzcd}
\]

where $i: A \hookrightarrow B$ is a cofibration in $\Ccal$, there exists applications $a:B\rightarrow X$ and $h:I_A B \rightarrow Y$ for some relative cylinder object $I_A B$ which extend the above square into the following commutative diagram:

\[
\begin{tikzcd}[ampersand replacement=\&]
A \arrow[hookrightarrow]{ddd}{i} \arrow[hookrightarrow]{dr}{i} \arrow{rr} \& \&  X \arrow{ddd} \\
\& B \arrow[hookrightarrow]{d} \arrow{ur}{a} \& \\
\& I_A B \arrow{dr}{h} \& \\
B \arrow{rr} \arrow[hookrightarrow]{ur}\& \& Y \\
\end{tikzcd}
\]

\end{itemize}

}

By the small object argument, every arrow can be factored as a trivial cofibration followed by a fibration, and trivial cofibrations are exactly the retract of transfinite iterated pushout of trivial cofibrations in $\Ccal$.

The filling condition in the definition of weak equivalence means that $a$ is a diagonal filler:

\[
\begin{tikzcd}[ampersand replacement=\&]
A \arrow[hookrightarrow]{d}{i} \arrow{r}  \&  X \arrow{d} \\
B \arrow{ur}{a} \arrow{r} \& Y \\
\end{tikzcd}
\]

Where the upper triangle is commutative and the lower triangle commute up to the homotopy $h$ relative to $A$. In particular, weak equivalences in $\Ccal$ are weak equivalences in $\widetilde{\Ccal}$ by proposition \ref{Prop_Semilifting2}. The converse is also true and will be proved in \ref{Prop_SameWe}.

}

\block{\Lem{The existence of a filling as in the definition of weak equivalences \ref{DefFibWeakEq} does not depends on the choice of the relative cylinder object $I_A B$.}

\Dem{It follows from the fact that any two cylinder objects $I_A B$ and $I'_A B$ are connected by maps compatible to the cofibrations from $B \coprod_A B$.}
}

\block{\label{Prop_trivFibInComp}\Prop{Trivial fibrations in $\widetilde{\Ccal}$ are exactly the fibrations that are also weak equivalences.}

\Dem{A trivial fibration is obviously a fibration, and satisfies the lifting condition defining weak equivalences (\ref{DefFibWeakEq}) striclty (not up to homotopy) hence is a weak equivalence. Conversely, if $f$ is both a weak equivalence and a fibration, then for any square:

\[
\begin{tikzcd}[ampersand replacement=\&]
A \arrow[hookrightarrow]{d}{i} \arrow{r}  \&  X \arrow[twoheadrightarrow]{d} \\
B \arrow{r} \& Y \\
\end{tikzcd}
\]

one can first use that $f$ is a weak equivalence to obtain the arrow $a$ and $h$ and then that $f$ is a fibration and $B \hookrightarrow I_A B$ is a trivial cofibration to obtain the arrow $b$ in the following commutative diagram:

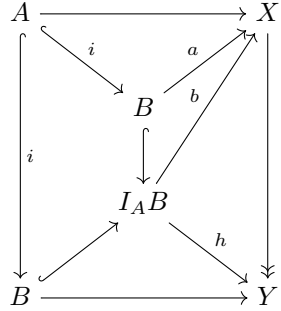
\[
\begin{tikzcd}[ampersand replacement=\&]
A \arrow[hookrightarrow]{ddd}{i} \arrow[hookrightarrow]{dr}{i} \arrow{rr} \& \&  X \arrow[twoheadrightarrow]{ddd} \\
\& B \arrow[hookrightarrow]{d} \arrow{ur}{a} \& \\
\& I_A B \arrow{dr}{h} \arrow{uur}{b}\& \\
B \arrow{rr} \arrow[hookrightarrow]{ur}\& \& Y \\
\end{tikzcd}
\]

}

}

\blockn{In order to study properties of weak equivalences it will be convenient to introduce the following definition.}

\block{\Def{Let $X$ be an object in $\tilde{\Ccal}$, $i:A \hookrightarrow B$ a cofibration in $\Ccal$ and $x: A \rightarrow X$ a morphism, then one defines:

\[ \pi_{B/A}(X,x) = \left\lbrace f : B \rightarrow X | f \circ i = x \right\rbrace / \sim_A \]

Where $\sim_A$ is the homotopy relation\footnote{Extended to this situation in the obvious way: $f \sim_A f'$ iff the map $(f,f'): B \coprod_A B \rightarrow X$ can be extended into a map $h:I_A B \rightarrow X$ for some relative cylinder object $I_A B$.} relative to $A$.

A morphism $f :X \rightarrow X'$ induces a map: 

\[\pi_{B/A}(f,x) :  \pi_{B/A}(X, x) \rightarrow \pi_{B/A}(X',f\circ x). \]

}

\Lem{ The extension of $\sim_A$ above is indeed an equivalence relation and does not depends on the choice of a relative cylinder object $I_A B$.}
\Dem{The proof is exactly the same as the proof given in \ref{PropHomotopyGen} that in a cylinder category the relative homotopy relation is an equivalence relation and does not depend on the choice of the relative cylinder object.}
}

\block{\label{PiCharacOfEquiv}\Prop{For a morphism $f : X \rightarrow X'$ the following conditions are equivalent:

\begin{enumerate}

\item For all cofibrations $A \hookrightarrow B$ in $\Ccal$, and all $x:A \rightarrow X $, the map: 
\[\pi_{B/A}(f,x) :  \pi_{B/A}(X, x) \rightarrow \pi_{B/A}(X',f\circ x) \]
is bijective.

\item for all cofibration $A \hookrightarrow B$ in $\Ccal$, and all $x:A \rightarrow X $, the map: 
\[\pi_{B/A}(f,x) :  \pi_{B/A}(X, x) \rightarrow \pi_{B/A}(X',f\circ x) \]
is surjective.

\item $f$ is a weak equivalence in the sense of definition \ref{DefFibWeakEq}.

\end{enumerate}

}
\Dem{the implication $1. \Rightarrow 2.$ is tautological. The equivalence between $2.$ an $3.$ is relatively easy:

A square:

\[
\begin{tikzcd}[ampersand replacement=\&]
A \arrow[hookrightarrow]{d}{i} \arrow{r}{x}  \&  X \arrow{d}{f} \\
B \arrow{r}{v} \& X' \\
\end{tikzcd}
\]

gives an element $v$ of $\pi_{B/A}(X',f(x))$ (and every element corresponds to such a square) and a filing as in definition of weak equivalences in \ref{DefFibWeakEq} is exactly the data of map $B \rightarrow X$ that would make the upper triangle commutes and the lower triangle commutes up to an homotopy relative to $A$, i.e. an element of $\pi_{B/A}(X,x)$ whose image by $f$ is $v$. Hence the claim that such a filling exists is exactly the claim that $\pi_{B/A}(f,x)$ is surjective.

It remains to prove that $3. \Rightarrow 1.$. Assume that $f$ is a weak equivalence. We just need to prove that $\pi_{B/A}(f,x)$ is a monomorphism. Let $u,v \in \pi_{B/A}(X,x)$ such that $f \circ u \sim_A f \circ v$, hence one has a square:

\[
\begin{tikzcd}[ampersand replacement=\&]
\displaystyle B \coprod_A B \arrow[hookrightarrow]{d}{i} \arrow{r}{(u,v)}  \&  X \arrow{d}{f} \\
I_A B \arrow{r} \& X' \\
\end{tikzcd}
\]

Using that $f$ is a weak equivalence one gets a commutative triangle:

\[
\begin{tikzcd}[ampersand replacement=\&]
\displaystyle B \coprod_A B \arrow[hookrightarrow]{d}{i} \arrow{r}{(u,v)}  \&  X \\
I_A B \arrow{ur}  \\
\end{tikzcd}
\]

which proves that $u=v$ in $\pi_{B/A}(X,x)$ and hence concludes the proof.

}

}

\block{\label{PiLeftTransport}\Prop{\begin{itemize}
\item For $A \in \Ccal$, $X \in \widetilde{\Ccal}$ and $x:A \rightarrow X$

\[ B \mapsto \pi_{B/A}(X,x) \]

is a contravariant functor from $\Ccal(A)$ to Sets which send weak equivalences to isomorphisms, in particular it factor as a functor on $h \Ccal(A)$.

\item If $A \hookrightarrow B$ is a cofibration in $\Ccal$, $v:A \rightarrow A'$ is a map in $\Ccal$, $B'$ is a pushout $B \coprod_A A'$ and $x:A' \rightarrow X$ is a map in $\widetilde{\Ccal}$ then the natural map:

\[ \pi_{B'/A'}(X,x) \rightarrow \pi_{B/A}(X,x\circ v) \]

is a bijection.

\item If one has a commutative square:

 \[
\begin{tikzcd}[ampersand replacement=\&]
A \arrow[hookrightarrow]{d}{i} \arrow{r}{f} \& A' \arrow[hookrightarrow]{d}{i'} \& \\
B \arrow{r}{g} \& B' \\
\end{tikzcd}
\]

where either $f$ and $g$ are weak equivalences, or the map $B \coprod_A A' \rightarrow B'$ is a weak equivalence, then for any $x: A' \rightarrow X$ in $\widetilde{\Ccal}$ one has a bijection:

\[ \pi_{B'/A'}(X,x) \rightarrow \pi_{B/A}(X,x \circ v) \]

\end{itemize}

Moreover all the construction done here are functorial in $X$ in the appropriate sense.

}

\Dem{\begin{itemize}

\item The functoriality follows from the fact that the (extended) relative homotopy relation is stable under precomposition (same proof as in proposition \ref{PropHomotopyGen}) and the proof that weak equivalences induce bijections is essentially the same as the proof of the fact that weak equivalences in $\Ccal$ are bijections in $h \Ccal$ given in the proof of proposition \ref{PropDefhC}. Functoriality in $X$ is also immediate.

\item A map from $B'$ to $X$ extending the map from $A'$ to $X$ is, because of the universal property of the pushout, the same as a map from $B$ to $X$ extending the map from $A$ to $X$. Moreover as $A' \coprod_{A} I_A B$ can be used as a relative cylinder object for $A' \hookrightarrow B'$ the homotopy relation is the same in the two descriptions, hence one has the bijection. Once again the functoriality in $X$ is clear.

\item If $f$ and $g$ are weak equivalences then the $2$-out-of-$3$ property implies that the comparison map from the pushout to $B'$ is a weak equivalence, the result follows then from the conjunction of the two points above.

\end{itemize}

}
}

\block{\label{PiLeftInvariance}\Cor{If $x,x':A \rightrightarrows X$ are two arrows, and $h: IA \rightarrow X$ is an homotopy between $x$ and $x'$ then there is a bijection:

\[ h^* : \pi_{B/A}(X,x) \rightarrow \pi_{B/A}(X,x') \]

Moreover, if $f: X \rightarrow Y$ is a morphism in $\tilde{\Ccal}$ the square:

\[
\begin{tikzcd}[ampersand replacement=\&]
\pi_{B/A}(X,x) \arrow{d}{\pi_{B/A}(f,x)} \arrow{r}{h^*}  \&  \pi_{B/A}(X,x') \arrow{d}{\pi_{B/A}(f,x')} \\
\pi_{B/A}(X',f \circ x) \arrow{r}{(f \circ h)^*} \& \pi_{B/A}(X',f \circ x') \\
\end{tikzcd}
\]

commutes.

In particular, in order to test whether a map $f :X \rightarrow X'$ is a weak equivalence, we only need to test that the maps $\pi_{B/A}(f,x)$ are surjective for one $x$ in each homotopy class of map $x:A \rightarrow X$.

}

\Dem{We choose a factorization as cofibration followed by a weak equivalence:

\[B \coprod B \hookrightarrow \left( B \coprod B \right) \coprod_{A \coprod A} IA \hookrightarrow IB \overset{\sim}{\rightarrow} B\]

$IB$ is indeed a cylinder object for $B$ and one has a diagram of cofibrations:

\[
\begin{tikzcd}[ampersand replacement=\&]
A \arrow[hook]{d}[swap]{\sim} \arrow[hook]{r} \&  B \arrow[hook]{d}{\sim} \\
IA \arrow[hook]{r} \& IB  \\
A  \arrow[hook]{u}{\sim} \arrow[hook]{r} \& B \arrow[hook]{u}[swap]{\sim}
\end{tikzcd}
\]

Applying the construction of \ref{PiLeftTransport}, if $h$ is an homotopy between two morphism $x,x' : A \rightrightarrows X $ one gets two isomorphisms:

\[ \pi_{B/A}(X,x ) \overset{\sim}{\leftarrow} \pi_{IB/IA}(X,h) \overset{\sim}{\rightarrow} \pi_{B/A}(X,x') \]

There composite give us the transport map $h^*$ we are looking for. The compatibility with transport along a map $f :X \rightarrow X'$ follow from the compatibility of the construction explained in \ref{PiLeftTransport}.

}
}

\block{\label{CompWE2out3}\Cor{Weak equivalences in $\tilde{\Ccal}$ satisfies the $2$-out-of-$3$ condition.}

\Dem{Let $X \overset{f}{\rightarrow} Y \overset{g}{\rightarrow} Z$ be two composable arrows in $\tilde{\Ccal}$.

if $f$ and $g$ are weak equivalences or if $g \circ f$ and $g$ are weak equivalences, then the third map also is because of the characterizations $1.$ of weak equivalences in proposition \ref{PiCharacOfEquiv}.

The last case need a little more care: if $g \circ f$ and $f$ are weak equivalences, then, by the same argument, one can conclude that for all $A \hookrightarrow B$ and for all $x : A \rightarrow X$, the map $\pi_{B/A}(g,f \circ x)$ is a bijection. But as $f$ is a weak equivalence the map $\pi_{A/0}(f,0)$ is a surjection, hence for any $y : A \rightarrow Y$ there exists a $x : A \rightarrow X$ such that $y$ is homotopic (without base) to $f \circ x$. Hence, by \ref{PiLeftInvariance} this proves that $g$ is a weak equivalence and concludes the proof.

}
}

\block{\label{Prop_SameWe}\Prop{A map in $\Ccal$ is a weak equivalence in $\Ccal$ if and only if it is a weak equivalence in $\widetilde{\Ccal}$.}

\Dem{We already saw that a weak equivalence in $\Ccal$ is a weak equivalence in $\widetilde{\Ccal}$ (this is basically proposition \ref{Prop_Semilifting2}). For the converse, using the cofibration/weak equivalence factorization in $\Ccal$ and the $2$-out-of-$3$ property for weak equivalences in $\widetilde{\Ccal}$ it is enough to prove the result for cofibrations in $\Ccal$. Let $i:A \hookrightarrow B$ be a cofibration in $\Ccal$ which is a weak equivalence in $\widetilde{\Ccal}$, then using that $\pi_{B/A}(i,Id_A)$ is surjective one can find a retraction $p:B \rightarrow A$ such that $i \circ p$ is homotopy equivalent to $Id_B$ relative to $A$, but an inverse of $i$ on $h\Ccal$ hence by proposition \ref{prop_WeAreHiso} this proves that $i$ is a weak equivalence and concludes the proof.

}}

\block{\Prop{If $f :X \rightarrow Y$ is an arrow in $\tilde{\Ccal}$ with a cofibrant domain then $f$ is a trivial cofibration if and only if it is a cofibration and a weak equivalence. }

\Dem{
Let $X$ be a cofibrant object, and let $x:C \rightarrow X$ be an arrow with $C \in \Ccal$, and $C \overset{\sim}{\hookrightarrow} D$ be a trivial cofibration in $\Ccal$. Let $X^+$ be the pushout $\displaystyle X \coprod_C D$.
We will first prove that the natural map $X \rightarrow X^+$ is a weak equivalence.

Consider a square:

\[
\begin{tikzcd}[ampersand replacement=\&]
A \arrow[hookrightarrow]{d}{i} \arrow{r}  \&  X \arrow{d} \\
B \arrow{r} \& X^+ \\
\end{tikzcd}
\]

Where $i: A \hookrightarrow B$ is a cofibration in $\Ccal$.

By assumption $X$ is retract of a cellular object $X'$ which is a directed co-limit of representable objects $X_i$ (with cofibrant transition map), the map from $C \rightarrow X$ lifts to a map $C \rightarrow X'$ and one defines $X'^+ = X' \coprod_C D$, and $X^+$ is easily seen to be a retract of $X'^+$, in particular, one can lift all the above square to $X'$ and $X'^+$

$X_0^+$ can then be represented by a colimits of $X_i \coprod_C D$ restricted to the $X_i$ with a map from $C$ compatible to the map to $X_0$, and hence using again that directed colimit are computed objectwise, there exists a representable $X_i$ with a map $X_i \rightarrow X$ such that the above square is lifted to:

\[
\begin{tikzcd}[ampersand replacement=\&]
A \arrow[hookrightarrow]{d}{i} \arrow{r}  \& X_i \arrow{d} \arrow{r} \& X' \arrow{d} \arrow{r} \&  X \arrow{d} \\
B \arrow{r} \& X_i\coprod_C D  \arrow{r}\& X'^+ \arrow{r} \& X^+ \\
\end{tikzcd}
\]

But as the map $X_i \rightarrow X_i \coprod_C D$ is a trivial cofibration in $\Ccal$ one has the lifting we are looking for directly in the first square by proposition \ref{Prop_Semilifting2}.

Moreover, because of the characterization of weak equivalence in terms of the $\pi$-sets (\ref{PiCharacOfEquiv}) and the fact that the $\pi$-sets clearly commutes to directed colimits one directly see that a transfinite composition of weak equivalences is a weak equivalences and that a retract of a weak equivalence is a weak equivalence. Hence this show that every trivial cofibration of cofibrant domain is a weak equivalences.

Conversely, let $f: X \rightarrow Y$ be a cofibration with a cofibrant domain which is a weak equivalence. One can factor $f$ as a trivial cofibration followed by a fibration, but the trivial cofibration part will be a weak equivalences by the above result and hence by the $2$-out-of-$3$ property the fibration part will be a weak equivalence, hence a trivial fibration. Using the left lifting property of $f$ (a cofibration) against this trivial fibration exhibit $f$ as a retract of a trivial cofibration which concludes the proof.
}
}

\blockn{At this point, we have proved that:}

\block{\Th{If $\Ccal$ is a cylinder category, $\widetilde{\Ccal}$ is a combinatorial semi-model category in which every object is fibrant. The generating (trivial) cofibrations are the (trivial) cofibrations in $\Ccal$. }
}

\blockn{It has moreover a very simple universal property:}

\block{\label{CompletionMorphism}\Prop{Let $\Vcal$ be a semi-model category in which every object is fibrant. Then the category $\Vcal^{cof}$ of cofibrant objects of $\Vcal$ is a cylinder category. If $\Ccal$ is a small cylinder category then any morphism $F$ from $\Ccal$ to $\Vcal^{cof}$ extend to a Quillen adjunction:

\[ \tilde{F}: \tilde{\Ccal} \rightarrow \Vcal : R_F. \]

Any left Quillen functor $\tilde{\Ccal} \rightarrow \Vcal$ is obtained this way.
}

\Dem{The fact that one gets an adjunction which extend $F$, and the fact that any such functor is obtained this way is basic Kan extension theory (it only require $\Vcal$ to be co-complete). Moreover, any such $\tilde{F}$ send the generating (trivial) cofibrations in $\tilde{\Ccal}$ to (trivial) cofibrations in $\Vcal$ by definition hence is a left Quillen functor.}
}

\block{\label{CompletionEquivalence}\Prop{Let $F: \Ccal \rightarrow \Dcal$ be a morphism of cylinder categories Then:

\begin{enumerate}

\item For every object $Y$ of $\tilde{\Dcal}$, every cofibration $A \hookrightarrow B$ in $\Ccal$ and every $y:F(A) \rightarrow Y$ the natural map:

\[ \pi_{B/A}(R_F(Y),y) \rightarrow \pi_{F(B)/F(A)} (Y,y) \]

induced by the adjunction is a bijection.

\item If $F$ is homotopy fully faithful then, for every cofibrant object $X$ of $\tilde{\Ccal}$, for every $A \hookrightarrow B$ a cofibration in $\Ccal$, and for every $x : A \rightarrow X$, the natural map:

\[ \pi_{B/A}(X,x) \rightarrow \pi_{F(B)/F(A)}(\tilde{F}(X),\tilde{F}(x) ) \]

induced by applying $\tilde{F}$ is a bijection.

\item If $F$ is acyclic, then the Quillen pair $(\tilde{F} ,R_F)$ is a Quillen equivalence.

\end{enumerate}

}

\Dem{ Note that if $I_A B$ is a relative cylinder object for $A \hookrightarrow B$ in $\Ccal$ then $F(I_A B)$ is a relative cylinder object for $F(A) \hookrightarrow F(B)$ in $\Dcal$.

\begin{enumerate}

\item A morphism $B \rightarrow R_F(Y)$ is the same as morphism from $F(B)$ to $Y$, and an homotopy $I_A B \rightarrow R_F (Y)$ is the same as a homotopy $F(I_A B) \rightarrow Y$ so the first point follows immediately.

\item The map is clearly well defined: any homotopy $I_A B \rightarrow X$ will produces a homotopy $F(I_A B) \rightarrow \tilde{F}(X)$ by applying $\tilde{F}$.

Fix $x:A \rightarrow X$ and let $v :F(B) \rightarrow \tilde{F}(X)$ be any morphism in $\tilde{\Dcal}$ which represents an element of $\pi_{F(B)/F(A)}(\tilde{F}(X),\tilde{F}(x))$ . Being cofibrant, $X$ can be written as (a retract of) a directed colimit of representable, and $\tilde{F}$ commutes to colimit, so there exists a representable object $X_i \rightarrow X$ such that $v$ can be factored as: $F(B) \overset{v'}{\rightarrow} F(X_i) \rightarrow \tilde{F}(X)$.

Applying the fact that $F$ is homotopy fully faithful, one can find a replacement $v'':B \rightarrow X_i$ such that $F(v'')$ will be homotopically equivalent to $v'$ relative to $F(A)$ and composing with the map $X_i \rightarrow X$ produces a map $B \rightarrow X$ in $\pi_{B/A}(X,x)$ whose image in $\pi_{F(B)/F(A)}(\tilde{F}(X),\tilde{F}(x))$ is $v$. This proves the surjectivity of this comparison map.

The exact same argument applied to $B \coprod_A B \hookrightarrow I_A B$ instead of $A \hookrightarrow B$ proves the injectivity.

\item Let $X$ be a cofibrant object in $\tilde{\Ccal}$ and $Y$ be any object in $\Dcal$, let $v: X \rightarrow R_F(Y)$ with adjoint transpose: $w: \tilde{F}(X) \rightarrow Y$. We need to prove that $v$ is a weak equivalence if and only if $w$ is. From the previous two points it is easy to see that $v$ is a weak equivalence if and only if for every $A \hookrightarrow B$, for every $x:A \rightarrow X$ the application:

\[ \pi_{F(B)/F(A)} (w, \tilde{F}(x)) : \pi_{F(B)/F(A)} (\tilde{F}(X), \tilde{F}(x)) \rightarrow \pi_{F(B)/F(A)} (Y, w \circ \tilde{F}(x)) \]

is a bijection. We need to show that this is enough to conclude that $w$ is a weak equivalence.

First applying the second point to $0 \hookrightarrow A$ one can see that any arrow $F(A) \rightarrow \tilde{F}(X)$ is homotopy equivalent to an arrow of the form $\tilde{F}(x)$ so by \ref{PiLeftInvariance} this shows that the above map is again a bijection when $\tilde{F}(x)$ is replaced with any map from $F(A)$ to $\tilde{F}(X)$, and finally, using that $F$ is acyclic, one can construct for any cofibration $i:A' \hookrightarrow B'$ in $\Dcal$ a square of the form:

\[
\begin{tikzcd}[ampersand replacement=\&]
A' \arrow[hookrightarrow]{d}{i'} \arrow{r}{\sim}  \&  F(A) \arrow[hookrightarrow]{d}{F(i)} \\
B' \arrow{r}{\sim} \& F(B) \\
\end{tikzcd}
\]

where the horizontal arrows are weak equivalences and using \ref{PiLeftTransport}, this concludes the proof.

\end{enumerate}
}

}

\block{Let us mention an example: If $\Ccal$ is the category of finite $CW$-complexes then $\widetilde{\Ccal}$ can be difficult to describe explicitly. But the full  subcategory of cellular objects is equivalent to the category of $CW$-complexes: indeed, this follows from lemma \ref{LemCellularObjectDirect} and the fact that a map from a finite $CW$-complex to an arbitrary $CW$-complex always factor into a finite sub-$CW$-complex by a compactness argument. As the general cofibrant objects are the retract of cellular objects, the full subcategory of cofibrant objects is equivalent to the category of retract of $CW$-complexes, i.e. of topological spaces which are cofibrant from the model structure on the category of spaces.
Moreover it is immediate that the weak equivalences, cofibrations and fibrations between them are the same as in the category of spaces.

It follows from that that the Quillen adjunction $\widetilde{\Ccal} \rightleftarrows Spaces$ induced by the inclusion of $\Ccal$ in $Spaces$ is a Quillen equivalence: it induces a categorical equivalence between the categories of fibrant-cofibrant objects which preserves weak equivalences, fibrations and cofibrations.

In particular, one has that:
}

\block{\label{EquivWithFinCW}\Prop{If $\Ccal$ is a cylinder category with an acyclic morphism to the cylinder category of finite $CW$-complexes, then the induced adjunction between $\widetilde{\Ccal}$ and the model category of spaces is a Quillen equivalence.}}

\blockn{Finally, in the case where the category $\Ccal$ comes with a class of generating cofibrations, one can use it to obtains smaller set of generating cofibrations and trivial cofibrations in $\widetilde{\Ccal}$.}

\block{\label{PropModelStrFromGenCof}\Prop{Let $\Ccal$ be a cylinder category with a set $A \hookrightarrow B$ of generating cofibrations. Then for a map $f:X \rightarrow Y$ in $\widetilde{\Ccal}$:

\begin{itemize}

\item $f$ is a trivial fibration if and only it has the right lifting property with respect to the generating cofibrations.

\item $f$ is a weak equivalence if and only if all $A \hookrightarrow B$ a generating cofibration, the maps $\pi_{B/A}(f,x):\pi_{B/A}(X,x) \rightarrow \pi_{B/A}(Y,f\circ x)$ is surjective, i.e. if all square of the form:

\[
\begin{tikzcd}[ampersand replacement=\&]
A \arrow[hookrightarrow]{d}{i} \arrow{r}  \&  X \arrow{d} \\
B \arrow{r} \& Y \\
\end{tikzcd}
\]

where $A \hookrightarrow B$ is a generating cofibration can be filled into:

\[
\begin{tikzcd}[ampersand replacement=\&]
A \arrow[hookrightarrow]{ddd}{i} \arrow[hookrightarrow]{dr}{i} \arrow{rr} \& \&  X \arrow{ddd} \\
\& B \arrow[hookrightarrow]{d} \arrow{ur}{a} \& \\
\& I_A B \arrow{dr}{h} \& \\
B \arrow{rr} \arrow[hookrightarrow]{ur}\& \& Y \\
\end{tikzcd}
\]

\item $f$ is a fibration if and only if it has the right lifting property with respect to at least one trivial cofibration of the form $B \hookrightarrow I_A B$ for each generating cofibration $A \hookrightarrow B$.

\end{itemize}
}

\Dem{

\begin{itemize}

\item The first claim is relatively immediate: Any cofibration in $\Ccal$ is a composite of pushout of the generating cofibrations, so something that has the right lifting property with respect to generating cofibrations has the right lifting property with respect to all cofibrations.

\item First take a cofibrant replacement $\widetilde{X}$ of the domain of $X$ the map $\widetilde{X} \rightarrow X$ satisfies the same condition as $f$ because the map $\widetilde{X} \rightarrow X$ preserves all $\pi$-sets. Hence one can freely assume that $X$ is cofibrant. Then consider a trivial cofibration/fibration factorization of $f$, as the domain of $f$ is cofibrant, the trivial cofibration part of the factorization is an equivalence hence by the exact same argument as in the proof of \ref{CompWE2out3}, the fibration part satisfies the same condition as $f$. But this is enough to make the proof of \ref{Prop_trivFibInComp} shows that this fibration has the right lifting property with respect to all generating cofibrations hence is a trivial fibration, hence a weak equivalence and this concludes the proof. 

\item Let $f$ be a trivial cofibration with cofibrant domain. Using the small object argument, one can factor $f$ as $p \circ i$ where $i$ is in the class generated by maps of the form $B \hookrightarrow I_A B$ and $p$ has the right lifting property against those maps. Both $f$ and $i$ are weak equivalences, hence $p$ is also a weak equivalence. Under those conditions, the proof of \ref{Prop_trivFibInComp} shows that $p$ is a trivial fibration. By the usual arguement on weak factorization systems, this suffices to prove that $f$ is a retract of $i$ and hence that the lifting property of fibrations can be tested against only the maps that we have used in the construction of $i$, i.e. o the form $B \hookrightarrow I_A B$.

\end{itemize}
}

One can also gives a longer proof (but more informative) of the last point by proving that any trivial cofibration $C \hookrightarrow D$ in $\Ccal$ is a retract of $D \hookrightarrow I_C D$ and if a cofibration $C \hookrightarrow D$ is obtained as $n$ pushout of generatng cofibration $A_i \hookrightarrow B_i$, the trivial cofibration $D \hookrightarrow I_C D$ can be obtained as $n$ pushout of of the trivial cofibrations $B_i \hookrightarrow I_{A_i} B_i$.
}

\section{Pre-cylinder categories and their algebraic properties}
\label{SecPreCylinder}

\blockn{We will now define and studies pre-cylinder categories, this notion corresponds roughly to the algebraic structure underlying the notion of cylinder categories, i.e. the part of the structure that is preserve by morphisms and subject to free constructions. In the next section we will construct a sort of weak model structure on the category of pre-cylinder categories such that the cylinder categories are exactly the fibrant objects.}

\subsection{Definition and generalities}

\block{\Def{A Pre-cylinder category $\Ccal$ is a cofibration category (definition \ref{DefCofibCat}) with an additional class of maps $\Wcal$ called weak equivalences for which the $2$-out-of-$6$ property and the ``cube lemma'', i.e. proposition \ref{gluingLemma} holds. }

Of course we will call ``trivial cofibrations'' the maps that both cofibrations and weak equivalences.
}

\blockn{The idea behind this definition is to isolate the ``algebraic part'' in the notion of cylinder category. In fact, up to some problem related to the fact that those are not ordinary category but a $2$-category, the category of pre-cylinder categories and the category of cofibrations categories would be category of models of a quasi-equational theory or partial Horn theory as in \cite{palmgren2007partial} (for a correct $2$-categorical version of this) and hence should have all limits, all co-limits and all the free constructions that we might want. Those free constructions will be discussed more precisely in the next subsection \ref{SectionFreeCofCat}. }

\blockp{Note that the ``cube lemma'' immediately implies some properties the we have for cylinder categories:}

\block{\Prop{In a pre-cylinder category, a pushout of a trivial cofibration is a trivial cofibration, a pushout of a weak equivalence along a cofibration is a weak equivalence.}

\Dem{If $A \hookrightarrow B$ is a trivial cofibration a map $C \rightarrow C \coprod_{A} B$ can be considered as a comparison map $C \coprod_A A \rightarrow C \coprod_{A} B$ and hence is a weak equivalence by the cube lemma. Similarly, if $A \hookrightarrow B$ is a cofibration and $A \rightarrow C$ is a weak equivalence, then the map $B \rightarrow B \coprod_A C$ is the comparison map $ B \coprod_A A \rightarrow B \coprod_A C$ and hence is a weak equivalence by the cube lemma.}

}

\block{\label{adjointsOfForget}\Prop{The forgetful functor from pre-cylinder categories to cofibrations categories have both a left and a right adjoint.

They are given by taking the same underlying cofibration category and deciding respectively that only the isomorphisms are weak equivalences or that all maps are weak equivalences.
}
The proposition is rather immediate to check: those trivially satisfies the cube axiom and the $2$-out-of-$6$ property, and they do have the correct universal property to define adjoint functors.

This proposition has a very interesting consequences: Computing limits, co-limit or free construction in the category of pre-cylinder category or in the category of cofibrations category would be essentially the same: the forgetful functor will preserve all sort of limits and co-limits, so when one wants to compute this kind of categorical constructions in the category of pre-cylinder category one can do it by first computing it in the category of cofibration categories and then constructing the class of weak equivalences as the smallest or largest class (depending if one perform an inductive or projective construction) satisfying some conditions.

}

\blockn{We will now give basic example of limits, co-limits and free construction before moving in the next subsection to a more general theory of free constructions (including colimits) in the category of pre-cylinder category.}

\block{\Prop{The product of two pre-cylinder categories (or of two cofibration categories) exists and is their categorical product with the cofibrations (resp. the weak equivalences) being the maps whose two components are cofibrations (resp. weak equivalences). 

This extend to infinite products.
}

\Dem{It is easy to see that it is indeed a pre-cylinder category (or a cofibration category) and that it indeed satisfies the universal property.}
}

\block{\Lem{The functor $0$ is always a morphism between any two pre-cylinder categories (or cofibration categories). If $F,G:\Ccal \rightrightarrows \Dcal$ are two morphisms between pre-cyinder categories (or cofibrations category) then:

 \[ X \mapsto F(X) \coprod G(X) \]

is also a morphism.}

Those two claims are immediate to check (the fact that the co-product preserve weak equivalences is a special case of the cube lemma).
}

\block{\Prop{The co-product of two pre-cylinder categories (or cofibration categories) $\Ccal$ and $\Dcal$ is isomorphic to their product $\Ccal \times \Dcal$. The co-projections are the maps $X \mapsto (X,0)$ and $X\mapsto (0,X)$. If $f:\Ccal \rightarrow \Xcal$ and $g:\Dcal \rightarrow \Xcal$ are two morphisms, then the corresponding morphism $(f,g):\Dcal \times \Ccal \rightarrow \Xcal$ is:

\[(f,g)(X,Y) = f(X) \coprod g(Y) \]

This extend to finite co-product.
}

One could also form infinite co-products: it will be the full subcategory of the infinite product of objects which are initial on all component except a finite number, but this would follow from the more general theory developed in the next subsection.

\Dem{It is not very hard to check directly. It also follows from the lemma above by the exact same categorical argument as the fact that the co-product (if it exists) is always isomorphic to the product in an additive category (with the addition being replaced by the co-product).}
}

\block{\label{PropFstar}
\Prop{ The free pre-cylinder category $F_*$ on one object $*$ is (equivalent to) the category of finite sets, with the monomorphisms as cofibrations and the isomorphisms as weak equivalences. The object $*$ being the singleton.
}

\Dem{It is easy to check that it is indeed a cylinder category. Moreover if $\Ccal$ is any pre-cylinder category and $X$ is an object of $\Ccal$, then one can form for $S$ a finite set the iterated co-product of $X$:

\[ X^{(S)} = \coprod_{i \in S} X \]

This defines a functor from finite sets to $\Ccal$, isomorphisms are sent to isomorphisms hence weak equivalences, injections are sent to the structural maps of a coproduct: $X^{(S_1)} \rightarrow X^{(S_1) \coprod (S_2)} = X^{(S_1)} \coprod X^{(S_2)}$ hence a cofibration, and pushout along cofibration are preserved by this functor. Hence it is a morphism from $F_*$ to $\Ccal$ which send the singleton to $X$ and such a morphism is obviously unique (up to unique isomorphism)
}
}

\subsection{Completion and computation of free pre-cylinder categories}
\label{SectionFreeCofCat}

\blockn{The goal of this sub-section is to develop a general procedure to compute all co-limits and free constructions within the category of pre-cylnder categories and the category of cofibrations categories.

\bigskip

As observed in \ref{adjointsOfForget}, the important part is to understand free constructions and colimits of cofibration categories. The case of pre-cylinder categories is then easily deduced from that by first making the construction at the level of cofibration categories and then just taking the smallest class of weak equivalences that makes the cofibration category obtained a pre-cylinder category and the structural morphisms, morphisms of pre-cylinder category.}

\blockn{The first step would be to precisely introduce the type of free constructions that one want to consider.}

\block{\Def{A system of generators for a cofibrations category is the data of:

\begin{enumerate}

\item A set of cofibrations categories $(\Ccal_i)_{i \in I}$.

\item A set $J$ of indices for objects $(X_j)_{j \in J}$.

\item A set $A$ of indices for arrows, to each indices for arrow is attached a ``domain'' and a ``co-domain'' which are formal co-product of a finite family of objects of the $\Ccal_i$ and/or of some of the $X_j$.

\item A subset $C$ of the arrows are ``cofibrations''.

\item A set of relations between the arrows.

\end{enumerate}

A system of generators for a pre-cylinder category is the same except that: $1.$ and $4.$ are replaced by:

\begin{itemize}
\item[1'.] A set of pre-cylinder categories $(\Ccal_i)_{i \in I}$.
\item[4'.] A subset $C$ of the arrows are ``cofibrations'' and a subset $W$ of the arrows are weak equivalences.
\end{itemize}

}

The points $4.$ ($4'.$) and $5.$ need to be clarified: we mean by a subset of arrows ($W$ or $C$) a subset of the set of all arrow that can be constructed out the arrows of the category $\Ccal_i$, the objects $X_j$ for $j\in\Jcal$ and the arrows in $A$  between finite coproducts of these objects. Similarly a relation between arrows is an equality between two such parallel arrows.

\bigskip

A reader which find this last paragraph unclear or not formal enough can follow a more precise approach: one first consider systems of generators with only the axioms $1,1',2,3$ and one defines the cofibration categories/pre-cylinder categories freely generated by such a system of generators (in the sense of the next definition). 

One can then return to the more general definition of system of generators: the point $4.$ asks for a subset of the sets of arrows of the category freely generated by the previous structure (points $1,1',2,3$), the set $W$ is the same things, and the point $5$ is a set of pairs of parallel arrows in this freely generated category. 

\bigskip

Moreover, one can iterate this procedures in order to also ask some properties or relations between arrows that are constructed using pushout along maps that are cofibrations because they are forced to be by the the choice made in point $4.$. But this does not need to be done at this level, it will be automatically covered by iterating free construction and everything that is discussed below will also apply to this situation.

}

\block{\Def{ 
A representation of a system of generator for a cofibration category in a cofibration category $\Dcal$, is the data of:

\begin{enumerate}
\item a morphism $\Ccal_i \rightarrow \Dcal$ of each $i \in I$,
\item an object $X_j \in \Dcal$ for each $j \in J$, 
\item for each $a \in A$ an arrow in $\Dcal$ between the corresponding co-product,
\item the arrow in $C$ are cofibrations,
\item the relations are satisfied.
\end{enumerate}

A morphism between two representations is the data of natural transformations between the morphisms specified in point $1.$ and of morphisms between the object specified in point $2.$, such that it satisfies the naturality conditions with respect to all the arrows specified in points $3.$.

\bigskip

A cofibrations category freely generated by a system of generators is a cofibrations category $\Ccal$ such that for any other cofibrations category $\Dcal$ the category of morphisms from $\Ccal$ to $\Dcal$ is equivalent, functorially in $\Dcal$, to the category of representations in $\Dcal$ of the system of generator.

\bigskip

One defines the same notions for pre-cylinder category, the only change is that in $1.$ one asks for morphisms of pre-cylinder categories and in $4.$ on also asks that the maps in $W$ are weak equivalences in $\Dcal$.

}}

\block{\Prop{For every system of generators, there is a essentially unique small cofibrations category/ pre-cylinder category generated by it.}

\Dem{The uniqueness property is clear. The existence would basically follow from a $2$-categorical version of theorem $29$ of \cite{palmgren2007partial}, as we do not know any reference developing this sort of result into a $2$-categorical context we will use some trick to deduce it from the $1$-categorical version. The proof for cofibrations categories and cylinder categories are exactly the same, we will only write it for cofibration categories.

\bigskip

Let us temporarily call a strict cofibration category to be a (small) cofibration category with a specified initial object and in which for every diagram of the form $C \leftarrow A \hookrightarrow B$ one has chosen a specific pushout $C \coprod_A B$. A morphism between strict cofibration category is called strict if it preserves the marked initial object and the chosen pushout.

\bigskip

The ($1$-)category of strict cofibration categories and strict morphism is the category of models of quasi-equational theory (in the temrinology of \cite{palmgren2007partial}), which would have one sort for objects, one sort for morphisms, one predicate for ``cofibrations''. Hence, because of theorem $29$ of \cite{palmgren2007partial} one can have models freely generated by a system of generator in the sense that one will have a strict cofibration category $\Ccal$ such that strict morphisms from $\Ccal$ to $\Dcal$ are the same as representation of the system of generators in $\Dcal$.

We will now show that this $\Ccal$ satisfies the correct universal property. The first observation is that for any (strict) cofibration category $\Dcal$, there is a (strict) cofibration category $\Dcal^{mor}$ whose objects are pair of objects of $\Dcal$ together with a morphism between them (cofibrations are the object wise cofibrations) and a morphism from $\Ccal$ to $\Dcal^{mor}$ is the same as a pair of morphisms from $\Ccal$ to $\Dcal$ together with a natural transformation between them. Hence morphism between two (strict) functors from $\Ccal$ to $\Dcal$ can be described as (strict) morphisms from $\Ccal$ to $\Dcal^{mor}$ hence as representations of the system of generator in $\Dcal^{mor}$ and those are exactly the morphisms of representations as defined above. Hence the category of strict morphism from $\Ccal$ to $\Dcal$ and natural transformation between them is indeed what we want. Any morphism (possibly non strict) form $\Ccal$ to $\Dcal$ can be made strict by eventually replacing $\Dcal$ by a equivalent cofibration category, so this concludes the proof.
}
}

\blockn{The argument above is extremely inexplicit, and the rest of this subsection is devoted to present a more efficient construction to produce those free categories. As mentioned above it is enough to deal with the case of cofibration categories, one can then deal with weak equivalences separately. The key observation is the following: }

\blockn{\Prop{Let $\setop$ be the opposite of the category of sets, with every map being a cofibration. Then $\setop$ is a cofibration category and for every cofibration category $\Ccal$, there is a natural equivalence:

\[ \text{Hom}(\Ccal,\setop) \simeq \widetilde{\Ccal}^{op} \]
 }

\Dem{This is relatively immediate: $\setop$ has all colimits because the category of sets has all limits, a morphism from $\Ccal$ to $\setop$ is hence exactly a contravariant functor from $\Ccal$ to sets which takes the initial object to the singleton and a pushout along cofibration to pullback of sets, hence an object of $\widetilde{\Ccal}$. The natural transformations between such morphisms corresponds to morphism of pre-sheaf in the opposite direction. } 
 
}

\block{\label{FreePrecylinderCatMain}This proposition allow to very easily obtain a simple description of the completion $\widetilde{\Dcal}$ of the cofibration category $\Dcal$ freely generated by a system of generators: it corresponds to a representation of the system of generator in $\setop$. More precesely:

\begin{enumerate}
\item For each $i \in I$, an object $D_i$ of $\widetilde{\Ccal_{i}}$.
\item for each $j \in J$ a set $D_j$.
\item For each arrow in $A$, an arrow in the other direction between the corresponding finite product.
\item Asking for some arrows to be a cofibrations (or a weak equivalence) does not have any influence at this level as every map in $\setop$ is a cofibration.
\item All the relations between the arrows should be satisfied.
\end{enumerate}

Morphisms being the obvious notion.

Next, one needs to compute what are the functors $\Ccal_i \rightarrow \widetilde{\Dcal}$ and what are the objects of $\widetilde{\Dcal}$ corresponding to the $X_j$. There is no formula for this, but they satisfy a universal property:

\begin{itemize}
\item For any $i \in I$ and for any $c \in \Ccal_i$ then there is an object $f_i(c) \in \Dcal$ such that functorially in $D \in \Dcal$ one has:

\[ \text{Hom}(f_i(c),D) \simeq D_i(c) \]

\item For any $j \in J$ there is an object $X_j$ in $\Dcal$ such that functorially in $D \in \Dcal$:

\[ \text{Hom}( X_j,D) \simeq D_j \]

\end{itemize}

Such objects always exist but there is in general no easy procedure to construct them and one need to find those objects ``by hand'' in each specific case. Although it appears that it is often easy and we will see some of examples in the rest of the article. 

Using the universal property of these objects, one easily extend this into a functor $f_i : \Ccal_i \rightarrow \Dcal$ which send coproducts along cofibrations to coproducts, and all for all arrows in $A$ one can construct a corresponding arrow between coproducts of the $f_i(c)$ and the $X_j$.

One can then say that maps in $\widetilde{\Dcal}$ is a \emph{generating cofibration}, if it is: either coming from a (generating) cofibration in $\Ccal_i$ though the construction above, or from an arrow in $C$ from the construction above, or it is the arrow $0 \rightarrow X_j$ for a $j \in J$.

One then has:

\Prop{\begin{itemize}

\item $\Dcal$ is the full sub-category of $\widetilde{\Dcal}$ of objects which can be obtained from $0$ as an iterated pushout of generating cofibrations.

\item Cofibrations in $\Dcal$ are the arrows that are iterated pushout of the generating cofibrations (i.e. the generating cofibrations are indeed generating cofibrations of $\Dcal$).
 
\end{itemize} 
 }

\Dem{Let us first observe that if we define $\Dcal'$ to be the full subcategory of object which can be obtained as iterated pushout of generating cofibration and the cofibrations between objects of $\Dcal'$ being the iterated pushout of generating cofibration, then $\Dcal'$ is a cofibraton category: isomorphisms are cofibration, and cofibrations are stable under composition, pushout along cofibrations exists (and are computed in $\widetilde{\Dcal}$).

Moreover, any cofibrations in $\Dcal'$ is also a cofibration in $\widetilde{\Dcal}$ (in the sense of definition \ref{CofInTildeC}).

$\Dcal'$ contains all the structure that generated $\Dcal$: the image of the $f_i$ the object $f_j$ and all the maps defined between them, moreover the maps $f_i : \Ccal_i \rightarrow \Dcal'$ are morphisms of cofibration categories and the maps in the set $C$ are cofibration in $\Ccal$ (because they are generating cofibration). Hence there is a morphism from $\Dcal$ to $\Dcal'$ preserving all these structure.

But if one look at the composition of this morphism with the inclusion $\Dcal' \rightarrow \widetilde{\Dcal}$ we obtain the Yoneda embeddings: indeed all the generators of $\Dcal$ are send to their Yoneda embeddings and this is a morphism from $\Dcal$ to the cofibration category of cofibrant object of $\widetilde{\Dcal}$ hence it is the Yoneda embeddings by the uniqueness in the universal property of $\Dcal$. This proves that $\Dcal$ and $\Dcal'$ are equivalent as category. Moreover, Their cofibrations are the same: any map which is a cofibration in $\Dcal$ is also a cofibration in$\Dcal'$ because the functor $\Dcal \rightarrow \Dcal'$ is a morphism of cofibration category, and any morphism which is a cofibration in $\Dcal'$ is a cofibration in $\Dcal$ because it is an iterated pushout of generating cofibrations and generating cofibration are cofibration in $\Dcal$.

}

}

\blockn{We finish with a simple example of the construction above:}

\block{\label{DefFhook}\Prop{The pre-cylinder category freely generated by a cofibration $A \hookrightarrow B$ between two objects is the category $F_{\hookrightarrow}$ such that:

\begin{itemize}
\item objects are pair of finites sets $X_A,X_B$ with a map:

\[ X_c: X_B \rightarrow X_A \]

\item Morphisms are morphisms of diagrams.

\item Cofibrations are objectwise injections.

\item Weak equivalences are isomorphisms.

\end{itemize}

Moreover, if $A \hookrightarrow B$ is a cofibration in a pre-cylinder category $\Ccal$, the corresponding morphism from $F_{\hookrightarrow}$ to $\Ccal$ send $X_B \rightarrow X_A$ to the object constructed by first taking the co-product of $X_A$ copies of $A$ and then gluing $X_B$ copies of $B$ along the copies of $A$ corresponding to the image in $X_A$ by $X_c$.
}

\Dem{The completion of $F_{\hookrightarrow}$ is the category of pair of sets with a map between them and morphisms are morphisms of diagram. Checking the universal property, the object $A$ is $ \emptyset \rightarrow \{ *\}$ and $B$ is $\{*\} \rightarrow \{*\}$ and the generating cofibrations are\footnote{As well as the map from the initial object to $B$, but it is the composite of the two others.} the map from the initial object to $A$ and the map from $A$ to $B$. Hence the finite cellular maps are the map that are obtained by iteratively: freely adding an element to $X_A$ or freely adding an element of $X_B$ whose image in $X_A$ is specified, i.e. they are exactly the object wise ``co-finite'' monomorphisms, it follows that the object of $F_{\hookrightarrow}$ are just the pair of finite set with a map, cofibrations are the objectwise monomorphism, and as nothing in the system of generators impose that certain arrow are weak equivalences, the weak equivalences are just is isomorphisms.

The last claim just translate that the morphism will preserve co-product along cofibrations.  }

}

\subsection{The monoidal closed structure on the category of pre-cylinder categories}
\label{SubSecmonoidalStr}

\blockn{Both the category of pre-cylinder categories and of categories with cofibrations have a very natural structure of symmetric monoidal closed category. We will explain it only of pre-cylinder categories, but everything carries over easly to the case of cofibrations categories by just forgetting everything concerning the weak equivalences. We will start by describing the ``Hom object''.}

\block{\Prop{Let $\Ccal$ and $\Dcal$ be pre-cylinder category, there is a pre-cylinder category $[\Ccal,\Dcal]$ such that:

\begin{itemize}

\item Objects are the morphisms from $\Ccal$ to $\Dcal$.

\item Morphisms are the natural transformations.

\item Weak equivalences are the objectwise weak equivalences.

\item Cofibrations are the natural transformations $\lambda:F \rightarrow G$ such that for any cofibration $c:A \hookrightarrow B$ in $\Ccal$, the arrow:

\[ G(A) \coprod_{F(A)} F(B) \rightarrow G(B) \]

induced by the naturality square (for $\lambda$ and $c$) is a cofibration.

\item Pushout along cofibrations are computed objectwise.

\end{itemize}

}
\Dem{Remark first that cofibrations are in particular objectwise cofibrations : this can be seen by considering the case of the cofibration $0 \hookrightarrow A$ in $\Ccal$, which gives that $F(A) \hookrightarrow G(A)$ is a cofibration for all $A$.

Let $i:A \hookrightarrow B$ and $f:A \rightarrow C$ be a cofibration and an arrow in $[\Ccal,\Dcal]$, as $i$ is in particular an objectwise cofibration one can form the objectwise pushout out $D= B \coprod_A C$ which is a functor from $\Ccal$ to $\Dcal$. We will first prove that it is a morphism of pre-cylinder categories:

\begin{itemize}

\item $D$ preserves cofibrations: let $i: X \hookrightarrow Y$ be a cofibration in $\Ccal$. Then:

\[ C(X) \coprod_{A(X)} B(X) \rightarrow C(Y) \coprod_{A(X)} B(X) \]

Is a cofibration because $C(X) \rightarrow C(Y)$ is a cofibration.

The square:

\[
\begin{tikzcd}[ampersand replacement=\&]
\displaystyle A(Y) \coprod_{A(X)} B(X) \arrow{d} \arrow{r} \& \displaystyle C(Y) \coprod_{A(X)} B(X) \arrow{d} \\
B(Y) \arrow{r} \& \displaystyle C(Y) \coprod_{A(Y)} B(Y) \\
\end{tikzcd}
\]

Is a pushout (this is easily observed by examining the universal properties) and the left map is a cofibration by assumption (because $i$ is a cofibration) hence the right map is also a cofibration and this shows that $D(i)$ is a cofibration.

\item $D$ preserves pushout along cofibrations because colimits commute to colimits (or by examining the universal property of all the pushout involved).

\item $D$ preserves weak equivalences exactly because of the cube axiom.

\end{itemize}

At this point one has proved that $[\Ccal,\Dcal]$ admit pushout along cofibrations and that they are computed objectwise, and almost all the axioms of pre-cylinder categories can hence now be checked objectwise in $[\Ccal,\Dcal]$ and are all valid for $[\Ccal,\Dcal]$.

The only non trivial axiom left is that if $A \overset{u}{\hookrightarrow} B \overset{v}{\hookrightarrow} C$ are two cofibrations in $[\Ccal,\Dcal]$ then $v \circ u$ is a cofibration. Let $i : X \hookrightarrow Y$ be a cofibration in $\Ccal$.

By assumption, one has a cofibration:

\[A(Y) \coprod_{A(X)} B(X) \hookrightarrow B(Y) \]

One can push it forward to a cofibration:

\[A(Y) \coprod_{A(X)} C(X) \hookrightarrow B(Y) \coprod_{B(X)} C(X) \]

and then compose it to the cofibration:
\[B(Y) \coprod_{B(X)} C(X) \rightarrow C(Y) \]
and this concludes the proof.

}
}

\blockn{Even when $\Dcal$ is a cylinder category, it is relatively rare that $[\Ccal,\Dcal]$ is a cylinder category. This makes this hom object a little more difficult to use in a correct homotopy theoretic fashion and will not study in the present paper the precise relationship between the monoidal structure and the weak model structure on the category of pre-cylinder category, but it is not hard too see that it is \emph{not} a monoidal model structure. This being said, conditions on $\Ccal$ under which $[\Ccal,\Dcal]$ is a cylinder category when $\Dcal$ is will be studied in subsection \ref{SubSecReedy}.}

\block{\label{CritForCofibinHom}\Prop{A morphism $\lambda: F \rightarrow G$ in $[\Ccal,\Dcal]$ is a cofibration if and only if for a class of \emph{generating} cofibrations $A_i \hookrightarrow B_i$ in $\Ccal$ the map:

\[ G(A_i) \coprod_{F(A_i)} F(B_i) \rightarrow G(B_i) \]

induced by the naturality square is a cofibration.}

\Dem{We will prove by induction on the number of gluing generating cofibrations that under the assumption of the theorem if $i:C \hookrightarrow D$ is a cofibration then the map:

\[ G(C) \coprod_{F(C)} F(D) \rightarrow G(D) \]
is a cofibration.

It holds if $i$ is an isomorphism. Assume it hold for a cofibration $i:C \hookrightarrow D$ then we have to prove that it also holds for $C \hookrightarrow D^+ = D \coprod_{A_i} B_i$ for some map $A_i \rightarrow D$.

\[ G(C) \coprod_{F(C)} F(D^+) \simeq G(C) \coprod_{F(C)} F(D) \coprod_{F(A_i)} F(B_i) \]
\[G(D^+) \simeq G(D) \coprod_{G(A_i)} G(B_i) \]

Hence we have a cofibration:

\[G(C) \coprod_{F(C)} F(D) \coprod_{F(A_i)} F(B_i) \hookrightarrow G(D) \coprod_{F(A_i)} F(B_i) \]

which is a pushout of $G(C) \coprod_{F(C)} F(D) \hookrightarrow G(D)$.

and a cofibration:

\[ G(D) \coprod_{F(A_i)} F(B_i) \hookrightarrow G(D) \coprod_{G(A_i)} G(B_i) \]

which is the pushout of \[ G(A_i) \coprod_{F(A_i)} F(B_i) \hookrightarrow G(B_i) \]

along the map $G(A_i) \rightarrow G(D)$. The composite of the two cofibrations with the isomorphism above is our map $G(C) \coprod_{F(C)} F(D^+) \hookrightarrow G(D^+)$ hence this concludes the proof.
}
}

\block{\Prop{Let $\Ccal$, $\Ccal'$ and $\Dcal$ be three pre-cylinder categories.
A morphism from $\Ccal'$ to $[\Ccal,\Dcal]$ is the same as a functor $F: \Ccal' \times \Ccal \rightarrow \Dcal$ such that:

\begin{itemize}

\item For each object $c \in \Ccal$, $F_c: x \mapsto F(x,c)$ is a morphism from $\Ccal'$ to $\Dcal$.

\item For each object $c' \in \Ccal$, $F_{c'}: x \mapsto F(c',x)$ is a morphism from $\Ccal$ to $\Dcal$.

\item If $i: a \hookrightarrow b$ is a cofibration in $\Ccal$ and $i:a' \hookrightarrow b'$ is a cofibration in $\Ccal'$ then the natural map induced by $i$ and $i'$:

\[ F(a',b) \coprod_{F(a',a)} F(b',a) \rightarrow F(b',b) \]

is a cofibration.

\end{itemize}

}

The proof is a routine check. Such a functor $F :\Ccal' \times \Ccal \rightarrow \Dcal$ will be called a bi-morphism.

}

\block{\label{Prop_TensorProdConstr}\Prop{For any two pre-cylinder categories $\Ccal$ and $\Ccal'$, there is a essentially unique pre-cylinder category $\Ccal \otimes \Ccal'$ endowed with a bi-morphism $\_ \otimes \_ : \Ccal \times \Ccal' \rightarrow \Ccal \otimes \Ccal'$ such that for any other pre-cylinder category $\Dcal$ morphisms from $\Ccal \otimes \Ccal'$ to $\Dcal$ are the same as bi-morphism $\Ccal \times \Ccal' \rightarrow \Dcal$. }

\Dem{One can consider the pre-cylinder category $\Ccal \otimes \Ccal'$ ``freely generated by a bi-morphism from $\Ccal \times \Ccal'$''. Indeed such a bi-morphism can be described as a choice of:

\begin{itemize}

\item An object $c \otimes c'$ for each pair of objects $c \in \Ccal$ and $c' \in \Ccal$.
\item A morphisms: $f \otimes f': c_1 \otimes c'_1 \rightarrow c_2 \otimes c'_2$ for all pair of arrow $f \in \Ccal$ and $f'\in \Ccal$.
\item relation between these morphisms translating composition in $\Ccal$ and in $\Ccal'$.
\item Some of the maps are imposed to be cofibrations/weak equivalence to translate the fact that $c \otimes \_$ and $\_ \otimes c$ should preserve weak equivalences and cofibration as well as the last conditions in the definition of bi-morphism.
\item The natural maps $0 \rightarrow 0 \otimes c$ and $0 \rightarrow c \otimes 0$ are forced to be isomorphisms.

\item Certain comparison maps from a pushout to a cone should be forced to be isomorphisms in order to impose that $c \otimes \_$ and $\_ \otimes c$ preserve pushout along cofibrations.

\end{itemize}

and this fits in the framework of \ref{SectionFreeCofCat}. The corresponding free pre-cylinder category $\Ccal \otimes \Ccal'$ satisfies the universal property mentioned in the proposition.
}

Note that in particular, the underlying cofibrations category of $\Ccal \otimes \Ccal'$ corresponds to the tensor product of the underlying cofibrations category of $\Ccal$ and $\Ccal'$. This follow from the fact that the underlying cofibration category of the free pre-cylinder category $\Ccal \otimes \Ccal'$ is obtained by just forgeting everything about weak equivalences in the system of generators and this corresponds exactly to the notion of a bi-morphism of cofibration categories from $\Ccal \times \Ccal'$.

}

\block{If $i:A \hookrightarrow B$ and $i' : A'\hookrightarrow B'$ are respectively cofibrations in $\Ccal$ and $\Ccal'$ then one denote by $ i' \square i'$ the cofibrations:

\[ A \otimes B' \coprod_{A \otimes A'} B \otimes A' \hookrightarrow B \otimes B' \]

in $\Ccal \otimes \Ccal'$.

}

\block{\label{GenCofTensProd}\Prop{If $\Ccal$ and $\Ccal'$ are two pre-cylinder categories with classes of generating cofibrations $I$ and $I'$ then the $i \square i'$ for $i \in I$ and $i' \in I'$ are a class of generating cofibrations for $\Ccal \otimes \Ccal'$.}

\Dem{ We will first show that if $j:A \hookrightarrow B$ and $j' : A' \hookrightarrow B'$ are respectively cofibrations in $\Ccal$ and $\Ccal'$ then $j \square j'$ is a iterated pushout of the $i \square i'$.

Let $B^+ = B \coprod_C D$ where $i:C \hookrightarrow D$ is one of the generating cofibration of $\Ccal$ and $C \rightarrow B$ is any map, let $j^+$ be the composite cofibration $A \hookrightarrow B^+$. Then in the following diagram the two square are pushouts and the total vertical arrow is $j^+ \square j'$:
\[\begin{tikzcd}[column sep={10pt},ampersand replacement=\&]
 \&  \displaystyle \left(B^+ \otimes A'\right) \coprod_{A \otimes A'} \left( A \otimes B'\right) \arrow[hook]{d} \& \displaystyle \left( B \otimes A' \right) \coprod_{A \otimes A'} \left( A \otimes B' \right) \arrow{l} \arrow[hook]{d}{j \square j'} \\
\displaystyle \left(D \otimes A' \right) \coprod_{C \otimes A' } \left( C \otimes B' \right) \arrow{r} \arrow[hook]{d}{i \square j'} \& \displaystyle \left(B^+ \otimes A' \right) \coprod_{B \otimes A'} \left( B \otimes B' \right) \arrow[hook]{d}   \& B \otimes B' \arrow{l} \\
 D \otimes B' \arrow{r}\& B^+ \otimes B' \& \\
\end{tikzcd}
\]

Hence $j^+ \square j'$ is an iterated pushout of $j \square j'$ and $i \square j'$, moreover if either $j$ and $j'$ is an isomorphism then $j \square j'$ is an isomorphism, hence one can proceeds by induction on gluing of generating cofibrations, we first show that if $j'$ is a generating cofibration and $j$ is any cofibration then $j \square j'$ is an iterated pushout of the $i \square j'$ for $i$ a generating cofibration and then reversing the role of $j$ and $j'$ that for any pair of cofibrations $j$ and $j'$ the cofibrations $j \square j'$ is an iterated pushout of the generating cofibrations $i \square i'$.

From the construction of $\Ccal \otimes \Ccal'$ as a free pre-cylinder categories given in \ref{Prop_TensorProdConstr} and the general discussion of \ref{FreePrecylinderCatMain}, the generating cofibrations of $\Ccal \otimes \Ccal'$ are the cofibrations of the form $i \square i'$ together with the cofibration in $\Ccal$ tensored by an object of $\Ccal'$ and conversely, but those are special case of $j \square j'$ when one of the two cofibrations $j$ and $j'$ is of the form $0 \hookrightarrow c$ hence this shows that all the generating cofibrations of $\Ccal \otimes \Ccal'$ are iterated pushout of the $i \square i'$ and hence that the $i \square i'$ form a system of generating cofibrations for $\Ccal \otimes \Ccal'$.

}
}

\block{\Prop{A morphism from $\Ccal \otimes( \Ccal' \otimes \Ccal'')$ to $\Dcal$ (with $\Ccal, \Ccal',\Ccal''$ and$ \Dcal$ four pre-cylinder categories) is the same as a functor $F: \Ccal \times \Ccal' \times \Ccal'' \rightarrow \Dcal$ such that:

\begin{itemize}
\item If one fixes any two of the variables one gets a morphism of pre-cylinder category in the third variables.

\item For any triplet $i:A_0 \hookrightarrow A_1$, $i' :B_0 \hookrightarrow B_1$ and $i'' :C_0 \hookrightarrow C_1$ of cofibrations respectively in $\Ccal$, $\Ccal'$ and $\Ccal''$, the morphism:

\[ \underset{(u,v,w) \in \{ 0,1\}^3-\{1,1,1\}}{\text{colim}} F(A_u,B_v,C_w) \rightarrow F(A_1,B_1,C_1) \]
 
is a cofibration, where $\{ 0,1\}^3-\{1,1,1\}$ is considered as an ordered set with the lexicographic order.

\end{itemize}

}

The only point of this proposition is that this description is completely symmetric in the three variables, hence proving that the tensor product is associative (up to canonical coherence isomorphisms) hence making the category  of pre-cylinder categories (as well as the category of cofibrations categories) into a symmetric monoidal $2$-category. The unit object of the tensor product is the free pre-cylinder category on one object $F_*$ because a morphism from $F_* \otimes \Ccal$ to $\Dcal$ is the same as a morphism from $F_*$ to $[\Ccal,\Dcal]$ is an object of $[\Ccal,\Dcal]$, i.e. a morphism from $\Ccal$ to $\Dcal$, hence $F_* \otimes \Ccal$ is isomorphic to $\Ccal$.

\Dem{Such a morphism corresponds to a morphism $\Ccal \rightarrow [\Ccal' \otimes \Ccal '',\Dcal]$.

Objects of $[\Ccal' \otimes \Ccal '',\Dcal]$ are bi-morphisms $\Ccal \times \Ccal' \rightarrow \Dcal$ and morphisms between them are the natural transformations.

Cofibrations in $[\Ccal' \otimes \Ccal '',\Dcal]$ can be described using proposition \ref{CritForCofibinHom} and proposition \ref{GenCofTensProd} together, as the natural transformations $\lambda:F \rightarrow G$ such that for each $i:B_0 \hookrightarrow B_1$ and $j:C_0 \hookrightarrow C_1$ cofibrations in $\Ccal'$ and $\Ccal''$ the comparison map between the pushout of the first three object to the last one in the diagram:

\[ \begin{tikzcd}[ampersand replacement=\&]
F(B_0,C_1) \coprod_{F(B_0,C_0)} F(B_1, C_0) \arrow[hook]{r} \arrow{d}{\lambda}\& F(B_1,C_1) \arrow{d}{\lambda} \\
G(B_0,C_1) \coprod_{G(B_0,C_0)} G(B_1, C_0) \arrow[hook]{r} \& G(B_1,C_1) \\
\end{tikzcd}
\]

is a cofibrations.

From this one immediately deduces that the second condition in the proposition above exactly says that cofibrations in $\Ccal$ are sent to cofibrations in $[\Ccal' \otimes \Ccal'',\Dcal]$.

The rest of the verifications are a long but easy routine check.
}
}

\blockn{At this point, the following theorem becomes evident:}

\block{\Th{The objects constructed in this sections endows the category of pre-cylinder categories with the structure of a symmetric closed monoidal $2$-category.
}}

\section{A weak model structure on the category of pre-cylinder categories}
\label{SecWeakModelStr}
\blockn{In this section we will construct some sort of weak model structure on the category of pre-cylinder categories. It is not a real model structure but it will be more than enough to construct the homotopically free objects that we need.

We will call acyclic (co)fibration the maps that are both (co)fibrations and weak equivalences. They might be different from trivial (co)fibrations which are the maps satisfying some lifting properties (although some comparison results between those classes will be proved showing that they are often the same).
}

\subsection{Fibrations and Cofibrations of pre-cylinder categories}
\label{SubSecFibCofib}

\block{\label{Def_FibPCylcat}\Def{
A morphism $F :\Ccal \rightarrow \Dcal$ between pre-cylinder categories will be called a pre-fibration if the retraction of trivial cofibrations and the factorization in cofibration followed by weak equivalences can be lifted along $F$. More precisely, $F$ is a pre-fibration if:

\begin{itemize}

\item If $s: A \overset{\sim}{\hookrightarrow} B$ is a trivial cofibration in $\Ccal$ and $r: F(B) \rightarrow F(A)$ is a retraction of $F(s)$ in $\Dcal$, then $s$ admit a retraction $r':B \rightarrow A$ such that $r=F(r')$.

\item If $\lambda: A \rightarrow B$ is a morphism in $\Ccal$ and:
 
 \[ F(A) \underset{i}{\hookrightarrow }P \overset{\sim}{\underset{w}{\rightarrow}} F(B) \]
 
is a factorization of $F(\lambda)$ into a cofibration followed by a weak equivalence then there is a factorization of $\lambda$ as a cofibration $i'$ followed by a weak equivalence $w'$ such that $w=F(w')$ and $i=F(i')$.

\end{itemize}

}}

\block{A pre-cylinder category $\Ccal$ is a cylinder category if and only the map $\Ccal \rightarrow *$ is a pre-fibration, indeed the pre-cylinder category $*$ has all the retraction and all cofibration/weak equivalence factorisations.}

\block{\Def{\begin{itemize}

\item A pre-cylinder category $\Ccal$ will be said to be fibrant if the map $\Ccal \rightarrow *$ is a pre-fibration, i.e. if $\Ccal$ is a cylinder category.

\item A trivial cofibration is a morphism which has the left lifting property with respect to all pre-fibrations between fibrant objects.

\item A cofibration is a morphism which has the left lifting property with respect to all acyclic pre-fibrations between fibrant objects.

\item A fibration is a morphism which has the right lifting property with respect to all trivial cofibrations.

\item A trivial fibration is a morphism which has the right lifting property with respect to all cofibrations.

\item An anodyne morphism is a morphism which has the left lifting property with respect to all pre-fibrations.

\end{itemize}
}}

\blockn{Here are some immediate (and very easy) consequences of these definitions:}

\block{\label{trivcof-fib-facto}\Prop{\begin{enumerate}

\item Anodyne morphisms and pre-fibrations form a weak factorization system.

\item Anodyne maps are in particular trivial cofibrations.

\item Fibrations are in particular pre-fibrations.

\item A map with fibrant co-domain, is a fibration if and only if it is a pre-fibration.

\item Any morphism with fibrant co-domain can be factored as a trivial cofibration followed by a fibration.

\end{enumerate}
}
\Dem{

\begin{enumerate}

\item The key observation here is that being a pre-fibration can be expressed as having the right lifting property with respect to two specific morphisms between free pre-cylinder category: the morphism from the pre-cylinder category freely generated by a trivial cofibration to the pre-cylinder category freely generated by a trivial cofibration with a retraction, and the morphism from the pre-cylinder category freely generated by a morphism $\lambda: A\rightarrow B$ to the pre-cylinder category freely generated by a diagram:

\[ A \hookrightarrow I \overset{\sim}{\rightarrow} B \]

Hence the result follow directly from the small object argument.

\item Anodyne maps have the left lifting property with respect to all pre-fibrations, in particular with respect to the pre-fibrations between fibrant objects and hence are trivial cofibrations.

\item Fibrations have the right lifting property with respect to all trivial cofibrations, hence in particular with respect to anodyne morphisms and hence (as anodyne morphisms and pre-fibrations form a weak factorization system) they are pre-fibrations.

\item Let $f:A \rightarrow B$ be a pre-fibration with fibrant domain, hence $A$ is fibrant as well, hence $f$ has the right lifting properties against all trivial cofibrations by definition of trivial cofibrations and hence it is a fibration. The converse is the previous point.

\item Factor it has an anodyne map followed by a pre-fibration, the anodyne map is a trivial cofibration and the pre-fibration has a fibrant co-domain hence it is a fibration.

\end{enumerate}

}

}

\block{\Def{A full subcategory $\Ccal \subset \Ccal'$ of a pre-cylinder category is said to be \emph{h-saturated} if it is stable under pushout by cofibrations and if any object in $\Ccal'$ which is weakly equivalent to an object in $\Ccal$ is in $\Ccal$.}

The following proposition is immediate:

\Prop{If $\Ccal \subset \Ccal'$ is h-saturated then $\Ccal$ is a pre-cylinder category for the induced structure and the inclusion functor is a pre-fibration.}

In particular, if $\Ccal'$ is a cylinder category then $\Ccal$ is also a cylinder category.

}

\block{\label{FibrationCoslice}\Prop{Let $F:\Ccal \rightarrow \Dcal$ be a pre-fibration, let $A \in \Ccal$ then the induced morphism $F_A:\Ccal(A) \rightarrow \Dcal(F(A))$ is also a pre-fibration. }

\Dem{Retraction of a trivial cofibration in $\Ccal(A)$ is the same as a retraction of the underlying map in $\Ccal$ and a cofibration/weak equivalence factorization of a map in $\Ccal(A)$ is the same such a factorization of the underlying map in $\Ccal$ hence the lifting of those for $F:\Ccal \rightarrow \Dcal$ immediately gives the same lifting property for $F_A$.}
}

\block{\label{MoreLiftingForPreFib}\Prop{Let $F : \Ccal \rightarrow \Dcal$ be a pre-fibration. Then $F$ has the following lifting properties:

\begin{enumerate}

\item If the solid diagram:
\[
\begin{tikzcd}[ampersand replacement=\&]
A \arrow{r}{f} \arrow[hook]{d}{\sim}[swap]{i} \& X \\
B \arrow[dashed]{ur}{h} \& \\
\end{tikzcd}
\]
Is a diagram in $\Ccal$ (with $i$ a trivial cofibration) such that the dashed arrow $h$ exists in $\Dcal$, then there is such a dashed arrow $h'$ in $\Ccal$ such that $F(h')=h$.

\item If $A \hookrightarrow X$ is a cofibration in $\Ccal$ and $I_{F(A)}F(X)$ is a relative cylinder object for $F(A) \hookrightarrow F(X)$ in $\Dcal$ then there exists a cylinder object $I_A X$  in $\Ccal$ such that $F(I_AX) = I_{F(A)}F(X)$

\item If:

\[
\begin{tikzcd}[ampersand replacement=\&]
A \arrow{r}{f} \arrow[hook]{d}{i} \& X \\
B \arrow{ur}{h} \& \\
\end{tikzcd}
\]

Is a diagram in $\Ccal$, and assume that there is another map $h_2$ which makes the diagram:

\[
\begin{tikzcd}[ampersand replacement=\&]
F(A) \arrow{r}{F(f)} \arrow[hook]{d}{F(i)} \& F(X) \\
F(B) \arrow{ur}{h_2} \& \\
\end{tikzcd}
\]

commutes in $\Dcal$, and such that $h_2 \sim_{F(A)} F(h) $ for some relative cylinder object (existing) in $\Dcal$. Then there is an arrow $h' : B \rightarrow X$ which can replace $h$ in the first diagram, such that $h' \sim_A h$ (for some lifting of the relative cylinder object in $\Dcal$) and $F(h')=h_2$.

Moreover if one fixes the homotopy between $h_2$ and $F(h)$ one can chose $h'$ and an homotopy between $h'$ and $h$ which is sent to the the chosen homotopy by $F$.

\end{enumerate}

 }

\Dem{

\begin{enumerate}

\item Such a diagonal filler is the same as a retraction of the trivial cofibration $X \hookrightarrow X \coprod_{A} B$, hence this follows from the fact that pre-fibrations lift retraction of trivial cofibrations.

\item A relative cylinder object is a special case of cofibration/weak equivalence factorization and those are lifted by pre-fibrations.

\item Let us first chose a relative cylinder object $I_A B$ in $\Ccal$ that lifts the one in $\Dcal$ used to define the homotopy $F(h) \sim_{F(A)} h_2$, we also fix such an homotopy $F(h) \sim_A h_2 : F(I_A B) \rightarrow F(X)$. The situation can then be presented as follow, one has a solid diagram in $\Ccal$:

\[
\begin{tikzcd}[ampersand replacement=\&]
B \arrow{r}{h} \arrow[hook]{d}{\sim} \& X \\
I_A B \arrow[dashed]{ur}[swap]{F(h) \sim_{F(A)} h_2} \& \\
\end{tikzcd}
\]

with a dashed arrow in $\Dcal$ given by the map $h_2$ together with the homotopy $F(h) \sim_{F(A)} h_2$ as the vertical map is a trivial cofibration this can be lifted to $\Ccal$ producing a map $h'$ over $h_2$ and an homotopy $h \sim_A h'$ over the homotopy between $F(h)$ and $h_2$.

\end{enumerate}

} 
}

\block{\label{Prop_LiftingOfWeakDiag}\Prop{Let $F: \Ccal \rightarrow \Dcal$ be a fibration of pre-cylinder categories.
Assume that one has a square in $\Ccal$ (with a cofibration and a weak equivalence):

\[\begin{tikzcd}[ampersand replacement=\&] 
A \arrow{r} \arrow[hook]{d}{i} \& X \arrow{d}{\sim} \\
B \arrow{r} \& Y \\
\end{tikzcd}
\]

And assume that its image by $F$ can be extended in $\Dcal$ into:

\[\begin{tikzcd}[ampersand replacement=\&] 
F(A) \arrow[hook]{rd}{F(i)} \arrow{rr} \arrow[hook]{ddd}{F(i)}\& \& F(X) \arrow{ddd}{\sim} \\
\& F(B) \arrow[hook]{d}{i_0} \arrow{ru}{f} \& \\
\& I_{F(A)} F(B) \arrow{rd}{h}\& \\
F(B) \arrow[hook]{ru}{i_1} \arrow{rr}\& \& F(Y) \\
\end{tikzcd}
\]

Where $I_A B$ is a relative cylinder object of $i$ in $\Dcal$, then the extension can be lifted to $\Ccal$.

}

\Dem{This can be expressed as the fact that all fibrations have the right lifting property with respect to a certain morphism between freely generated pre-cylinder categories (freely generated by the diagram above and by a relative cylinder object). Hence the proposition boils down to saying that this map is a trivial cofibration, and hence it is enough to check the lifting property when $\Ccal$ and $\Dcal$ are both cylinder categories.

The fact that the relative cylinder object can be lifted has been proved in proposition \ref{MoreLiftingForPreFib}.

We proved in \ref{Prop_Semilifting2} that such filing always exists in a cylinder category, so there is such a filing in $\Ccal$, using the lifted relative cylinder object. We call $f'$ and $h'$ the map obtain in $\Ccal$ this way. We need to prove that one can replace $f'$ and $h'$ by maps that are send to $f$ and $h$ and this will be done mostly by using the last point of proposition \ref{MoreLiftingForPreFib}.

Let us call $v$ the equivalence $X \rightarrow Y$, then  $F(v) \circ F(f')$ and $F(v) \circ f$ are both homotopic relative to $A$ to the map $F(k):F(B) \rightarrow F(Y)$ corresponding to the bottom arrow of the square, the homotopy being respectively given by $h$ and $F(h')$, in particular they are homotopic relative to $A$ by composing these two homotopies, but as $F(v)$ is a weak equivalence it shows that $F(f')$ and $f$ are homotopic relative to $F(A)$, we will call $h_f$ this homotopy, chosen such that $F(v) \circ h_f$ is homotopic relative to $F(B) \coprod_{F(A)} F(B)$ to the composite of $F(h')$ and $h$. Applying the last point of proposition \ref{MoreLiftingForPreFib} one can lift $f$ and $h_f$ to maps in $\Ccal$, which we will also call $f$ and $h_f$ (the maps in $\Dcal$ formerly known as $f$ and $h_f$ are now $F(f)$ and $F(h_f)$.

It remains to show that the homotopy $h$ between $F(v \circ f)$ and $F(k)$ can also be lifted to a homotopy in $\Ccal$, but one already has a homotopy between $v \circ f$ and $k$ in $\Ccal$: it can be obtained as the composite of the homotopy $h_f$ between $f$ and $f'$ composed with $v$ and the homotopy $h'$ between $k$ and $v \circ f'$, moreover it is not very hard to see (manipulating composition of homotopies and inverse of homotopies which satisfies the usual groupoids axiom up to higher homotopy) that the image of this homotopy in $\Dcal$ is homotopy equivalent, relative to $F(B) \coprod_{F(A)} F(B)$, to $h$ hence proving (applying the last point of proposition \ref{MoreLiftingForPreFib} one more time) that $h$ can be lifted to $\Ccal$ which concludes the proof.
}
}

\blockn{In order to obtain a similar (partial) weak factorization system for cofibrations and trivial fibrations one needs to observe that acylic fibrations between fibrant objects can also be detected by the right lifting property with respect to a set of maps. This is achieved by the following proposition:}

\block{\label{LiftingForAcyclicFib}\Prop{A morphism $F: \Ccal \rightarrow \Dcal$ of pre-cylinder category with fibrant co-domain is an acyclic fibration if and only if it satisfies the following condition:

\begin{itemize}

\item It detects weak equivalences.

\item If $X$ is an object in $\Ccal$ and $i:F(X) \hookrightarrow Y$ is a cofibration in $\Dcal$ then there is a cofibration $i':X \hookrightarrow Y'$ such that $F(i')=i$.

\item If one has a diagram in $\Ccal$ of the form:

 \[
\begin{tikzcd}[ampersand replacement=\&]
A \arrow{r}{f} \arrow[hook]{d}{i} \& X \\
B \& \\
\end{tikzcd}
\]

with $i$ a cofibration, such that it admits a filling in $\Dcal$:

 \[
\begin{tikzcd}[ampersand replacement=\&]
F(A) \arrow{r}{F(f)} \arrow[hook]{d}{F(i)} \& F(X) \\
F(B) \arrow{ur}{h}\& \\
\end{tikzcd}
\]

Then there is such a filling $h'$ in $\Ccal$ such that $h=F(h')$.

\end{itemize}

Moreover, if the domain of $F$ is also fibrant the first conditions can be omitted.

 }
 
Note that the first condition is also a right lifting property: it corresponds to the right lifting property against the map from the free pre-cylinder category on one arrow to the free pre-cylinder category on one weak equivalence.
 
\Dem{ We first assume that $F$ satisfies the lifting property of the proposition.

If we know that the domain is fibrant, then those conditions imply that $F$ is acyclic: it is homotopy surjective because of the first lifting property and homotopy fully faithful because of the second condition. This in particular implies that $F$ detect weak equivalences (because $hF$ is an equivalence of category and hence detect isomorphisms) hence in this case the first conditions can indeed be omitted.

In the general case, if $F(A) \hookrightarrow P \overset{\sim}{\rightarrow} F(B)$ is a factorization as a cofibration followed by a weak equivalence in $\Dcal$ of a map in $\Ccal$, then one can lift $P$ and the cofibration $F(A) \hookrightarrow P$ using the second lifting property, then one can use the third lifting property to lift the map $P \rightarrow F(B)$, and finally as $F$ detect weak equivalences this implies that the lifting is a weak equivalence.

Finally, the lifting conditions immediately implies the lifting of retractions of trivial cofibrations (it is a special case of the second lifting property). So this proves that a map satisfying those lifting conditions is a fibration. In particular the domain is fibrant hence it is acyclic.

\bigskip 

Conversely, let $F :\Ccal \rightarrow \Dcal$ be an acyclic fibration between cylinder categories. As $F$ is acyclic it detect weak equivalences and it satisfies the weak form of the two lifting conditions that we need to prove. The fact that it is a fibration will allow to make then strict. More precisely: 

If $A \in \Ccal$ and $F(A) \hookrightarrow B$ is a cofibration then as $F_A : h\Ccal(A) \rightarrow h\Dcal(A)$ is essentially surjective, there is an object $B'$ of $\Ccal$ with a cofibration $A \hookrightarrow B'$ such that $F(B')$ is equivalent to $B$ as objects of $\Dcal(A)$. One can hence see $F(A) \hookrightarrow B \overset{\sim}{\rightarrow} F(B')$ as a cofibration/weak equivalences factorization of the cofibration $F(A \hookrightarrow B')$ which as such can be lifted along the fibration $F$ and hence produces in particular a lifting of the cofibration $F(A) \hookrightarrow B$ as desired.

The third lifting property is obtained similarly: if one has such a diagram in $\Ccal$ with a filling in $\Dcal$ then as $F$ is homotopy fully faithful there is a filling in $\Ccal$ which is homotopy equivalent (relative to the the domain of the cofibration) in $\Dcal$ to the filing that one already has, hence by the point $3.$ of proposition \ref{MoreLiftingForPreFib} one can find a lifting of the original filling.
}  
}

\block{\label{CorFactoCofTrivFib}\Cor{Any morphism of pre-cylinder categories $f : \Ccal \rightarrow \Dcal$ with fibrant co-domain can be factored as a cofibration followed by a trivial fibration.}

\Dem{This is exactly the same as proposition \ref{trivcof-fib-facto} but replacing pre-fibration with the class of map satisfying the lifting properties of \ref{LiftingForAcyclicFib}. We will call those maps ``trivial pre-fibration''\footnote{The proof of proposition \ref{trivcof-fib-facto} show that they are in particular pre-fibrations.}. Trivial pre-fibrations are characterized by the right lifting property with a set of morphism: the one from the free pre-cylinder category on one arrow to the free pre-cylinder category on one weak equivalences, and the two free constructions corresponding to the two other lifting properties of proposition \ref{LiftingForAcyclicFib}.

Applying the small object argument any map $f$ can be factored in a map that satisfies the left lifting property with respect to trivial pre-fibration followed by a trivial pre-fibration. The first map is in particular a cofibration: it has the left lifting property with respect to all acyclic fibration with fibrant co-domain because they are trivial pre-fibrations, and if the map $f$ has a fibrant co-domain then the trivial pre-fibration will be an acylic fibration by proposition \ref{LiftingForAcyclicFib}.
}
}

\subsection{Reedy type categories}
\label{SubSecReedy}

\block{\Def{A morphism $f :\Ccal \rightarrow \Ccal'$ of pre-cylinder category is said to be a \emph{Reedy extension} if for any fibration $F:\Dcal \rightarrow \Dcal'$ between cylinder categories, the morphism:

\[  [\Ccal',\Dcal] \rightarrow [\Ccal, \Dcal] \times_{[\Ccal,\Dcal']} [\Ccal',\Dcal'] \]

is a pre-fibration.

A pre-cylinder category is said to be of Reedy type if $\emptyset \rightarrow \Ccal$ is a Reedy extension.

}

The point of this notion is that if $\Ccal$ is a Reedy type pre-cylinder category then for all cylinder category $\Dcal$, the functor category $[\Ccal,\Dcal]$ is a cylinder category, and for all fibrations between cylinder category $\Dcal \rightarrow \Dcal'$, the morphism $[\Ccal,\Dcal] \rightarrow [\Ccal,\Dcal']$ is a fibration. We will think of Reedy type pre-cylinder categories as being ``sketches'' for some sort of theory whose category of models in a cylinder category will be a cylinder category. 

}

\block{\Prop{Reedy extension are stable under pushout, retract and transfinite composition.}

\Dem{ Manipulating the adjunction between internal hom objects of pre-cylinder categories and the tensor product of pre-cylinder categories constructed in subsection \ref{SubSecmonoidalStr}, one can see that for morphisms $f : \Ccal \rightarrow \Ccal'$, $F: \Dcal \rightarrow \Dcal'$, and $i : \Acal \rightarrow \Acal'$ saying that the map  
\[ [\Ccal',\Dcal] \rightarrow [\Ccal, \Dcal] \times_{[\Ccal,\Dcal']} [\Ccal',\Dcal'] \]
induced by $f$ and $F$ has the right lifting property against $i$ is the same as saying that the map 

\[ [\Acal',\Dcal] \rightarrow [\Acal, \Dcal] \times_{[\Acal,\Dcal']} [\Acal',\Dcal'] \]

induced by $i$ and $F$ has the right lifting property with respect to $f$. Hence $f$ is a Reedy extension if and only it has the left lifting property with respect to all the map of the form:

\[ [\Acal',\Dcal] \rightarrow [\Acal, \Dcal] \times_{[\Acal,\Dcal']} [\Acal',\Dcal'] \]

for $F: \Dcal \rightarrow \Dcal'$ a fibration and $i:\Acal \rightarrow \Acal'$ a anodyne morphism (or just one of the generating morphisms of the class of anodyne morphisms).

In particular, the class of Reedy extension is stable under retract, pushout and transfinite composition simply because it is characterized by a left lifting property.
}
}

\blockn{So we just have to give a few examples of Reedy extension, and then the proposition above will allows to produce plenty of examples by transfinite iterated pushout.}

\block{\label{KeyExemplesReedy}\Prop{The following maps are Reedy extension:

\begin{itemize}

\item The map from $F_*$ the free pre-cylinder category on one object to $F_{\hookrightarrow}$ the free pre-cylinder category on one cofibration, which send the generating object to the source of the generating cofibrations is a Reedy extension.

\item The map from the free pre-cylinder category on one arrow to the free pre-cylinder category on one equivalences. 
\end{itemize}

}

Note that $0 \rightarrow F_*$ is a pushout of the first example (along the unique morphism $F_* \rightarrow 0$) hence it is also a Reedy type pre-cylinder category.

\Dem{
\begin{itemize}

\item  Let $\Dcal$ be a cylinder category, $[F_*,\Dcal]$ is just $\Dcal$, and $[F_{\hookrightarrow},\Dcal]$ is the category that we will denote $\Dcal_{\downarrow}$ whose objects are pairs  $A=(A_0,A_1)$ of objects of $\Dcal$ together with a cofibration $A_0 \hookrightarrow A_1$, morphisms are the diagram morphisms,  weak equivalences are objectwise weak equivalences and cofibrations are the map $f:A \rightarrow B$ such that $f_0:A_0 \rightarrow B_0$ is a cofibration and the comparison map $A_1 \coprod_{A_0} B_0 \rightarrow B_1$ is a cofibration (this follows from proposition \ref{CritForCofibinHom} and the description of the generating cofibrations of a free pre-cylinder categories given in \ref{FreePrecylinderCatMain}).

Let now $F:\Dcal \rightarrow \Dcal'$ be a fibration between cylinder categories, we need to show that the morphism:

\[ H: \Dcal_{\downarrow} \rightarrow \Dcal'_{\downarrow} \times_{\Dcal'} \Dcal \]

is a prefibration.

Let $v :X \rightarrow Y$ be a morphism in $\Dcal_{\downarrow}$ such that one has a cofibration/weak equivalence factorization:

\[ X_0 \underset{i}{\hookrightarrow} Z_0 \overset{\sim}{\underset{w}{\rightarrow}} Y_0 \]

in $\Dcal$, such that its image by $F$ can be extended into a cofibration/weak equivalences factorization of $F(v):F(X) \rightarrow F(Y)$ in $\Dcal'_{\downarrow}$, this means that one has in $\Dcal'$:

\[ F(X_1 \coprod_{X_0} Z_0 ) \hookrightarrow Z_1 \overset{\sim}{\rightarrow} F(Y_1)\]

lifting this cofibration/weak equivalence factorization to $\Dcal$ (using that $F$ is a fibration) produces exactly what we need to complete the factorization $X_0 \hookrightarrow Z_0 \overset{\sim}{\rightarrow} Y_0$ into a factorization in $\Dcal_{\downarrow}$. Hence $H$ lifts cofibration/weak equivalence factorizations.

If one has a trivial cofibration $i:X \hookrightarrow Y$ in $\Dcal^{\downarrow}$ with a retraction $r_0 : Y_0 \rightarrow X_0$ of $i_0$, then finding a retraction of $i$ that extend $r_0$ is the same as finding a retraction to the trivial cofibration:

\[ X_1 \coprod_{X_0} Y_0 \hookrightarrow Y_1 \]

Hence if one can find such a retraction in $\Dcal'$ it can be lifted to $\Dcal$ and this shows that $H$ lift retract of trivial cofibrations.

\item Let $F_{\rightarrow}$ and $F_{\sim}$ be the two free categories mentioned in the proposition. 

$[F_{\rightarrow},\Ccal]$ is the category of arrows in $\Ccal$ with objectwise weak equivalences and objectwise cofibrations, $[F_{\sim},\Ccal]$ is the full subcategory of weak equivalences.

If $F:\Dcal \rightarrow \Dcal'$ is a fibration then the map:

\[ [F_{\sim},\Dcal] \rightarrow [F_{\rightarrow},\Dcal] \times_{[F_{\rightarrow},\Dcal']} [F_{\sim},\Dcal'] \]

Is the inclusion of the subcategory of equivalences in $\Dcal$ into the subcategory of arrows in $\Dcal$ whose image in $\Dcal'$ are weak equivalences, hence it is a h-saturated subcategory this concludes the proof

\end{itemize}

}

}

\blockn{The second example show that adding new equivalences (while keeping the same underlying cofibration category) to a Reedy type category gives a Reedy category. The first examples in the proposition above is exactly what we need to produce all the ``Reedy type structure'' for directed category as they were defined in \cite{cisinski2010categories} or \cite{radulescu2006cofibrations}, which is what we will now explain:}

\block{\Def{A locally finite directed category is a category $\Ical$ endowed with a ``height'' function $h:\Ical \rightarrow \mathbb{N}$ such that any non-identity arrow $f:i \rightarrow j$ in $\Ical$ is such that $h(i)<h(j)$ and such that all the slice category $\Ical_{/x}$ are finite.
}
}

\block{\label{ReedyStructure}\Prop{Let $\Ical$ be a locally finite directed category, then:

\begin{itemize}

\item The category $\widehat{\Ical}_f$ of finite presheaves\footnote{We mean by that presheaves which are non empty on only a finite number of element of $\Ical$ and which are finite on those elements} over $\Ical$, with monomorphisms as cofibration and isomorphisms as weak equivalence is a pre-cylinder category.

\item $\widehat{\Ical}_f$ is a Reedy type pre-cylinder category.

\item The completion of $\widehat{\Ical}_f$ is the category $\widehat{\Ical}$ of pre-sheaves over $\Ical$.

\item Morphisms from $\widehat{\Ical}_f$ to any pre-cylinder category $\Ccal$ are exactly the ``Reedy cofibrant $I$-diagram''\footnote{This can be read as ``certain functors from $\Ical$ to $\Ccal$ with cofibrancy conditions that we will define in the proof'', but it do corresponds to the usual and well known notion of Reedy cofibrant object.} in $\Ccal$. Moreover the functor corresponding to a diagram is its left Kan extension.

\item If $\Jcal \subset \Ical$ is a full subcategory such that if $ i\rightarrow j$ is an arrow in $\Ical$ and $j \in \Jcal$ then $i \in \Ical$, then inclusion morphism:

\[ \widehat{\Jcal}_f \rightarrow \widehat{\Ical}_f \]

is a Reedy extension.

\end{itemize}

}

\Dem{The first claim is immediate: taking monomorphisms between finite sheaves clearly produces a cofibration category and taking isomorphisms as the weak equivalences makes it into a pre-cylinder category.

All the other results are obviously true for the case of $\Ical= \emptyset$, or maybe more clearly for $\Ical=\{* \}$.

We will prove the other claims for any finite category $\Ical$, by induction on the number of objects of $\Ical$.

Let $x$ be a maximal object in $\Ical$ (either such that $h(x)$ is maximal, or more generally such that there is no non-identity arrow in $\Ical$ whose domain is $x$), then $\Jcal = \Ical-\{x\}$ is a category satisfying the condition of the last point, which is also directed and with strictly less object than $\Ical$.

Let $\partial x$ be the object of $\widehat{\Jcal}_f$ ``represented by $x$''.

Let $\Ccal$ be the pre-cylinder category freely generated from $\Jcal$ by freely adding an object $x$ endowed with a cofibration $\partial x \hookrightarrow x$.

A morphism from $\Ccal$ to any other pre-cylinder category is given by first a morphism $F$ from $\widehat{\Jcal}_f$, i.e. a Reedy cofibrant $\Jcal$-diagram, together with an object $F_x$ and a cofibration $F(\partial x) \hookrightarrow F_x$. But as $F$ is a left Kan extension, a map from $F(\partial x)$ to $F_x$ is the same as a compatible collection of map $F(y) \rightarrow F_x$ for each $f:y \rightarrow x$ in $\Ical$, which makes $F$ into a functor from $\Ical$ with $F(x) = F_x$, the condition ``Reedy cofibrant'' can be described as $F$ is Reedy cofibrant when restricted to $\Jcal$ and the map $F(\partial x) \rightarrow F(x)$ is a cofibration.

We then use the general procedure of \ref{SectionFreeCofCat} to compute free categories, the completion of $\Ccal$ is the category of morphisms from $\Ccal$ to $\setop$, but by the above observations these are just presheaves over $\Ical$, $\Ccal$ is then the full subcategory of iterated pushout of: monomorphisms between finite $\Jcal$-presheaves (which are identified with certain $\Ical$-presheaves by taking $F(x)=\emptyset$ and by the cofibration $\delta x \rightarrow x$ which allow to add any finite number of elements we want to $F(x)$, hence $\Ccal$ is indeed $\widehat{\Ical}_f$.

As $\widehat{\Jcal}_f \rightarrow \Ccal= \widehat{\Ical}_f$ is a Reedy extension this concludes the proof for a finite category.

Indeed, for the last claim if $\Jcal$ is a subcategory satisfying this condition, then if $\Jcal = \Ical$ then the claim is immediate and if not one could have chosen $x$ not in $\Jcal$ and concludes by induction on the difference of cardinality between the cardinal of $\Jcal$ of $\Ical$.

If $\Ical$ is not finite but locally finite then it can be written as an increasing union of finite full subcategory $\Jcal$ satisfying the condition of the last claim, and $\widehat{I}_f$ appears as the colimit of the $\widehat{\Jcal}_f$ and hence is also a Reedy type category.

}
}

\blockn{We will also need the following proposition, which unfortunately does not seem to follow from our framework:
}

\block{\Prop{Let $\Ical$ be a locally finite directed category, $F:\Ical \rightarrow \Ccal$ a functor to a cylinder category, then there is a Reedy cofibrant functor $\tilde{F}:\Ical \rightarrow \Ccal$ endowed with a natural transformation $\tilde{F} \rightarrow F$ which is objectwise a weak equivalence.}

\Dem{We will construct $\tilde{F}$ gradually on objects of $\Ical$ following the same induction process as in the proof of proposition \ref{ReedyStructure}.

Let $\Jcal \subset \Ical$ be a full subcategory such that if $i \rightarrow j$ with $j \in \Jcal$ then $i \in \Ical$, and let $x$ an object of $\Ical - \Jcal$ such that all the (non identity) arrow with co-domain $x$ have domain in $\Jcal$. Assume that $\tilde{F}$ has already been constructed on $\Jcal$. We just need to extend $\tilde{F}$ to $x$.

Seeing $\tilde{F}$ as a morphism from $\widehat{\Jcal}_f$ to $\Ccal$, one can form the object $\tilde{F}(\partial x)$ in $\Ccal$. Because of the natural transformation $\tilde{F} \rightarrow F$, $F(X)$ is endowed with a compatible family of maps $\tilde{F}(y)$ for all arrow $y \rightarrow x$ in $\Ical$, hence there is a map $\tilde{F}(\partial x) \rightarrow x$, taking a cofibration/weak equivalence factorization of this maps exactly gives what we need for $\tilde{F}(x)$ and this concludes the proof.

}

}

\block{\Prop{Let $f : T \rightarrow T'$ be a Reedy extension of Reedy type categories which is also a trivial cofibration. Then for any cylinder category $\Ccal$, the  map:

\[ [f,\Ccal]: [T',\Ccal] \rightarrow [T,\Ccal] \]

is a trivial fibration.
}

\Dem{We already know that this map is a fibration between fibrant objects, so we just need to show that it is acyclic. Let $A \rightarrow B$ be an other Reedy extension of Reedy type pre-cylinder category category, and consider a diagram of the form:

\[
\begin{tikzcd}[ampersand replacement=\&]
A \arrow{r} \arrow{d} \& \lbrack T',\Ccal \rbrack \arrow{d}{\lbrack f,\Ccal \rbrack} \\
B \arrow{r} \& \lbrack T,\Ccal \rbrack \\
\end{tikzcd}
\]

the two horizontal arrow can be turned into arrow $T' \rightarrow [A,\Ccal]$ and $T \rightarrow [B,\Ccal]$ and the commutativity of the above diagram is equivalent to the commutativity of:

\[
\begin{tikzcd}[ampersand replacement=\&]
T \arrow{r} \arrow{d}{f} \& \lbrack B,\Ccal \rbrack \arrow{d} \\
T' \arrow{r} \& \lbrack A,\Ccal \rbrack \\
\end{tikzcd}
\]

but in this square the left arrow is a trivial cofibration and the right arrow is a fibration, hence there is a diagonal filler: $T' \rightarrow [B,\Ccal]$, which can be reinterpreted as a diagonal filler $B \rightarrow [T',\Ccal]$ in the first square.

This proves that our map has the right lifting property with respect to all Reedy extensions between Reedy type category. We will prove that $F=[f,\Ccal]$ is acyclic using the second characterization of proposition \ref{Charac_homotopyequiBetweenCylinder}. The case of the map $F_* \rightarrow F_{\hookrightarrow}$ of \ref{KeyExemplesReedy} show that for all object $A$, $F_A$ is homotopy surjective. The case of the morphism from $F_{\hookrightarrow}$ to the free category on a trivial cofibration shows that $F$ detect trivial cofibrations among cofibrations. Using the cofibration/weak equivalence factorization this is enough to deduce that $F$ detect weak equivalences and hence this concludes the proof.
} 
}

\blockn{We finish this section by introducing some Reedy type categories and Reedy extension which will be very useful in the next section.}

\block{\label{CeqCI}We now consider the following directed categories:

The category $C_{eq}$ corresponding to the diagram:

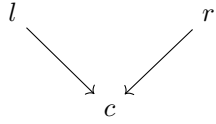
\[\begin{tikzcd}[ampersand replacement=\&]
l \arrow{rd} \&  \& r \arrow{ld}\\
\& c \& \\
\end{tikzcd} \]

The category $C_I$ corresponds to the ``co-equalizer'' diagram:
\[ 0 \rightrightarrows I \rightarrow 1 \]

We let $R_{eq}$ and $R_I$ be the Reedy type category of finite presheaves over those directed categories in which we force all the arrow of the two diagrams to be weak equivalences.

Finally if $\Ccal$ is a cylinder category we define: $\Ccal^{eq}:=[R_{eq},\Ccal]$ and $\Ccal^I := [R_I,\Ccal]$.

$\Ccal^{eq}$ is the category of triple of objects $X=(X_l,X_l,X_c)$ endowed with a cofibration $X_r \coprod X_l \hookrightarrow X_c$ whose two components are weak equivalences.

$\Ccal^I$ is a category of triple of objects $X_0,X_I,X_1$ endowed with cofibrations:

\[ X_0 \coprod X_0 \hookrightarrow X_I \]
\[ X_I \coprod_{X_0 \coprod X_0} X_0 \hookrightarrow X_1 \]

such that the two legs $X_0 \rightarrow X_I$ and the map $X_I \rightarrow X_1$ are weak equivalences.

One has morphisms $\Ccal^{eq} \rightarrow \Ccal \times \Ccal$ which send $X$ to $(X_r,X_l)$ it corresponds to the inclusion of the two minimal objects in $\Ccal^{eq}$, hence it is a fibration.

Moreover one easily see that the morphism $F_* \rightarrow R_{eq}$ corresponding to any of the minimal element is a trivial cofibration, hence the legs $\Ccal^{eq} \rightarrow \Ccal$ of the above map are trivial fibrations.

The inclusion of the minimal object in $C_{I}$ induces a Reedy extension $F_* \rightarrow R_I$, which is easily seen to be a trivial cofibration. Hence the morphism $b:\Ccal^I \rightarrow \Ccal$ which send $X$ to $X_0$ is a trivial fibration. 

Finally the functor which send $C_{eq} \rightarrow C_I$ sending $r$ and $l$ to $0$ and $c$ to $I$ induces a morphism $\Ccal^I \rightarrow  \Ccal^{eq}$ which send $(X_0,X_I,X_1)$ to $(X_0,X_0,X_I)$, hence the composite with the two morphisms $\Ccal^{eq}\rightarrow \Ccal$ are equal to $b$. The morphism $\Ccal^I \rightarrow \Ccal^{eq}$ is in particular acyclic by $2$-out-of-$3$.

}

\blockn{The intuitive idea behind those construction is that we want to construct a sort of ``path object'' for the cylinder category $\Ccal$. To do that we want to construct a category whose objects are pair of objects "with a weak equivalence" between them. This will not produce a cylinder category as such, our solution is to use instead this object $X_c$ to represent an equivalence between the two objects $X_r$ and $X_l$. The new problem is that if one want a ``identity map'' from $\Ccal$ to $\Ccal^{eq}$ on will need a functorial cylinder object, which in general is not available. The solution to this problem is the category $\Ccal^I$ which is in some sense the category of ``cylinder object''. Specifying the two maps $X \rightarrow IX$ is not enough to define a cylinder object, one also need the map $IX \rightarrow X$, adding the maps would take us outside of the realm of Reedy category, hence we replace it by this map $X_I \rightarrow X_1$ for $X_1$ an object equivalent to $X$. One hence obtain a category $\Ccal^I$ over $\Ccal$ which is equivalent to $\Ccal$ and which have a map $\Ccal^I \rightarrow \Ccal^{eq}$ playing the role of this `'ìdentity map'' we were trying to construct. These constructions will be key in the construction of the rest of the model structure in the next section: they will allows us to show that two fibrant replacement of a same cofibrant pre-cylinder categories are equivalents.}

\subsection{The weak model structure on pre-cylinder categories}
\label{SubSecWeakModelStrs}

\block{\label{Prop_CofPCylCatUnit}\Prop{Let $\Ccal$ be cofibrant pre-cylinder category, for any morphism $f:\Ccal \rightarrow \Dcal$ to a cylinder category, there is a morphism $\Ccal \rightarrow \Dcal^{eq}$ whose composition with the two maps $\Dcal^{eq} \rightarrow \Dcal$ are equal to $f$.

}

\Dem{As $\Ccal$ is cofibrant and $\Dcal^I \rightarrow \Dcal$ is a trivial fibration (see \ref{CeqCI}), the map $f$ to $\Dcal$ can be lifted to $\Dcal^I$, one can then compose to $\Dcal^I \rightarrow \Dcal^{eq}$ to obtain the desired map.}

}

\block{\label{PropUniqFibRep}\Prop{Let $\Ccal$ be a cofibrant pre-cylinder category, then any two fibrant replacement of $\Ccal$, i.e. a cylinder category $\widetilde{\Ccal}$ endowed with a trivial cofibration $\Ccal \rightarrow \widetilde{\Ccal}$ are equivalent under $\Ccal$, in the sense that there is an acyclic map between them compatible to the map from $\Ccal$.}

\Dem{ If $\Ccal_1$ and $\Ccal_2$ are two fibrant replacements of $\Ccal$ then using the lifting of trivial cofibrations against fibrant object, one obtain two maps $f,g$ under $\Ccal$ between $\Ccal_1$ and $\Ccal_2$, we will show that $g \circ f$ and $f \circ g$ are equivalence, which, because of the $2$-out-of-$6$ property, implies that $f$ and $g$ are weak equivalences between $\Ccal_1$ and $\Ccal_2$

So we need to prove that any endormorphisms $f:\Ccal_1 \rightarrow \Ccal_1$ under $\Ccal$ is a weak equivalence. Consider the (commutative) square:

\[
\begin{tikzcd}[ampersand replacement=\&]
\Ccal \arrow{r} \arrow[hook]{d}{\sim} \&  \Ccal_1^{eq} \arrow[two heads]{d} \\
\Ccal_1 \arrow{r}{(Id,f)} \& \Ccal_1 \times \Ccal_1 \\
\end{tikzcd}
\]

where the upper map is the map produced by proposition \ref{Prop_CofPCylCatUnit} using that $\Ccal$ is cofibrant. A lifting in this square produces a map $\Ccal_1 \rightarrow \Ccal_1^{eq}$ which can be thought of as ``a homotopy between $f$ and the identity''. The map $\Ccal_1 \rightarrow \Ccal_1^{eq}  $ is it self acyclic (because when composed to the projection to first componenent one get the identity,and the projection $\Ccal_1^{eq} \rightarrow \Ccal$ is acyclic, and hence the map $f$ (which is the composite of this map to $\Ccal_1^{eq}$ with the second projection $\Ccal_1^{eq}$ is an equivalence. )

}
}

\block{\label{DefExtendedAcyclic}\Def{A morphism between fibrant \emph{or} cofibrant pre-cylinder categories is said to be acyclic if its lifting to fibrant replacements of the cofibrant pre-cylinder categories involved is acyclic.}

Those definitions are consistant: they do not depends on the choice of the fibrant replacements because of proposition \ref{PropUniqFibRep} above, and if an object is both fibrant and cofibrant then it can be used as its own fibrant replacement.
Moreover any trivial cofibration from a cofibrant object is a weak equivalence as a fibrant replacement of the co-domain will also be a fibrant replacement of the domain.

Also this extended notion of acyclic morphisms still satisfies $2$-out-of-$6$ because any chain of morphisms can be lifted to a chain of morphism between some fibrant replacements (not changing the fibrant objects) and the $2$-out-of-$6$ property for the replacement implies the $2$-out-of-$6$ property for the original maps.

\bigskip

As far as we know, there is no reason for trivial cofibrations with fibrant domains to be weak equivalences (assuming their target is fibrant or cofibrant so this has a meaning) in particular one should be careful to not use fibrant replacement of an already fibrant but non cofibrant pre-cylinder category in this definition.

}

\block{\label{PropTrivCofAcyclicCof}\Prop{A cofibration of cofibrant domain is acyclic if and only it is a trivial cofibration.}

\Dem{One direction is already known, so consider an acyclic cofibration with cofibrant domain. Assume first that the co-domain is fibrant. Then one can factor our map as a trivial cofibration followed by a fibration. By $2$-out-of-$3$, the fibration part is a trivial fibration and then the lifting of cofibrations against trivial fibration gives that our cofibration is a retract of the trivial cofibration part of its factorization, hence it is a trivial cofibration.

In the general case one post compose our acyclic cofibration with a trivial cofibration whose target is fibrant. The composite is an acyclic cofibration with fibrant target hence it is a trivial cofibration. But because trivial cofibration are charaterized by the left lifting property with respect to map between fibrant object this is enough to conclude that the map is a trivial cofibration, indeed if $i$ and $i \circ u$ are trivial cofibration then for any solid square where the right arrow is a fibration between fibrant objects:

\[
\begin{tikzcd}[ampersand replacement=\&]
\Acal \arrow{r} \arrow{d}{u} \& \Xcal \arrow[two heads]{d}\\
\Bcal \arrow{r} \arrow{d}{i} \& \Ycal \\
C \arrow[dashed]{ru}[description]{l_1} \arrow[dashed]{ruu}[description]{l_2} \\
\end{tikzcd}
\]

One can construct the dashed arrow $l_1$ because $\Ycal$ is fibrant and $i$ is a trivial cofibration. Once this is done, one can construct the dashed arrow $l_2$ using the fact that the composite $i \circ u$ is a trivial cofibration, and this shows that $u$ has the lifting property and hence that it is a trivial fibration.

}

}

\block{\label{ModelStrSumUp}We now sum up the properties of this weak model structure:

\begin{itemize}

\item One knows what are weak equivalences (acyclic morphisms) only between objects that are either fibrant or cofibrant (by definition \ref{DefExtendedAcyclic}).

\item We have four class of maps: fibration, cofibrations, trivial fibration, trivial cofibration, which are all characterized by the lifting property they are supposed to have, for examples, fibrations are exactly the maps that have the right lifting properties with respect to all trivial cofibrations and cofibrations are the map which have the left lifting property with respect to all trivial fibrations (see subsection \ref{SubSecFibCofib} and in addition one defines trivial fibrations as the maps having the right lifting properties with respect to cofibrations).

\item For maps with fibrant co-domain, being a trivial fibration is the same as being an acyclic fibration. (because of proposition \ref{LiftingForAcyclicFib})

\item For maps with cofibrant domain, being a trivial cofibration or an acyclic cofibration is the same (this is proposition \ref{PropTrivCofAcyclicCof}).

\item Any morphism with a fibrant co-domain can be factored as a trivial cofibration followed by a fibration (proposition \ref{trivcof-fib-facto}) or as a cofibration followed by a trivial fibration (corollary \ref{CorFactoCofTrivFib}). 

\end{itemize}

}

\subsection{Homotopy slice of cylinder categories}
\label{SubSecSlice}

\blockn{Finally, we want to consider ``slice cylinder category''. In general, if $\Ccal$ is a cylinder category, the ordinary slice category $\Ccal_{/X}$ has no reason to be a cylinder category, so will replace this with a ``homotopy slice category'' which we will denote $\Ccal^X$.}

\blockn{Strictly speaking, this subsection has nothing to do with the weak model structure on pre-cylinder category, but it is very much in line with the technique used here (especially with the theory of Reedy extension) and it provides new examples of fibrations between cylinder categories.}

\block{\Prop{Let $\Ccal$ be a cylinder category and $X \in \Ccal$ an object. There is a cylinder category $\Ccal^X$ such that:

\begin{itemize}
\item Objects of $\Ccal^X$ are the triple $A=(A_0,A_f,A_i)$ where $A_0$ and $A_f$ are objects of $\Ccal$ and $A_i$ is a cofibration $\left( A_0 \coprod X \right) \hookrightarrow A_f$ such that the restriction to $X$ is a trivial cofibration $A \overset{\sim}{\hookrightarrow} A_f$.

\item Morphism are the obvious diagram morphism: couple $(v_0,v_f)$ , $v_0:A_0 \rightarrow B_0$ , $v_f:A_f \rightarrow B_f$ such that the naturality square with $A_i$ and $B_i$ commutes.

\item Weak equivalences are the morphisms $v$ such that $v_0$ is a weak equvalence ($v_f$ is always a weak equivalence by $2$-out-of-$3$).

\item Cofibrations are the ``Reedy cofibrations'':  $v:A \rightarrow B$ is a cofibration if $v_0:A_0 \hookrightarrow B_0$ is a cofibration and the map $ B_0 \coprod_{A_0} A_f \rightarrow B_f$ is a cofibration.

\end{itemize}

Moreover the functor $\Ccal^X \rightarrow \Ccal$ which send $A$ to $A_0$ is a fibration of cylinder categories.

}

\Dem{Let $\Ccal'$ be the category with the same description as $\Ccal^X$ but without the requirement on objects that  $X \hookrightarrow A_a$ is a weak equivalence. Then $\Ccal'$ can be described as $\Ccal^{\downarrow}( 0 \hookrightarrow X)$ in particular it is a cylinder category and the morphism $\Ccal^{\downarrow}( 0 \hookrightarrow X) \rightarrow \Ccal(0)= \Ccal$ is a fibration by proposition \ref{FibrationCoslice}.

As $\Ccal^X$ is h-saturated in $\Ccal'$ this concludes the proof.

} 
}
\blockn{More generally:}
\block{\label{PropFibSlice}\Prop{If $F: \Ccal \rightarrow \Dcal$ is a fibration between cylinder categories and $X$ is an object of $\Ccal$, the natural morphism:

\[ F^X: \Ccal^X \rightarrow \Dcal^{F(X)} \]

is a fibration.
}

\Dem{$F$ induces a morphism $\Ccal_{\downarrow} \rightarrow \Dcal_{\downarrow}$ which is a fibration because $\Ccal_{\downarrow} = [F_{\hookrightarrow},\Ccal]$ with $F_{\hookrightarrow}$ a Reedy type pre-cylinder category. Hence by proposition \ref{FibrationCoslice}, $F$ induces a fibration: 

\[ \Ccal_{\downarrow}(0 \hookrightarrow X) \rightarrow \Dcal_{\downarrow}(0 \hookrightarrow F(X))\]

$\Ccal^{X}$ and $\Dcal^{F(X)}$ are h-saturated sub-categories of those two cylinder categories and this $F$ send $\Ccal^{X}$ to $\Dcal^{F(X)}$ and hence it is a fibration. 

}

}
\block{\Def{An object $X$ in a cylinder category $\Ccal$ is said to be h-terminal if any solid diagram in $\Ccal$ of the form:

\[
\begin{tikzcd}[ampersand replacement=\&]
A \arrow[hookrightarrow]{d} \arrow{r}  \&  X \\
B \arrow[dashed]{ru} \&  \\
\end{tikzcd}
\]

admit a dashed filling.

}
Note that this is equivalent to the fact that the map $X \rightarrow 1$ is a weak equivalence in $\widetilde{\Ccal}$. In particular, any object weakly equivalent to $X$ is also h-contractible, and it is not necessary to test the definition of h-contractible on all cofibrations $A \hookrightarrow B$ and all maps $A \rightarrow X$ but only on what would be necessary to test equivalences in $\widetilde{\Ccal}$ (see subsection \ref{SubSecCompletion}). 
}

\block{\label{PropCosliceAndHTerm}\Prop{\begin{itemize}
\item The object $(X,IX,X \coprod X \hookrightarrow IX)$ of $\Ccal^X$ is h-terminal.
\item The fibration $\Ccal^X \rightarrow \Ccal$ is acyclic if and only if $X$ is $h$-terminal in $\Ccal$.
\end{itemize}
}

\Dem{

\begin{itemize}

\item Let:

\[
\begin{tikzcd}[ampersand replacement=\&]
A \arrow[hookrightarrow]{d} \arrow{r}{u}  \&  (X,IX) \\
B \&  \\
\end{tikzcd}
\]
be a diagram in $\Ccal^{X}$.

The map $A_f \hookrightarrow B_f$ is a trivial cofibration (indeed they are both equivalent to $X$) so it amit a retraction $r : B_f \hookrightarrow A_f$, and $r$ clearly preserve the inclusion of $X$.

Composing this with the map $u_f:A_f \rightarrow IX$ and $IX \rightarrow X$, one obtains a map $B_f \rightarrow X$ which extend the map $u_0:A_0 \rightarrow X$. One defines $v_0:B_0 \rightarrow X$ to be the restriction of $B_f \rightarrow X$ to $B_0$, one then has a commutative square:

\[
\begin{tikzcd}[ampersand replacement=\&]
A_f \coprod_{A_0} B_0 \arrow[hookrightarrow]{d} \arrow{r}{(u_f,v_0)}  \&  IX \arrow{d}{\sim} \\
B_f \arrow{r} \& X \\
\end{tikzcd}
\]

and a diagonal filler (making only the upper triangle commute) produces the map $v_f$ such that $(v_0,v_f)$ is the filler we are looking for.

\item If $f :\Ccal \rightarrow \Dcal$ is a acyclic morphism then it takes h-terminal object to h-terminal object. Indeed any diagram $\Dcal$ testing whether $F(X)$ is h-terminal can be replaced by an equivalent one where the cofibration is of the form $F(c)$, then one an construct a filler in $\Ccal$ using that $X$ is h-terminal and its image is a filler in $\Dcal$.

We now assume that $X$ is h-terminal. Consider a solid diagram of the form:

\[
\begin{tikzcd}[ampersand replacement=\&]
A \arrow[hookrightarrow]{d} \arrow{r}{u}  \&  Y \\
B \arrow[dashed]{ru}{v} \&  \\
\end{tikzcd}
\]

in $\Ccal^{X}$ and assume that the dashed arrow $v$ exists in $\Ccal$, i.e. that $v_0 : B_0 \rightarrow Y_0$ exists. Then, as $Y_f$ is homotopy equivalent to $X$, it is also $h$-terminal, and hence one can find a dashed arrow in:

\[
\begin{tikzcd}[ampersand replacement=\&]
A_f \coprod_{A_0} B_0 \arrow[hookrightarrow]{d} \arrow{r}{(u_f,v_0)}  \&  Y_f \\
B_f \arrow[dashed]{ru}{v_f} \& \\
\end{tikzcd}
\]

That produces our desired arrow $v$ and concludes the proof.

\end{itemize}

}

}

\section{Cylinder coherators and examples}
\label{SecCylindCoh}

\subsection{Cylinder coherators}
\label{SubSecCylindCoh}

\blockn{We can now introduced the main concept of the present paper:}

\block{\label{DefCylindCoh}\Def{A cylinder coherator $\Ccal$ is a fibrant replacement of the free pre-cylinder category on one object $F_*$. I.e it is a cylinder category with a marked object $*$, such that the morphism $F_* \rightarrow \Ccal$ corresponding to this object is a trivial cofibration.

If $\Ccal$ is a cylinder coherator a $\Ccal$-groupoid is an object of the category $\widetilde{\Ccal}$.
}
}

\blockn{It follows\footnote{The fact that the functor $\pi_{\infty}$ mentioned in the last claim is an equivalence follows from a result in the next subsection.} from the theory developed in the paper that:}

\block{\label{MainThCylindCoh}\Th{Let $\Ccal$ be a cylinder coherator, then:

\begin{enumerate}

\item The category of $\Ccal$-groupoids is ``algebraic'', i.e. it is a category of models of a quasi-equational theory (see \cite{palmgren2007partial}).

\item The category of $\Ccal$-groupoids is endowed with a semi-model structure.

\item For any two cylinder coherator $\Ccal$ and $\Dcal$ the category of $\Ccal$-groupoid and $\Dcal$-groupoid are connected by Quillen equivalences (in the two directions). 

\item The Quillen equivalences mentioned above identifies the homotopy categories of $\Ccal$-groupoids and $\Dcal$-groupoids in a canonical way.

\item There is a ``fundamental $\infty$-groupoid'' functor $\pi_{\infty} : \text{Spaces} \rightarrow \Ccal$-groupoids. It is a right Quillen equivalence, the left adjoint is called the geometric realization. This functor is unique when considered as a functor on the homotopy categories.

\end{enumerate}

}

\Dem{
\begin{enumerate}

\item The definition of $\widetilde{\Ccal}$ as the category of functors $\Ccal^{op} \rightarrow \texttt{Sets}$ commuting to certain finite limits shows that $\widetilde{\Ccal}$ is the category of models of a finite limits sketches as in \cite[Chapter 4]{barr1985toposes}, and those are equivalents (see the remarks at the very end of section 9 of \cite{palmgren2007partial})  to the various kind of theories studied in \cite{palmgren2007partial}, including quasi-equational theories and partial Horn theories.

\item It is the semi-model structure on $\widetilde{\Ccal}$ constructed in subsection \ref{SubSecCompletion}.

\item Any two cylinder coherators are connected by an acyclic morphism because of proposition \ref{PropUniqFibRep}, and such acyclic morphism induces a Quillen equivalence on the categories of groupoids because of the third point of proposition \ref{CompletionEquivalence}.

\item The functor $\Ccal \mapsto h \widetilde{\Ccal}$ send acyclic morphisms of cylinder categories to equivalences of categories, hence it induces a functor from the homotopy category of cylinder categories to the homotopy category of categories. But as $F_*$ is cofibrant, the weak model structure on pre-cylinder categories (see \ref{ModelStrSumUp}) is sufficient to prove that all cylinder coherator are canonically isomorphic in the homotopy category.

\item The category of cofibrant objects in $Spaces$ is a cofibration category, hence one can construct a morphism $F$ from $\Ccal$ to it which send $*$ to the point. $\pi_{\infty}$ is then just the restriction functor $R_F$ discussed in subsection \ref{SubSecCompletion}. The uniqueness on the functor on the homotopy categories follows from the same argument as above. In order to show that it is an equivalence it suffices to prove it for a single cylinder coherator and a single such morphism and this will be done in \ref{SimpHomotopyHypothesis}, see also \ref{EquivWithFinCW}.

\end{enumerate}

}

}

\block{We would like to conclude this section by showing that (contrary to the present situation with Grothendieck $\infty$-groupoids) one can define a $\Ccal$-groupoids of weak functors between two $\Ccal$-groupoids relatively simply:

Let $X$ be a cofibrant $\Ccal$-groupoid and $Y$ be a general $\Ccal$-groupoids, then there is a morphism $F$ of cylinder categories from $\Ccal$ to the category of cofibrant $\Ccal$-groupoid which send $* \in \Ccal$ to $X$. One can then define:

\[ [X,Y] = R_F(Y) \]

i.e. if $A$ is an object of $\Ccal$, $[X,Y](A)$ is defined as the set of morphism from $F(A)$ to $Y$.

If $X$ is not cofibrant we just replace it by a cofibrant object, following the idea (justified by the existence of our semi-model structure on $\Ccal$-groupoids) that weak functors from $X$ to $Y$ are strict functors (i.e. morphisms) from a cofibrant replacement of $X$ to $Y$.
}

\subsection{Semi-simplicial Kan complex and the homotopy hypothesis for cylinder coherators}
\label{SubSecSimplicialCoh}

\block{Let $\Delta$ denotes the semi-simplcial category whose objects are the finite inhabited ordered sets and whose map ares the injective order preserving map. $\Delta$ is a locally finite directed category. One denotes $\Ccal^{\Delta}_0$ the Reedy type pre-cylinder category of finite presheaves in this locally finite directed category (following proposition \ref{ReedyStructure}) but in which all the arrows between the representable objects are forced to be weak equivalences. }

\block{We now introduce some objects of importance in $\Ccal^{\Delta}_0$:

\begin{itemize}

\item $\Delta_n$ which is the representable object corresponding to $[n]=\{0, \dots,n \}$

\[ \Delta_n( [v] ) = \left\lbrace f: [v] \rightarrow [n] \in \Delta \right\rbrace \]

\item $\partial \Delta_n$  which is equal to $\Delta_n$ in dimension $<n$ and equal to $\emptyset$ in dimension $n$.

\[ \partial \Delta_n( [v] ) = \left\lbrace f: [v] \rightarrow [n] \in \Delta \left\vert \begin{array}{c}
\text{such that $f$ avoid at least} \\ \text{one value. }
\end{array} \right. \right\rbrace \]

\item $\Lambda_k^n$ which is defined by:

\[ \Lambda^n_{k}( [v] ) = \left\lbrace f: [v] \rightarrow [n] \in \Delta \left\vert \begin{array}{c}
\text{such that $f$ avoid at least} \\ \text{one value which is not $k$. }
\end{array} \right. \right\rbrace \]

\end{itemize}

One has tautological inclusions (hence cofibrations) :

\[ \Lambda^n_k \hookrightarrow \partial \Delta_n \hookrightarrow \Delta_n \]

The cofibration $\partial \Delta_n \hookrightarrow \Delta_n$ are generating cofibrations because of the construction done in \ref{ReedyStructure} of this sort of pre-cylinder categories and the description of the generating cofibrations of a free pre-cylinder categories given in \ref{FreePrecylinderCatMain}.

}

\block{\Prop{In $\Ccal^{\Delta}_0$ the cofibrations $\Lambda^n_k \hookrightarrow \Delta_n$ are trivial cofibrations. }

\Dem{We will show instead that the map $s_k : \Delta_0 \rightarrow \Lambda^n_k$  is a trivial cofibration, as the map $s_k \Delta_0 \hookrightarrow \Delta_n$ is a weak equivalence by definition, this will concludes the proof. Moreover this will be done by induction on $n$. For $n=1$, $\Lambda_k^1$ is $\Delta_0$ hence the result is immediate.

assume the results is proved for all $j<n$. We introduce:

\[ X_i([v]) =\{ f: [v] \rightarrow [n] | \texttt{Im } f \cup \{k\} \text{ has cardinal $\leqslant i$}  \} \]

$X_1$ is $\Delta_0$ embedded by $s_k:\Delta_0 \rightarrow \Delta_n$, the $X_i$ are increasing, and $X_{n-1}$ is $\Lambda^n_{k}$. We will show that for all $i$ the inclusion $X_i \rightarrow X_{i+1}$ is a trivial cofibration. One start with $X_i$, and let $P$ be a subset of $\{0, \dots, n\}$ which contains $k$ and has cardinal $i+1$, it corresponds to a cofibration $\Delta_i \hookrightarrow \Delta_n$ into which the maps $s_k : \Delta_0 \hookrightarrow \Delta_n$ factors. The intersection of $P$ with $X_i$ (as sub-presheaves of $\Delta_n$) is isomorphic to $\Lambda^i_{j}$ for some $j$ (corresponding to the position of $k$ in $P$) and its inclusion to $P$ is the natural map $\Lambda^i_{j} \rightarrow \Delta_i$ taking the pushout:

\[  X_i \coprod_{\Lambda^i_j} \Delta_i \]

produces a sub-object of $[n]$, where in comparison to $X_i$, $P$ and $P-\{k\}$ have been added as possible images for $f$. One can iterate this for all such subset $P$, every times it corresponds to a pushout along $\Lambda^i_j \rightarrow \Delta_i$ for some $j$ and in the end one gets $X_{i+1}$. Hence $X_i \rightarrow X_{i+1}$ is a trivial cofibration by induction, and this concludes the proof.
}
}

\block{\Def{Let $\Ccal^{\Delta}$ be the pre-cylinder category freely obtained from $\Ccal^{\Delta}_0$ by freely  adding one retraction to all the trivial cofibration $\Lambda^n_k \hookrightarrow \Delta_n$ for $n>0$.}

The natural inclusion $\Ccal^{\Delta}_0 \rightarrow \Ccal^{\Delta}$ is a trivial cofibration as it is obtained by freely adding retractions to some trivial cofibrations.

$\Ccal^{\Delta}$ can be described using the theory developed in subsection \ref{SectionFreeCofCat}. Its completion is the category of semi-simplicial set $X$, endowed with maps $X(\Lambda^n_k) \rightarrow X(\Delta_n)$ which are sections of the natural restriction map $X(\Delta_n) \rightarrow X(\Lambda_k^n)$. They somehow are a semi-simplicial version of the category of ``algebraic Kan complexes'' introduced in \cite{nikolaus2011algebraic}. The category $\Ccal^{\Delta}$ it self is the full subcategory of objects which can be obtained by freely adding cells gradually.
}

\blockn{In order to prove the homotopy hypothesis, it will be convenient to consider a very specific (and very classical) geometric realisation morphism: $\Ccal^{\Delta} \rightarrow$ finite CW-complex.}

\block{\Def{Let $R$ be the geometric realization morphism defined first on $\Ccal^{\Delta}_0$ be sending $\Delta_n$ to the ordinary $n$-simplex in the cyinder category of finite CW-complex and then extending it $\Ccal^{\Delta}$ by chosing any image of the freely added retraction of $\Lambda_k^n \hookrightarrow \Delta_n$.}}

\block{\Lem{Let $X$ be an object of $\Ccal^{\Delta}$, one can forget its additional structure and consider it as a (no longer finite) semi-simplicial set, which we will denote by $X_{ss}$. There is a natural comparison map from the geometric realization $|X_{ss}|$ of $X$ as a semi-simplicial set to $R(X)$ and it is a homotopy equivalence.  }

\Dem{
Any object $X$ in $\Ccal^{\Delta}$ is obtained by a finite sequence of iterated pushout of $\partial \Delta_n \hookrightarrow \Delta_n$. We will prove the result by induction on this cell gluing process.

The case $X= \emptyset$ is trivial. Let $X$ be any object, and assume that the map $|X_{ss}| \rightarrow R(X)$ is constructed and is a homotopy equivalence. 

Let $a : \partial \Delta_n \rightarrow X$ and let:

 \[  X^+ = X \coprod_{\partial \Delta_n} \Delta_n.\]
 
We need to understand better what is $X^+_{ss}$ and its geometric realization in order to compare $|X^+_{ss}|$ to $R(X^+)=R(X) \coprod_{\partial \Delta_n} \Delta_n$, i.e. how the colimit describing $X^+$ is computed concretely.

One can first compute the colimit $Y =X_{ss} \coprod_{\partial \Delta_n} \Delta_n$ in the category of semi-simplicial set. It boils down to simply add one $n$-dimensional cell to $X_{ss}$. The geometric realization of $Y$ as a semi-simplicial set is $|X_ss| \coprod_{\partial \Delta_n} \Delta_n$. As any colimit in a category of models of a quasi-equational theory, one then need to iteratively freely add all the missing structures, in this case it boils down to iteratively take pushout along $\Lambda_k^n \hookrightarrow \Delta_n$, in particular the geometric realization of the map of semi-simplicial set $Y \rightarrow X^+$ is a trivial cofibration of CW-complexes. Moreover this map has a retraction, either simply because it is a trivial cofibration, but there is even a canonical (and functorial) such retraction using the fact that all the $\Lambda_k^n \hookrightarrow \Delta_n$ have already a chosen retraction (fixed by the morphism $R$).

Now $R(X^+)$ is $R(X) \coprod_{\partial \Delta_n} \Delta_n$ and one has a commutative diagram:

\[
\begin{tikzcd}[ampersand replacement=\&]
\vert X_{ss} \vert \arrow[hookrightarrow]{r} \arrow{d}{\sim} \& \displaystyle \vert X_{ss} \vert \coprod_{\partial \Delta_n} \Delta_n \arrow{d}{\sim} \arrow[hookrightarrow]{r}{\sim} \& \vert (X^+)_{ss} \vert \arrow{dl}{\sim}  \\
R (X) \arrow[hookrightarrow]{r} \& \displaystyle R( X^+ ) = R( X ) \coprod_{\partial \Delta_n} \Delta_n \& \\
\end{tikzcd}
\]

hence this concludes the proof.

}

}

\block{\Lem{Let $X \in \Ccal^{\Delta}$, let $a:\partial \Delta_n \rightarrow X$ be a map in $\Ccal^{\Delta}$ and assume that there is commutative triangle:

\[
\begin{tikzcd}[ampersand replacement=\&]
 \partial \Delta_n \arrow[hookrightarrow]{d} \arrow{r}{R(a)} \& R(X) \\
\Delta^n \arrow{ru}{f} \& \\
\end{tikzcd}
\]

in the category of finite CW-complexes.

Then $f$ can be replaced by an arrow of the form $R(x)$ for $x$ an $n$-cell de $X$ whose boundary is $a$ and such that $R(x)$ is homotopy equivalent to $f$ relative to $\partial \Delta_n$.

}

\Dem{As $R(X)$ and $|X_{ss}|$ are homtopy equivalent, one can replace $R(X)$ by $|X_{ss}|$ in the above diagram: it will not change the existence or not of a diagonal filling nor the question of homotopy between two such filling.

$X_{ss}$ is a semi-simplicial Kan complex, so by \cite{rourke1970delta} it can be endowed with degeneracies making it into a simplicial set and hence a simplicial Kan complex (a more direct combinatorial proof of this can also be found in \cite{mcclure2013semisimplicial}), in particular, there exists a simplicial Kan complex $Y$ such that $X_{ss}$ is equal to $Y_{ss}$, i.e. the semi-simplicial set obtained by forgetting the face maps of $Y$. There is a natural homotopy equivalence from $|X_{ss}|$ to $|Y|$ (because $|X_{ss}|$ is the so called fat geometric realization of $Y$). Our map $a'$ and $f$ can be pushed forward to $|Y|$ were it is well known that the map $f$ will be represented (up to homotopy relative to $|\partial \Delta_n|$) by a cell of $Y$ (this follows from the Quillen equivalence between simplicial set and spaces and the fact that $Y$ is a Kan complex, hence a fibrant object). This cell of $Y$ is in particular a cell of $X_{ss}$ and a cell of $X$ (they all have the same underlying semi-simplicial set), and as the map from $|X_{ss}|$ to $|Y|$ is an homotopy equivalence, the homotopy between $f$ and $|x|$ in $|Y|$ can be lifted to $|X_{ss}|$ and then pushed to $R(X)$ and this concludes the proof.
}
}

\block{\label{SimpRealFullyFaithful}\Prop{The geometric realization from $\Ccal^{\Delta}$ to finite CW-complexes is homotopy fully faithful, more precisely: 
for every cofibration $A \hookrightarrow B$ in $\Ccal^{\Delta}$ for every map $a:A \rightarrow X$ in $\Ccal^{\Delta}$ and for every commutative triangle in the category of finite CW-complex:

\[
\begin{tikzcd}[ampersand replacement=\&]
R(A) \arrow[hookrightarrow]{d} \arrow{r}{|a|} \& R(X) \\
R(B) \arrow{ru}{f} \& \\
\end{tikzcd}
\]

there exists a map $f':B \rightarrow X$ making the triangle commutes in $\Ccal^{\Delta}$ and such that $|f'|$ is homotopy equivalence to $|f|$ relative to $|A|$.
}

\Dem{This is an immediate induction using the above lemma and the fact that any cofibration inf $\Ccal^{\Delta}$ is an iterated pushout of maps $\partial \Delta_n \rightarrow \Delta_n$.}

}

\block{\Prop{$\Ccal^{\Delta}$ is a cylinder category.}

\Dem{We will first prove that it admit cylinder objects. Any object is obtained by iterated gluing of boundary inclusion: $\partial \Delta_n \hookrightarrow \Delta_n$, so we will work by induction on those gluing. The initial object always has a cylinder object (the initial object).

Assume that $X \in \Ccal^{\Delta}$ has a cylinder object $IX$. Let $X^+ = X \coprod_{\partial \Delta_n} \Delta_n$ for some integer $n$ and some map $a:\partial \Delta_n \rightarrow \Delta_n$. 

Let:

 \[ Y = IX \coprod_{X \coprod X} \left( X^+ \coprod X^+ \right) = IX \coprod_{\partial \Delta_n \coprod \partial \Delta_n } \left( \Delta_n \coprod \Delta_n \right) \]

There is a cofibration $X ^+ \coprod X^+ \rightarrow Y$ and a map $Y \rightarrow X^+$ which send $IX$ to $X$ by the natural map and the $X^+$ to $X^+$ by the identity. 

Let also:

\[ Z =  IX \coprod_{X \coprod X} \left( X \coprod X^+ \right) = IX \coprod_{\partial \Delta_n } \Delta_n \]

where we only add a new cell to the second component. The natural map $Z \rightarrow Y$ is a cofibration and the map $Z \rightarrow X^+$ is a weak equivalence (it is a pushout of the weak equivalence $IX \rightarrow X$ along $\partial \Delta_n \hookrightarrow \Delta_n $).

The map $i_0 \circ a : \partial \Delta_n \rightarrow Z$, i.e. the one which has not been filled yet, actually already have a filling: indeed it has one when it is composed with the weak equivalence $Z \rightarrow X^+$, hence it has one at the level of geometric realization and hence it has one by proposition \ref{SimpRealFullyFaithful}. Let $b : \Delta_n \rightarrow Z$ be this map.

Then one has a map $\Delta_n \coprod_{\partial \Delta_n} \Delta_n \rightarrow Y$ which send $\partial \Delta_n$ to $i_0 \circ a$, the first $\Delta_n$ to the cell freely added in $Y$ and the second $\Delta_n$ is send into $Y$ though the map $b$ defined above. Then let:

\[ IX^+ = Y \coprod_{\Delta_n \coprod_{\partial \Delta_n} \Delta_n} I_{\partial \Delta_n} \Delta_n = Z \coprod_{\Delta_n} I_{\partial \Delta_n} \Delta_n\]

$Z \rightarrow IX^+$ is a weak equivalence because it is a pushout of a tribial cofibration, and there is a natural map $IX^+ \rightarrow X^+$ which is hence a weak equivalence by $2$-out-of-$3$. The map $X^+ \coprod X^+ \hookrightarrow Y \hookrightarrow IX^+$ is a cofibration hence this concludes the proof that $X^+$ admit a cylinder object and hence that every object admit a cylinder object.

It remains to show that trivial cofibration admit retraction, but this follow from proposition \ref{SimpRealFullyFaithful} and the fact that trivial cofibration are cofibration and admit retraction on their geometric realization.
}
}

\block{\label{SimpHomotopyHypothesis}\Th{The geometric realization functor from $\Ccal^{\Delta}$ to finite CW-complexes is acyclic.}

\Dem{We have already shown that these two categories are cylinder category and that the functor is homotopy fully faithful, hence we only have to show that it is homotopy full i.e. that each finite $CW$-complexes is homotopy equivalent to something in the image of the geometric realization. This is immediately proved by cellular induction in the category of CW-complexes: it is obvious for the empty CW-complex. If a CW-complex $X$ is homotopy equivalent to $|Y|$, and $X^+$ is obtained from $X$ by gluing one extra cell then $X^+$ is homotopy equivalent to something of the form $|Y|$ with one extra cell glued. This cell is glued along a map $S^n \rightarrow |Y|$, but $S^n$ can be seen as the geometric realization of $\partial \Delta_n$ hence by proposition \ref{SimpRealFullyFaithful} this map is homotopy equivalent to the geometric realization of a map $\partial \Delta_n \rightarrow Y$ hence $X^+$ is homotopy equivalent to the geometric realization of $Y^+ = Y \coprod_{\partial \Delta_n} \Delta_n$.

}
}

\block{\Prop{The morphism $F_* \rightarrow \Ccal^{\Delta}$ and $F_* \rightarrow \Ccal^{\Delta}_0$ corresponding to the object $\Delta_0$ are trivial cofibrations.}

\Dem{It is enough to prove it for $F_* \rightarrow \Ccal^{\Delta}_0$, as the map $\Ccal^{\Delta}_0 \rightarrow \Ccal^{\Delta}$ is clearly a trivial cofibration.

Let $F :X \rightarrow Y$ be a fibration between cylinder categories, let $h: \Ccal^{\Delta}_0 \rightarrow Y$ and let $h_0$ be a lifting of $\Delta_0$ to $X$, one need to extend this to a lifting of $h$. As $Y$ is a cylinder category, $h$ can be extended to a functor from $\Ccal^{\Delta}$ to $Y$. Now using one can easily construct in $\Ccal^{\Delta}$ a maps $c_n:\Delta_n \rightarrow \Delta_0$ whose restriction to the face are the map $C_k : \Delta_k \rightarrow \Delta_0$: this is done by induction on $n$ by applying \ref{SimpRealFullyFaithful} to the cofibration $\partial \Delta_n \hookrightarrow \Delta_n$. In particular, $\partial \Delta_n \rightarrow \Delta_n \rightarrow \Delta_0$ is a factorization as cofibration followed by a weak equivalence. Such a factorization can always be lifted trough a trivial fibration, hence one can lift $h$ to $X$ going though the induction procedure that construct the pre-cylinder category attached to a a locally finite directed category: for each $n$, we will be able to lift the cofibration $\partial \Delta_n \hookrightarrow \Delta_n$, hence producing a lifting of $h$ at least as a morphism $\widehat{\Delta}_f \rightarrow X $. Moreover as at every step the lifting comes with a weak equivalence $\Delta_n \overset{\sim}{\rightarrow} \Delta_0$ it ensure that the newly lifted object $\Delta_n$ will indeed be homotopicaly equivalent to the previously lifted object and hence that it is indeed a morphism from $\Ccal^{\Delta}_0$.
}
}

\subsection{Grothendieck weak $\infty$-grouppoids}
\label{SubSecGrothGroupoid}

\blockn{We refer to \cite{ara2013homotopy} for the definition and all results about Grothendieck $\infty$-groupoid. We fix a coherator $C$ (see \cite[2.11]{ara2013homotopy}) and when we say $\infty$-groupoid is means $\infty$-groupoid of type $C$.}

\block{We denote by $D_n$ the free $\infty$-groupoid generated by one cell in dimension $n$ and by $S_n$ the free $\infty$-groupoid generated by two parallel cell of dimension $n$. $S_{-1}$ will be the initial object (the empty groupoid) $\emptyset$ by convention.

What we call the natural map $S_{n-1} \rightarrow D_n$ is the one that send the two generators of $S_{n-1}$ to the source and the target of the generator of $D_n$.
}

\block{\Def{One say that an $\infty$-groupoid is finitely cellular if it is obtained from the empty groupoid by freely adding cells, i.e. by iterated pushout of the maps $S_{n-1} \rightarrow D_n$ for $n \geqslant 0$. }}

\block{\label{conjecture}\Conjecture{Let $X$ be a finitely cellular $\infty$-groupoid, and let $a$ be an $n$-cell of $X$. Consider $X^+$ the (finitely cellular) $\infty$-groupoid obtained by freely adding a cell $a'$ parallel to $a$ and a cell $i$ between $a$ and $a'$. i.e. $X^+$ is the pushout:

\[ X^+ := X \coprod_{D_n} D_{n+1} \]

Then $X \rightarrow X^+$ is a weak equivalence.
}

``Weak equivalences'' is taken in the sense of \cite[4.17 and 4.18]{ara2013homotopy}.

}

\blockn{In all the rest of the subsection we investigate consequences of this conjecture.}

\block{\Lem{If conjecture \ref{conjecture} holds, then it also holds when $X$ is a retract of an arbitrary cellular object, i.e. is a retract of an object obtained from $\emptyset$ as a transfinite composition of iterated pushout of $S_n \rightarrow D_{n+1}$. }

\Dem{This follow easily from the fact that directed limits in the category of $\infty$-groupoids are computed level wise, and that $\pi$-groups commute to those directed limits.}
}

\block{\label{ModelCatForWeakGpd}\Th{Under conjecture \ref{conjecture} there is a ``canonical'' semi-model structure on the category of $\infty$-groupoid such that:

\begin{itemize}

\item Weak equivalences are the weak equivalences of \cite[4.17 and 4.18]{ara2013homotopy}.

\item Fibrations are the map with the right lifting property with respect to $s:D_n \rightarrow D_{n+1}$.

\item Trivial fibrations are the maps with the right lifting property with respect to the natural maps $S_n \rightarrow D_{n+1}$.

\end{itemize}

}

Note that as the maps $D_n \rightarrow D_{n+1}$ admit retraction, every object is fibrant.

\Dem{By the small object argument one has two factorization systems in cofibration/trivial fibration and trivial cofibration/fibration with the fibrations and trivial fibrations defined above.

The characterization of \cite[4.18.(4)]{ara2013homotopy} of weak equivalences can be reformulated as the fact that every square of the form:

\[
\begin{tikzcd}[ampersand replacement=\&]
S_{n-1} \arrow{d} \arrow{r}{(u,v)} \&  X \arrow{d}{f} \\
D_n \arrow{r} \& Y  \\
\end{tikzcd}
\]

admit a filing of the form:

\[
\begin{tikzcd}[ampersand replacement=\&]
S_{n-1} \arrow[hookrightarrow]{ddd}{i} \arrow[hookrightarrow]{dr}{i} \arrow{rr} \& \&  X \arrow{ddd} \\
\& D_n \arrow[hookrightarrow]{d} \arrow{ur}{a} \& \\
\& D_{n+1} \arrow{dr}{h} \& \\
D_n \arrow{rr} \arrow[hookrightarrow]{ur}\& \& Y \\
\end{tikzcd}
\]

hence the exact same proof as in proposition \ref{Prop_trivFibInComp} show that trivial fibrations are exactly the fibrations that are weak equivalences. 

The conjecture, and its extension by the lemma above shows that transfinite iterated pushout of the $D_n \rightarrow D_{n+1}$ out of a cofibrant objects are weak equivalences, hence by the small object argument any map out of a cofibrant object can be factored as ``trivial cofibration'' which is a weak equivalence followed by a fibration. As weak equivalences of $\infty$-groupoid satisfies $2$-out-of-$3$ this concludes the proof.
}
}

\block{\Prop{Assuming conjecture \ref{conjecture}. The category $\fgpd$ of finitely cellular $\infty$-groupoids is a cylinder category with the $S_{n-1} \hookrightarrow D_n$ as generating cofibrations and the weak equivalences of $\infty$-groupoids as weak equivalences.}

\Dem{The existence of retraction to trivial cofibrations follow from the fact that those will be trivial cofibrations in the  canonical model structure on $\infty$-groupoids and hence will admit a retraction as every object is fibrant. The only non trivial things is the existence of cylinder object. One will need the following fact easily proved by induction: if $X \hookrightarrow Y$ is an iterated pushout of some finite set of cofibration $A_i \hookrightarrow B_i$ in a cylinder category then one can chose a relative cylinder object $I_X Y$ such that $Y \hookrightarrow I_X Y$ is an iterated pushout of the trivial cofibrations $B_i \hookrightarrow I_{A_i} B_i$.
As in the cylinder category of cofibrant $\infty$-groupoids the cofibration $S_{n-1} \hookrightarrow D_n$ have relative cylinder object which are themselves finitely cellular ($D_{n+1}$) this implies that any finitely cellular object admit a finitely cellular cylinder object. This concludes the proof.
}
}

\block{\label{fgpdIsFree}\Prop{(Assuming conjecture \ref{conjecture}) The functor $F_* \rightarrow \fgpd$ corresponding to the object $D_0$ is cofibrant and the completion $\widetilde{\fgpd}$ is the category of $\infty$-groupoids with the semi-model structure of \ref{ModelCatForWeakGpd}.}

\Dem{Consider the following free construction of pre-cylinder category. One start with the globular category:

 \[ D_0 \rightrightarrows D_1 \rightrightarrows \dots \rightrightarrows D_n \rightrightarrows \dots \]

It is a directed locally finite category, hence one get a Reedy type pre-cylinder category of finite globular sets with monomorphism as cofibrations and we ask that the structural map between the $D_i$ are weak equivalence. This pre-cylinder category $\fglob$ is endowed with a morphism $F_* \rightarrow \fglob$ which is cofibrant, and the completion of $\fglob$ is the category of globular set. The map ``$S_{n-1} \rightarrow D_n$'' (the free globular set respectively generated by two parallel $n-1$ cells and one $n$-cell) are generating cofibrations.

One then follow the construction of the coherator $C$ we have fixed in the beginning, and for each lifting of a parallel pair that is freely added during the construction of $C$ one free add the same arrow freely to our pre-cylinder category. The completion of the resulting category $\Ccal$ is the category of $\infty$-groupoid and as we have not added any cofibration the pre-cylinder category we have generated is (at least as a cofibration category) $\fgpd$. The weak equivalences in $\Ccal$ are generated by the $D_n \hookrightarrow D_{n+1}$ (it is the smallest class that contains those maps and satisfies the axioms for pre-cylinder categories), in particular, weak equivalences in $\Ccal$ are weak equivalences in $\fgpd$, but as $\fgpd$ is a cylinder categories any trivial cofibrations in $\fgpd$ is a retract of an iterated pushout of the $D_n \hookrightarrow D_{n+1}$ (which are the legs of the relative cylinder object of the generating cofibrations) and hence $\Ccal$ is the same as $\fgpd$.

In particular $F_* \rightarrow \fgpd$ is cofibrant.

Moreover, the description given in proposition \ref{PropModelStrFromGenCof} of the semi-model structure on the completion in terms of the generating cofibrations shows that the model structure on the completion is exactly the one of proposition \ref{ModelCatForWeakGpd}.
}
}

\block{\label{DefDGPDwithoutConj}If we do not assume conjecture \ref{conjecture}, the proof of the above proposition still construct a pre-cylinder category $\Ccal$ that we will denote $\fgpd$ which is defined by this free construction. It is still true that it is the category of finitely cellular $\infty$-groupoid and that the $S_n \rightarrow D_{n+1}$ are generating cofibrations, but if we do not assume the conjecture it might not be a cylinder categories and its weak equivalences can be different from the ones defined in \cite[4.17 and 4.18]{ara2013homotopy}, in fact any of these two facts would implies the conjecture. This will be our definition of $\fgpd$ if we do not assume conjecture \ref{conjecture}. }

\block{\label{LeFGPDLiftingContractible}\Lem{Let $f:X\rightarrow Y$ be a fibration, let $o$ be a h-terminal object in $X$ such that $f(o)$ is h-terminal. Assume that there is a morphism $h:\fgpd \rightarrow Y$ such that $h(D_0)=f(o)$ then one can lift $h$ to a morphism to $X$.}

This lemma holds without assuming conjecture \ref{conjecture} as long as $\fgpd$ is defined using the free construction in the proof of proposition \ref{fgpdIsFree} as explained in \ref{DefDGPDwithoutConj}.

\Dem{We follow the "free" description of $\fgpd$ given in the proof of proposition \ref{fgpdIsFree} to construct the lifting. The first step is to lift the co-globular object. In $\fgpd$ one has a cofibration/weak equivalences:

\[ D_n \coprod_{S_{n-1}} D_n \hookrightarrow D_{n+1} \overset{\sim}{\rightarrow} D_{n} \]

Lifting this factorization gradually produces a co-globular object in $X$ whose zero object is $o$, which is a lifting of the one in $Y$ and whose structural maps are all weak equivalences. We now need to gradually lift the map corresponding to the liftings of parallel pairs freely added in the construction of the coherator.

Each such pair has the form $S_n \rightarrow S$ where $S$ is an $\infty$-groupoid obtained as a globular sum of the $D_n$ (see \cite[1.3]{ara2013homotopy}), and a lifting is an extension of $S_n \hookrightarrow D_{n+1} \rightarrow S$. The map $S_n \rightarrow G$ is already lifted to $X$ and a lifting exists in $Y$, one need to find a lifting in $X$ which is over the lifting in $Y$. One easily see that the objet $S$ is weakly equivalent to $o$ hence it is itself h-terminal hence such a lifting $S_n \hookrightarrow D_{n+1} \rightarrow S$ automatically exists in $X$, its image in $Y$ might not be equal to the lifting already existing but it is going to be homotopy equivalent relative to $S_n$ because $f(S)$ is also $h$-terminal. Hence one can apply the last point proposition \ref{MoreLiftingForPreFib} to find a lift in $X$ which is sent to the correct lift in $Y$ and this concludes the proof.
}
}

\block{The lemma above is enough to show (even without assuming conjecture \ref{conjecture}) that any object in any path category can be ``endowed with a structure of Grothendieck $\infty$-groupoid'', or more precisely, the globular object given by any iterated relative path object can be endowed with the operation making it into a Grothendieck $\infty$-groupoid.

Indeed, dually it corresponds to the statement that for any object $X$ in a cylinder categories $\Ccal$, there is a morphism $\fgpd \rightarrow \Ccal$ which send $D_0$ to $X$, and this is proved by applying the lemma the the fibration $\Ccal^{X} \rightarrow 0$. Indeed $X$ is $h$-contractible in $\Ccal^{X}$ because of proposition \ref{PropCosliceAndHTerm} hence one can lift the unique morphism $\fgpd \rightarrow 0$ to a morphism $\fgpd \rightarrow \Ccal^{X}$, composing with the fibration $\Ccal^{X} \rightarrow \Ccal$ produces the desired morphisms.

\bigskip

We have not been able to extend this argument to prove that $F_* \rightarrow \fgpd$ is a trivial cofibration without assuming conjecture \ref{conjecture}. Indeed it is not clear that any morphism $\fgpd \rightarrow \Ccal$ is obtained by this technique (and hence can be lifted by this technique). Fortunately, if one assume conjecture \ref{conjecture} this became an immediate consequence of the general machinery developed in the paper, and it is the lemma below.
}

\block{\label{LemAllActionsOfFGPDAregood}\Lem{(assuming conjecture \ref{conjecture}) Let $f : \fgpd \rightarrow \Ccal$ be a morphism into a cylinder category, then $f$ can be lifted to $\Ccal_{/f(D_0)} \rightarrow \Ccal$. The lifting can be chosen to send $D_0$ to any h-terminal object with the correct image in $\Ccal$. }

\Dem{In the category of $\infty$-groupoid the map from $D_0$ to the terminal groupoid $*$ is a trivial fibration. Hence $D_0$ is h-terminal in $\fgpd$, hence by proposition \ref{PropCosliceAndHTerm} the map $\fgpd_{/D_0} \rightarrow \fgpd$ is a trivial cofibration, and as $\fgpd$ is cofibrant one can construct a morphism $\fgpd \rightarrow \fgpd_{/D_0}$ which send $D_0$ to the object $(D_0,D_1)$ and whose composite with the natural fibration $\fgpd_{/D_0}$ is the identity. Composing this to the natural map $\fgpd_{/D_0} \rightarrow \Ccal_{/f(D_0)}$ gives us the desired lifting, and the image of $D_0$ is $f(D_0),f(D_1)$ which is a h-terminal object. }
}

\block{\Cor{(assuming conjecture \ref{conjecture}) The morphism $F_* \rightarrow \fgpd$ corresponding to $D_0$ is a trivial cofibration, hence $\fgpd$ is a cylinder coherator.}

\Dem{ Consider a square:

 \[
\begin{tikzcd}[ampersand replacement=\&]
F_* \arrow{d} \arrow{r}{o} \& X \arrow[two heads]{d}{f} \\
\fgpd \arrow{r}{h} \& Y  \\
\end{tikzcd}
\]
with $f$ a fibration between fibrant object. $h$ send $D_0$ to $f(o)$, one can hence construct using lemma \ref{LemAllActionsOfFGPDAregood} a first factorization:

 \[
\begin{tikzcd}[ampersand replacement=\&]
F_* \arrow{d} \arrow{rr}{o}\& \& X \arrow[two heads]{d}{f} \\
\fgpd \arrow{r}{h'} \& Y_{/f(o)}\arrow{r} \& Y  \\
\end{tikzcd}
\]

the image of $h'$ being of the form $f(o),If(o)$ for some cylinder object for $f(o)$. As $f$ is a fibration, this cylinder object can be lifted to $X$ and one can hence extend this into a factorization:

 \[
\begin{tikzcd}[ampersand replacement=\&]
F_* \arrow{d} \arrow{r}{(o,Io)}\& X_{/o} \arrow[two heads]{d} \arrow{r} \& X \arrow[two heads]{d}{f} \\
\fgpd \arrow{r}{h'} \& Y_{/f(o)}\arrow{r} \& Y  \\
\end{tikzcd}
\]

The new vertical arrow is a fibration because of \ref{PropFibSlice}. One can hence construct our lifting in the left square by lemma \ref{LeFGPDLiftingContractible} and this concludes the proof.
}

}

\block{\Cor{Conjecture \ref{conjecture} implies Grothendieck homotopy hypothesis (i.e. conjecture $2.8$ of \cite{maltsiniotis2010grothendieck}). }

\Dem{This follows from the fact that it is a cylinder coherator and theorem \ref{MainThCylindCoh}.}

}

\subsection{The opposite of the syntactic categories of (weak) intentional type theory}
\label{SubSecSyntactic}

\blockn{The main goal of this last subsection is to reconnect to theory developed here with type theory. For a reader not interested in type theory, it still contains an interesting result: there exists a notion of weak globular $\infty$-groupoid (i.e. a globular set with certain ``composition operations'' defined between finite projective limit for the set of $n$-cells) which fits into our framework and hence which satisfies the homotopy hypothesis. If the reader is only interested in this results he can take the discussion of \ref{MorphismOutOfsyntax} as a definition of the cylinder coherator $\Ccal(\T^*_{w.i.e.})^{op}$ as a free pre-cylinder category and avoid almost all mention of type theory.}

\blockn{Giving a precise and general enough definition of what is a type theory is a difficult task and subject to debate that we will not try to answer here. We will just give some examples of things that we will call ``type theory'' which defines a syntactical category and which are of the vague form described just below.}

\block{A type theory will be a formal system (i.e.  a set of rules) for deriving statement of the form

\[ \Gamma \vdash \Jcal \]

Where $\Jcal$ is a judgement and $\Gamma$ is a context.

A judgement can have one of the following form:

\[ A \texttt{ Type} \text{ or } x : A \]

Which are supposed to be interpreted respectively as the fact that $A$ is a well formed expression denoting a type and that $x$ is a well formed expression denoting an inhabitant of the type $A$.}

\block{Contrary to most reference (for example \cite{berg2016path}), one does not need equality of types or of variables to be judgements.
This is because all the computation rules of the type theory that we will consider will holds up to propositional equality (and not definitional equality), hence these types theory would have no rules proving definitional equality of types or variables aside from the obvious structural rules, and hence can only proves definitional equality when there is actually  ``syntactical'' equality of the expression. For this reason it is considerably simpler to remove judgemental equality all together and only consider propositional equality and ``syntactical equality''.
}

\block{A ``context'' is a device whose purpose is to make explicit the type of all the free variables that appears in a judgement $\Jcal$. More precisely it is a finite sequence:

\[ \left( x_0 \in A_0 ; x_1 \in A_1; \dots ; x_n \in A_n \right) \]

Such that the $x_i$ are free variables, no free variables appears in $A_0$ and the only free variables that appears in $A_i$ are $x_0, \dots, x_{i-1}$. A context is said to be valid (with respect to a type theory) if the following judgement can be inferred from the rules of the type theory:

\[ \vdash A_0 \texttt{ Type} \]
\[ x_0 \in A_0 \vdash A_1 \texttt{ Type} \]
\[ \dots \]
\[ x_0 \in A_0; x_1 \in A_1; \dots; x_{n-1} \in A_{n-1} \vdash A_n \texttt{ Type} \]
 
}

\block{A rules for a type theory will be of the form:

\[\begin{array}{c c c}
\Gamma_1 \vdash \Jcal_1 & \dots & \Gamma_n \vdash \Jcal_n \\
\hline
\multicolumn{3}{c}{\Gamma \vdash \Jcal}
\end{array}
\]

Such a rule means that if all the $\Gamma_i \vdash \Jcal_i$ can be inferred from the type theory, one can apply this rules to deduces that $\Gamma \vdash \Jcal$ can also be inferred. There should be some of the rules (called axiom) for which $n=0$, and one says that a judgement can be inferred from a type theory if it can be obtained this way by a finite number of applications of the rules (and axioms).

}

\block{We start by what we will call the ``syntactical rules'' which are the rules that in some sense ensure that the judgements and context behave in the way their intuitive interpretation suggest:

\begin{itemize}
\item[(S1)] - (Axiom): if $a : A$ appears in the context $\Gamma$ then:

\[\begin{array}{c c c}
 & & \\
\hline
\multicolumn{3}{c}{\Gamma \vdash a : A}
\end{array}
\]

\item[(S2)] - (Weakening):  

\[\begin{array}{c c c}
 \Gamma, \Delta \vdash \Jcal & & \Gamma \vdash A \texttt{ type} \\
\hline
\multicolumn{3}{c}{\Gamma,a:A,\Delta \vdash \Jcal }
\end{array}
\]

\item[(S3)] - (Substitution):

\[\begin{array}{c c c}
 \Gamma,x:A, \Delta \vdash \Jcal & & \Gamma \vdash  a: A  \\
\hline
\multicolumn{3}{c}{\Gamma,\Delta[a/x] \vdash \Jcal[a/x] }
\end{array}
\]

where $T[a/x]$ mean that all the occurrence of the free variable $x$ in $T$ have been replaced by $a$.

\end{itemize}

Those rules does not allow to prove anything interesting by themselves, but they ensure a good behaviour of the type theory and allow to construct its ``syntactic category''. 

}

\block{\label{Lem_Morph_Subs}If $\T$ is a type theory with at least the rules $(S1),(S2)$ and $(S3)$. And if $\Gamma$ and $\Delta$ are two valid context within a type theory $\T$, then a context morphism $t$ from $\Gamma$ to $\Delta$ with $\Delta=(x_1: A_1, \dots, x_n : A_n)$ is a collection of terms $t_1, \dots, t_n$ such that:

\[ \Gamma \vdash t_1 : A_1 \]
\[ \Gamma \vdash t_2 :A_2[t_1/x_1]\]
\[ \dots \]
\[ \Gamma \vdash t_n :A_n[t_1/x_1, \dots, t_{n-1}/x_{n-1}]\]

The following lemma show in particular that this definition is meaningful, i.e. that: 

\[ \Gamma \vdash A_i[t_1/x_1, \dots, t_{i-1}/x_{i-1}] \text{ Type} \]

\Lem{If $t:\Gamma \rightarrow \Delta$ is a context morphism and $\Delta \vdash \Jcal$ then $\Gamma \vdash \Jcal[t]$ where $\Jcal[t]$ denotes $\Jcal[t_1/x_1][t_2/x_2]\dots[t_n/x_n]$.}

\Dem{By iterated application of the weakening axioms $(S2)$ one has that $\Gamma,\Delta \vdash \Jcal$ then by iterated application of the substitution axiom one can gradually remove the variable $x_i$ (starting with $x_1$ and going up) for the context $\Delta$ and replacing it every where by $t_i$ to obtain the result once they have all been replaced.}

This lemma also allows to define composition of context morphisms $f \circ g$ with $f=(f_1, \dots, f_n) $ as $(f_1[g], \dots, f_n[g])$. And one easily see that $K[f][g]=K[f \circ g ]$ and that the composition is associative. This defines a category $\Ccal(\T)$ called the syntactic category of $\T$ whose object are the context and whose morphisms are the context morphisms.
}

\block{\Prop{If $\T$ is a type theory with at least the rules $(S1),(S2)$ and $(S3)$. Then its syntactic category $\Ccal(\T)$ is a fibration category (the opposite of a cofibration category) with fibrations (cofibration in the opposite category) being the map that are isomorphic to a morphism $ p:\Gamma, \Delta \rightarrow \Gamma$ which is of the form $(x_1,\dots,x_n)$ for $\Gamma=(x_1 : A_1, \dots, x_n : A_n)$
}
\Dem{Those morphisms clearly contains isomorphisms. We will first show that they are stable under pullback, for this it is enough to show that the ``dependant projection'', i.e. the map of the form $ p:\Gamma, \Delta \rightarrow \Gamma$ are stable under pullback.

Let $f :\Gamma' \rightarrow \Gamma$ be a context morphism, then because of the lemma \ref{Lem_Morph_Subs} $\Gamma',\Delta[f]$ is a valid context which fits into a square:

\[
\begin{tikzcd}[ampersand replacement=\&]
\left( \Gamma',\Delta[f] \right) \arrow[two heads]{d}{i} \arrow{r} \& \left( \Gamma,\Delta \right) \arrow[two heads]{d} \\
\Gamma' \arrow{r}{f} \&  \Gamma \\
\end{tikzcd}
\]

where the upper arrow if $f$ on the component corresponding to $\Gamma'$ and is ``the identity'' (at least syntactically) on the variable corresponding to $\Delta[f]$. If one has a morphism $\Gamma'' \rightarrow \Gamma',\Delta[f]$ then the first componenents should be a morphism $\Gamma'' \rightarrow \Gamma'$ and the remaining component should be an extenion of the composition with $f$ to $\Gamma,\Delta$, hence $\Gamma',\Delta[f]$ is indeed the pullback we are looking for, this proves that pullback of a dependant projection exists and can be represented by dependant projection, hence pullback of fibrations exist and are fibrations. Another consequence of this ``pullback stability'' of dependant projections is that if $p$ is a dependant projection and $i$ is an isomorphism then $p \circ i$ is equal to $i' \circ p'$ for $i'$ an isomorhisms and $p'$ a dependant projection hence fibrations are closed under compositions.

Finally, the empty context is a terminal object and all the map to the empty context are dependant projection, hence fibrations.
}

}

\block{We now introduce the three examples of type theory we are concerned with:

The theory $\T^*_{Id}$ has the three rules $(S1), (S2),(S3)$ and in addition:

\begin{itemize}
\item[(T1)] (Free object)

\[\begin{array}{c c c}
& &   \\
\hline
\multicolumn{3}{c}{\vdash * \texttt{ Type}}
\end{array}
\]

\item[(T2)] (Identity type)

\[\begin{array}{c c c}
   \Gamma \vdash A \text{ Type} & \Gamma \vdash a: a & \Gamma \vdash b: A \\
\hline
\multicolumn{3}{c}{\Gamma, \vdash Id_A(a,b) \texttt{ Type}}
\end{array}
\]

\end{itemize}

The first axioms is what makes all those type theory ``freely generated on one object $*$''.

The type theory $\T^*_{refl}$ has all the rules of $\T^*_{Id}$ with in addition the rule:

\begin{itemize}

\item[(T3)] (Identity intorduction) or (Reflexivity):

\[\begin{array}{c c c}
  \Gamma \vdash A \text{ Type} & & \Gamma \vdash a:A   \\
\hline
\multicolumn{3}{c}{\Gamma \vdash r(a) : Id_A(a,a) }
\end{array}
\]

\end{itemize}

Finally, the theory\footnote{$w.i.e.$ stands for ``weak identity elimination''.} $\T^*_{w.i.e.}$ has all the rules above, with in addition following two rules:

\begin{itemize}

\item[(T4)] (identity elimination)

\[\begin{array}{c}
 \Gamma,a:A,b:A,\gamma:Id(a,b),\Delta \vdash C \text{ Type}   \\
 \Gamma, a :A,\Delta[a/b,r(a)/\gamma] \vdash d :C[a/b,r(a)/\gamma] \\
\hline
\Gamma,a:A,b:A,\gamma:Id(a,b),\Delta  \vdash \J(a,b,\gamma) : C 
\end{array}
\]

\item[(T5)] (weak identity computation)

\[\begin{array}{c}
 \Gamma,a:A,b:A,\gamma:Id(a,b),\Delta \vdash C \text{ Type}   \\
 \Gamma, a :A,\Delta[a/b,r(a)/\gamma] \vdash d :C[a/b,r(a)/\gamma] \\
\hline
\Gamma,a:A,\Delta[a/b,r(a)/\gamma]  \vdash \Hb(a) :  Id(\J(a,a,r(a)), d ) 
\end{array}
\]

\end{itemize} 

A few things need to be clarified in this last two axioms about what are $\J$ and $\Hb$ syntactical. Our point is that we want to make those terms ``as different as possible'' from one others (i.e. minimalize the syntactical equality between $\J$ terms and $\Hb$ terms). They will both have free varibables: $\J$ has all the variables appearing in $\Gamma,a:A,b:A,\gamma:Id(a,b),\Delta$ as free variables and $\Hb$ has all the variables appearing in $\Gamma,a:A,\Delta[a/b,r(a)/\gamma]$ as free variables. They also have a subscript to distinguishes them which contains all the assumption of the rules that produces them (one consider two such subscript equal if one can go from one set of assumption to another by just renaming the free variables).

It is probable that those conventions can be modified to have more equality between $\J$ terms and $\Hb$ terms but one need to be very careful here as it can change the fact that the obtained syntactical category will or will not be a trivially cofibrant extension of $F_*$.

}

\block{The idea of having ``weak identity type'' as in our theory $\T^*_{w.i.e}$ where the computation rules is a propositional identity instead of a definitional identity comes from B. van den Berg in \cite{berg2016path}.}

\block{Let us gives some examples of context that one can form in those different type theories.

In $\T^*_{Id}$ (and hence also in all the other type theory mentioned) $\emptyset$ is a context, one has $\emptyset \vdash * \texttt{ Type}$ by the rule $(T1)$, hence $(x : *)$ is also a valid context. One can iterate this and $(x :*; y:*)$ and $(x:*,y:*,z:*)$ and so one are all contexts.

\bigskip

By rule $(S1)$, one has $x:*,y:* \vdash x:*$ and $x:*,y:* \vdash y:*$ hence, by rule $(T2)$ one has $x:*,y:* \vdash Id(x,y) \texttt{ Type}$ and hence $(x:*,y:*,t:Id_*(x,y) )$ is again a valid context. Similarly:

 \[ (x:*,t:Id_*(x,x))\] \[ (x:*,y:*,z:*, t:Id_*(x,y) ,w :Id_*(x,z) ) \] \[ (x : * ,y:*,t:Id_*(x,y), s: Id_*(x,y))\]
 
are valid context. Going one level higher:
  
   \[ (x : * ,y:*,t:Id_*(x,y), s: Id_*(x,y), k:Id_{Id_*(x,y)}(t,s)) \]
   
is also a valid context, at this point (and especially if one want to go even higher) it might be preferable to stop writing the indices of the identity types: $(x : * ,y:*,t:Id(x,y), s: Id(x,y), k:Id(t,s))$.

\bigskip

In $\T^*_{refl}$ and $\T^*_{w.i.e.}$ one can form additional terms, which in turn (because of rules $(T2)$) allow to form new context by using $Id(t_1,t_2)$ for $t_1$ and $t_2$ terms involving $r$,$\J$ or $\Hb$.

}

\blockn{The convention we chose about equalities between $\J$ terms and $\Hb$ terms produces an essential difference between the rules $(T4)$ and $(T5)$ on one hand and the rules $(T2)$ and $(T3)$ on the other hand. Indeed, while (by convention) each application of the rules $(T4)$ and $(T5)$ produces syntactically different terms, different applications of the rules $(T3)$ and $(T2)$ can produces the same terms. For example if one has: $\Gamma \vdash \texttt{ Type A} , a:A, b:A$ then by application of $(T2)$ one has $\Gamma \vdash Id_A(a,b) \texttt{ Type}$, but one also has: $\Gamma,x:A,y:A \vdash Id_A(x,y)$ and one can substitute $x$ and $y$ by $a$ and $b$ using the rules $(S3)$ to obtain again $\Gamma \vdash Id_A(a,b)$. In this example we obtained twice the same judgement by applying different rules, which shows that the syntactic category is not ``Freely generated by the rules'' which will make it considerably more complicated to use and suggest that the corresponding cylinder category that we will obtain might not be cofibrant. There is several solutions to solve this problem, the first would be to take the same conventions as we did for the rules $(T4)$ and $(T5)$, this solution would work but the resulting cylinder category would have considerably more object and its completion would not be formed of ``globular groupoids'' anymore. Instead we will modify the rules $(T2)$ and $(T3)$ in order to restrict them to a ``minimal'' set of applications.}

\block{Consider the following weakening of the rules $(T2)$ and $(T3)$:
\Prop{In the three theory above, the rules $(T2)$ can be replaced by the following infinite family of axioms:
\begin{itemize}
\item[$(T2'_1)$]  $s_0:* ,t_0:* \vdash Id_1(s_0,t_0) \texttt{ Type}$
\[\dots\]
\item[$(T2'_i)$]  $s_0:* ,t_0:*,s_1:Id_1(s_0,t_0),T_1:Id_1(s_0,t_0), \dots, s_{i-1}: Id_{i-1}(s_{i-2}, t_{i-2}), t_{i-1}: Id_{i-1}(s_{i-2},t_{i-2}) \vdash Id_i(s_{i-1},t_{i-1}) \texttt{ Type}$
\end{itemize}

and the rule $(T3)$ can be replaced by the following infinite family of axioms:
\begin{itemize}
\item[$(T3'_i)$] $s_0:*,t_0:*, \dots, s_{i-2}: Id_{i-2}(s_{i-3}, t_{i-3}), t_{i-2}: Id_{i-2}(s_{i-3},t_{i-3}), a:Id_{i-1}(s_{i-2},t_{i-2}) \vdash r(a) : Id_i(a,a) $
\end{itemize}}

Those new rules do not suffer from the problem with the rules $(T2)$ and $(T3)$ mentioned above: each rules produces a syntactically different object, and you can read from the syntax which rules has been used to produce it.

\Dem{All these rules are special case of $(T2)$ and $(T3)$ so we need to show that the general form of the rules $(T2)$ and $(T3)$ can be deduced from those specific cases and manipulation involving the structural rules $(S1),(S2)$ and $(S3)$.

In the three theories above, there is no other rules allowing to form type than $(T1)$ and $(T2)$, hence any type is either of the form $*$ or $Id_A$ for some other type $A$. In particular one can define by induction the ``height'' of a type: $*$ has height $0$ and $Id_A$ has height the height of $A$ plus one.
One define the height of a variable to be the height of its type. A variable $x$ of height $n>0$ has a "source" and a "target" of height $n-1$: indeed $x$ is of type $Id_A(a,b)$ for some $a$ and $b$ of type $A$ of height $n-1$ we call $a$ the source of $x$ and $b$ the target of $x$.

Specifying a variables of heights $i$ definable in some context $\Gamma$ is the same as specifying a context morphism:

\[ \Gamma \rightarrow (s_0:*,t_0:*,s_1:Id_1(s_0,t_0), \dots, v:Id_i(s_{i-1},t_{i-1}) \]

Where (at this point) $Id_i$ is a short cut for: $Id_1 = Id_*$, $Id_2=Id_{Id(s_0,t_0)}$ $\dots$

Indeed such a morphism $v$ will corresponds to some variable of height $i$ and then all the $s_k$ and $t_k$ will be the iterated source and target of $v$. And similarly a pair of variables of a same type $A$ of height $(i-1)$ is the same as a morphism:

\[ \Gamma \rightarrow (s_0:*,t_0:*,s_1:Id_1(s_0,t_0), \dots,s_{i-1} :Id_{i-1}(s_{i-2},t_{i-2}),t_{i-1} :Id_{i-1}(s_{i-2},t_{i-2}) \]

Hence if $\Gamma \vdash A \texttt{ Type}, \Gamma \vdash a :A$ and $\Gamma \vdash b:A$ and if $i$ denotes the height of $A$ there such a morphism from $\Gamma$ to the context above and $Id_A(a,b)$ can be defined as $Id_{i}(a,b)$. Hence any application of the rules $(T2)$ can be deduce as an application of a rule $(T2'_i)$ and an application of the substitution rule.

Similarly, if $\Gamma \vdash A \texttt{ type}$ and $\Gamma \vdash a :A$ with $A$ of height $i$, one has a unique morphism

\[ \Gamma \rightarrow (s_0:*,t_0:*,s_1:Id_1(s_0,t_0), \dots, v:Id_i(s_{i-1},t_{i-1}) \]
which send $a$ and $v$ and one can define $r(a)$ as $r(v)$ after substitution along this morphism and this concludes the proof.
}
}

\blockn{We will now explain how to construct morphisms from a syntactical category as $\Ccal(\T_{w.i.e.}^*)$ to another fibration category (or equivalently from $\Ccal(\T_{w.i.e.}^*)^{op}$ to a cofibration category), and how to detect if such morphisms are morphisms of (pre-)cylinder category/(pre-)path category. }

\block{\label{MorphismOutOfsyntax}If $F : \Ccal(\T) \rightarrow \Fcal$ is a morphism of fibrations category whose domain is a syntactical category of a type theory $\T$, then each judgement of $\T$ can be ``interpreted'' in $\Fcal$:

\begin{itemize}

\item Each valid context $\Gamma$ in $\T$ produce an object $F(\Gamma)$ in $\Fcal$.

\item Each type in context: $\Gamma \vdash A \texttt{ Type}$ produce a fibration $F(\Gamma,A) \twoheadrightarrow F(\Gamma)$ which is the image of the fibration $(\Gamma,x:A) \rightarrow \Gamma$.

\item Each terms $\Gamma \vdash a :A$ produce a sections of the fibration $F(\Gamma,A) \twoheadrightarrow F(\Gamma)$ which corresponds to the map $\Gamma \rightarrow \Gamma,x:A$ which is the identity on all variables except the last ones, which substitute $x$ for $a$.

\end{itemize}

In the special case where $\T$ is generated by ``nice'' rules, it is enough to interpret the rules. In order to construct inductively a morphism from the syntactical category of one of our type theory one just need to explain how each rules $(T1),(T2'_i),(T3'_i),(T4),(T5)$ is interpreted:

\begin{itemize}

\item[(T1)] One has a chosen object $X$ in $\Fcal$.
\item[$(T2')_1$] One has a chosen object $Id^X_1$ in $\Fcal$ endowed with a fibration $Id^X_1 \twoheadrightarrow X \times X$
\item[$(T2')_i$] One has a chosen object $Id^X_i$ endowed with a fibration $p_i: Id^X_i \twoheadrightarrow Id^X_{i-1} \times_{Id^X_{i-2} \times Id^X_{i-2}} Id^X_{i-1}$ where the maps $Id^X_{i-1}  \rightarrow Id^X_{i-2} \times Id^X_{i-2}$ involved in the fiber product are both $p_{i-1}$. (by convention $Id^X_0 =X$ and $Id^X_{-1}$ is the terminal object.)
 
 \item[$(T3')_i$] One has a chosen map $r_{i-1}: Id^{X}_{i-1} \rightarrow Id^X_{i}$ which factor the diagonal maps:
 
 \[  Id^{X}_{i-1}\rightarrow Id^X_{i} \twoheadrightarrow Id^X_{i-1} \times_{Id^X_{i-2} \times Id^X_{i-2}} Id^X_{i-1} \]
 
 (one has in particular for $i=0$, $r_0 : X \rightarrow Id^X_1$)
 
\item[(T4)] Assume one has judgements as in the premise of $(T4)$ which have already been interpreted. It means that one has a diagram of the form:

\[\begin{tikzcd}[ampersand replacement=\&]
 \&  F(\Gamma,A,A,Id_A,\Delta,C) \arrow[two heads]{d}{p} \\
F(\Gamma,A,\Delta[a/b,r(a)/\gamma]) \arrow[two heads]{d} \arrow{ur}{d} \arrow{r}{r'} \& F(\Gamma,A,A,Id_A,\Delta) \arrow[two heads]{d} \\
F(\Gamma,A) \arrow{r}{r} \& F(\Gamma,A,A,Id_A)
\end{tikzcd} \]

where the square is a pullback. Then one need to construct a section $\J$ of the fibration $p$.

\item[(T5)] In the situation above, one has a diagram:

\[ d,\J \circ r' : F(\Gamma,A,\Delta[a/b,r(a)/\gamma]) \rightrightarrows F(\Gamma,A,A,Id_A,\Delta,C) \overset{p}{\twoheadrightarrow} F(\Gamma,A,A,Id_A,\Delta) \] 

We also need to a construct a ``relative homotopy'' between them:
\[ \Hb :F(\Gamma,A,\Delta[a/b,r(a)/\gamma])  \rightarrow F(\Gamma,A,A,Id_A,\Delta,C,C,Id_C)) \]

\end{itemize}

Finally:

\Lem{A morphism of fibration categories  $F:\Ccal(\T_{w.i.e.}^*) \rightarrow \Ccal$ for $\Ccal$ a pre-path category\footnote{i.e. the opposite of a pre-cylinder category.} is a morphism of pre-path category if and only if all the maps in $\Ccal$:

\[r: Id^{X}_{i-1} \rightarrow Id^X_{i}  \]

Produces by rules $(T3'_i)$ above are weak equivalences.

}

\Dem{All those maps $r$ are weak equivalences in $\Ccal(\T_{w.i.e.})$ hence they should be sent to weak equivalences by any morphism of pre-cylinder categories.

Conversely, if all these maps are sent to weak equivalences then any maps of the form $r_A: (\Gamma,A) \rightarrow (\Gamma,A,A,Id_A)$ produced by rules $(T3)$ will also be a weak equivalence (because it is a pullback of one of the maps produced by $(T3'_i)$ in a situation where the cube lemma applies). Finally, An arrow $f:X \rightarrow Y$ is a weak equivalence in $\Ccal(\T_{w.i.e.}^*)$ if and only if it has a ``homotopy inverse'', i.e. a map $g:Y \rightarrow X$ such that $f \circ g$ and $g \circ f$ are ``homotopic to the identity'' where homotopic is defined using the path object constructed in \cite[section 5 and 6]{berg2016path} from the identity types. Hence it is enough to forces the identity type to be sent to relative path object to have that $F$ preserve the equivalences. }
}

\block{\Th{$\Ccal(\T_{w.i.e.}^*)^{op}$ is a cylinder coherator (with marked object $(x:*)$).}

\Dem{The fact that it is a cylinder category is the main result of \cite{berg2016path}. So we just have to prove that the map $F_* \rightarrow \Ccal(\T_{w.i.e.}^*)^{op}$ corresponding to the marked object is a trivial cofibration, i.e. that it has the left lifting property with respect to all fibrations, but this follow immediately from the description of morphisms of pre-cylinder categories $\Ccal(\T_{w.i.e.}^*)^{op} \rightarrow \Ccal$ given in \ref{MorphismOutOfsyntax}:

The rule $(T1)$ corresponds to the choice of an object, which corresponds to the arrow from $F_*$ hence does not need to be lifted. The rules $(T2'_i)$ and $(T3'_i)$ together corresponds to choosing a cofibration/weak equivalences factorization of a morphism (in fact a relative cylinder object) and such a choice can always be lifted along a (pre-)fibrations (definition \ref{Def_FibPCylcat}) and finally axiom $(T4)$ and $(T5)$ together corresponds to the choice of a weak filing of the square (in the fibration category $\Ccal(\T_{w.i.e.}^*)$):

\[\begin{tikzcd}[ampersand replacement=\&]
F(\Gamma,A,\Delta[a/b,r(a)/\gamma]) \arrow{d}{r'}[swap]{\sim} \arrow{r}{d}s \& F(\Gamma,A,A,Id_A,\Delta,C) \arrow[two heads]{d}\\
F(\Gamma,A,A,Id_A,\Delta) \arrow{r}{Id} \& F(\Gamma,A,A,Id_A,\Delta) \\
\end{tikzcd}
\]

exactly as in (the dual of) proposition \ref{Prop_LiftingOfWeakDiag} and hence it can always be lifted through a fibration. 

}
}

\block{\Prop{\begin{itemize}

\item The completion of $\Ccal(\T^*_{Id})^{op}$ is the category of globular sets.
\item The completion of $\Ccal(\T^*_{refl})^{op}$ is the category of reflexive globular sets (i.e. globuar set with additional maps $r:G_n \rightarrow G_{n+1}$ such that $s\circ r =t \circ r = Id$).
\item The completion of $\Ccal(\T^*_{w.i.e.})^{op}$ can be described as the category of globular sets $(G_{n})$ endowed with additional operation which are maps between finite projective limits of the $G_n$ cells satisfying certain relations, the precise definition of those maps being controlled by the type theory $\T^*_{w.i.e.}$.

\end{itemize}
}

\Dem{In all case it is a direct application of the description of morphisms out of these syntactical category given in \ref{MorphismOutOfsyntax} to understand the completion as the category of morphisms to the pre-cylinder category $\setop$. In the first case, one just need to pick a set $G_0$ and then a series of set $G_n$ endowed with a map $G_n \rightarrow G_{n-1} \times_{G_{n-2} \times G_{n-2}} G_{n-1}$ and this corresponds exactly to a globular set, in the second case one need to add the reflexivity map $r:G_n \rightarrow G_{n+1}$ and this corresponds exactly to a reflexive globular set, and in the third case one need to add all the operations that can be produces using the reflexivity terms and the rules $(T4)$ and $(T5)$.
}
}

\block{Finally we would like to point out that for the cylinder coherator $\Ccal = \Ccal(\T^*_{w.i.e.})$ the notion of $\Ccal$-groupoid should be definable within intentional type theory using the same strategy of syntactic groupoid as G.Brunerie in \cite{brunerie2013syntactic}, we just have to replace the type theory introduced in \cite{brunerie2013syntactic} by the more classical type theory $\T^*_{w.i.e.}$ introduced here.
}

\bibliography{Biblio}{}
\bibliographystyle{plain}

\end{document}